\documentclass[11pt]{article}

\usepackage{amsmath,amssymb,amsfonts,amsthm,mathrsfs}
\usepackage{bbm}
\usepackage{graphicx}
\usepackage{xcolor}
\usepackage{ifpdf}
\usepackage{todonotes}
\allowdisplaybreaks[4]

\newtheorem{theorem}{Theorem}[section]
\newtheorem{lemma}{Lemma}[section]
\newtheorem{corollary}{Corollary}[section]
\newtheorem{proposition}{Proposition}[section]
\newtheorem{remark}{Remark}[section]

\newtheorem{definition}{Definition}[section]
\numberwithin{equation}{section}

\title{Emergent behaviors of relativistic thermodynamic flocks with Synge energy}

\author{
Ziming Bian\thanks{Innovation Academy for Precision Measurement Science and Technology, Chinese Academy of Sciences, Wuhan 430071, China, \texttt{bianziming23@mails.ucas.ac.cn}} 
\and
Seung-Yeal Ha\thanks{Department of Mathematical Sciences and Research Institute of Mathematics, Seoul National University, Seoul 08826, Republic of Korea, \texttt{syha@snu.ac.kr}}
\and
Tommaso Ruggeri\thanks{Department of Mathematics, University of Bologna, Bologna, Italy and Accademia Nazionale dei Lincei, Roma, Italy, \texttt{tommaso.ruggeri@unibo.it}}
\and
Qinghua Xiao\thanks{Innovation Academy for Precision Measurement Science and Technology, Chinese Academy of Sciences, Wuhan 430071, China, corresponding author, \texttt{xiaoqh@apm.ac.cn}} 
}

\date{} 

\begin{document}

\maketitle

\begin{abstract}
Collective motion and self-organization of interacting particles, such as flocking and swarming, can be viewed as nonequilibrium analogues of collective dynamics in gases. Motivated by the analogy between gas mixtures and Cucker--Smale models, we introduce a polyatomic classical model and its relativistic counterpart based on the Synge energy, and analyze their large-time behavior. The relativistic formulation provides a physically consistent setting for multi-species systems where inertia and internal energy depend on temperature, as occurs in astrophysical plasmas or relativistic fluids. Using the entropy principle, we derive uniform lower bounds for temperature and establish asymptotic flocking under various communication kernels. For nearly constant interactions, flocking emerges from arbitrary initial data. The results clarify how thermodynamic effects and relativistic corrections modify the emergence of coherent motion in particle systems, bridging kinetic theory, relativistic fluid mixtures, and collective dynamics.
\end{abstract}


%



\maketitle		
	
	
\section{Introduction}\label{sec:1}
\setcounter{equation}{0}

Emergent behaviors of self-propelled particles are ubiquitous in nature: examples include aggregation of bacteria \cite{T-B, T-B-L}, flocking of birds \cite{C-S,R,T-T,V-C-B-C-S}, synchronization of fireflies and pacemaker cells \cite{B-B, Ku2,Pe, Wi2}, and swarming of fish \cite{Ao, D-M2}. For general overviews, we refer to \cite{A-B-F, B-H, D-B1, P-R,VZ, Wi1}. Motivated by Reynolds’s pioneering work \cite{R}, Vicsek and co-workers \cite{V-C-B-C-S} introduced a stochastic discrete model in which each particle moves at unit speed and aligns its direction with the average of its neighbors. Their numerical results revealed distinct macroscopic patterns depending on local density and noise. A rigorous analysis was first provided in \cite{J-L-M}, and later Cucker and Smale \cite{C-S} proposed a deterministic mean-field model that allows for analytical characterization of flocking and its dependence on system parameters and initial configurations.

Ha and Ruggeri later observed that the mechanical Cucker–Smale (CS) model can be identified with the momentum equations for a mixture of gases in the isothermal case. Extending this analogy, they incorporated the energy balance of each species, obtaining the thermodynamic CS (TCS) model \cite{Ha-Ruggeri-ARMA-2017}, where temperature effects couple to mechanical alignment. This viewpoint opened a path connecting flocking models to kinetic theory and nonequilibrium thermodynamics. Building on this connection, Ha, Kim, and Ruggeri \cite{Ha-Kim-Ruggeri-ARMA-2020} proposed a relativistic version of the TCS model, where the internal energy was assumed proportional to temperature as in the ideal gas. However, at the kinetic level, relativistic energy has a more complex temperature dependence given by the Synge equation, involving ratios of modified Bessel functions.

From a physical perspective, the study of relativistic mixtures is relevant in various settings—such as astrophysical plasmas, relativistic jets, or the early universe—where several species interact at relativistic velocities and exchange energy and momentum. In such systems, relativistic corrections to inertia and temperature significantly affect collective motion and equilibration. The relativistic extension of the TCS model therefore provides a consistent framework for describing how thermodynamic and relativistic effects shape the emergence of organized behavior in multi-species particle systems.

In this work, we introduce a polyatomic classical model together with its relativistic counterpart governed by the Synge energy. We establish asymptotic flocking results under various communication kernels and derive uniform lower bounds for the temperature via the entropy principle. For nearly constant communication, flocking emerges from arbitrary initial configurations. We further analyze how relativistic effects modify the collective behavior for both monatomic and polyatomic gases, offering new insight into the interplay between kinetic temperature, inertia, and alignment in relativistic regimes.

\section{Preliminaries}
\setcounter{equation}{0}
In this section, we study some preliminaries on the TCS model and review the historical development on the thermodynamic counterparts of the CS model and present the new, more refined model, highlighting its classical limit in the case of polyatomic gas mixtures. \newline

Let $\mathbf{x}_a$, $\rho_a$, $\mathbf{v}_a$, and $\varepsilon_a$ be the position, mass density, velocity, and internal energy of the $a$-th particle, respectively, for $a\in [N]:=\{1,\dots,N\}$. In \cite{Ha-Ruggeri-ARMA-2017},  Ha and Ruggeri proposed a generalization of the Cucker--Smale model within the framework of thermo-mechanics, based on homogeneous solutions for mixtures of gases with multiple temperatures \cite{R-S}:
\begin{equation}\label{TCS-complete}
	\begin{cases}
		\displaystyle \frac{\mathrm{d}\mathbf{x}_a}{\mathrm{d} t}=\mathbf{v}_a,\quad  a \in [N], \\
		\displaystyle \frac{\mathrm{d}\rho_a}{\mathrm{d} t}=0, \\
		\displaystyle \frac{\mathrm{d} \rho_a \mathbf{v}_a}{\mathrm{d} t}
		=\frac{1}{N}\sum_{b=1}^{N}\phi_{ab}\left(\frac{ \mathbf{v}_{b}}{T_b}-\frac{\mathbf{v}_{a}}{T_a}\right),\\
		\displaystyle \frac{\mathrm{d}}{\mathrm{d} t}\left(\rho_a \varepsilon_a +\frac{\rho_a}{2}{|\mathbf{v}_a|^2}\right)
		=\frac{1}{N}\sum_{b=1}^{N}\zeta_{ab}\left(\frac{1 }{T_a}-\frac{1 }{T_b}\right),
	\end{cases}
\end{equation}
where we assume that the center-of-mass velocity vanishes without loss of generality thanks to the Galilean invariance:
\[
\langle \mathbf{v} \rangle =\frac{1}{N}\sum_{a=1}^{N}\mathbf{v}_{a} = 0,
\]
and $|\cdot|$ denotes the Euclidean $\ell^2$-norm, and $\phi_{ab}$ and $\zeta_{ab}$ are communication weight functions depending on the spatial distance between particles:
\begin{align*}
	&\phi_{ab}=\phi(|\mathbf{x}_{a}-\mathbf{x}_{b}|),\quad
	\zeta_{ab}=\zeta(|\mathbf{x}_{a}-\mathbf{x}_{b}|), \quad a, b \in [N].
\end{align*}
Since the mass densities are constants, we normalize them as the unity. For a rarefied polyatomic gas, the internal energy of a single species is proportional to its temperature $T_a$:
\begin{equation}\label{cvac}
	\varepsilon_a = c_{V_a} T_a, \quad a \in [N].
\end{equation}
Here $c_{V_a}$ is the specific heat at constant volume, which is related to the number of degrees of freedom $D_a$ for species $a$, given by
\begin{equation}\label{calore}
	c_{V_a} = \frac{D_a}{2}\frac{k_B}{m_a},
\end{equation}
$k_B$ is the Boltzmann constant and $m_a$ is the atomic mass of species $a$.
In the previous works, the authors assumed that all species have the same number of degrees of freedom and the same molecular mass, and  we set $c_{V_a} = 1$. Therefore, the Thermodynamic Cucker–Smale (TCS) model was reduced to the study of the following ODE system:
\begin{equation*}
	\begin{cases}
		\displaystyle \frac{\mathrm{d}\mathbf{x}_a}{\mathrm{d} t}=\mathbf{v}_a,\quad  a \in [N], \\
		\displaystyle \frac{\mathrm{d} \mathbf{v}_a}{\mathrm{d} t}
		=\frac{1}{N}\sum_{b=1}^{N}\phi_{ab}\Big(\frac{ \mathbf{v}_{b}}{T_b}-\frac{\mathbf{v}_{a}}{T_a}\Big),\\
		\displaystyle \frac{\mathrm{d}}{\mathrm{d} t}\left( T_a +\frac{1}{2}{|\mathbf{v}_a|^2}\right)
		=\frac{1}{N}\sum_{b=1}^{N}\zeta_{ab}\Big(\frac{1 }{T_a}-\frac{1 }{T_b}\Big).
	\end{cases}
\end{equation*}
Note that the first two equations with constant temperature $T_a\equiv 1,~~a \in [N]$ reduce to the classical CS model. The TCS and CS models have been extensively investigated in literature from diverse aspects, e.g., emergence of deterministic and stochastic asymptotic flocking \cite{A-H, C-F-R-T, C-S, C-H-H-J-K, Ha-Kim-Ruggeri-SIMA-2018, H-L, Ha-Ruggeri-ARMA-2017, H-T}, kinetic and hydrodynamic descriptions \cite{H-T}, uniform stability and mean-field limits \cite{A-H-K-S-S, H-K-Z, H-L}, etc. We also refer to survey articles \cite{C-F-T-V,C-H-L}.

If we want to characterize a particular polyatomic gas, it is convenient---also for what concerns  the relativistic limit---to take into account the internal degrees of freedom of the particles. Hence, even if we set $k_B/m_a = 1$, we cannot, in general, set $c_{V_a} = 1$. Assuming that each species has the same number of degrees of freedom, $D_a = D$, then by \eqref{calore} we may write
\[
c_{V_a} = \frac{2\chi + 1}{2}, \quad \chi = \frac{D - 1}{2}.
\]
When $\chi$ takes values in $\{1,2,3,4\}$, we recover a mixture of monatomic, diatomic, triatomic, and tetratomic gases, respectively.
Taking this aspect into account, we write the TCS model as
\begin{equation}\label{A-3}
	\begin{cases}
		\displaystyle \frac{\mathrm{d}\mathbf{x}_a}{\mathrm{d} t}=\mathbf{v}_a, \quad a \in [N], \\
		\displaystyle \frac{\mathrm{d} \mathbf{v}_a}{\mathrm{d} t}
		=\frac{1}{N}\sum_{b=1}^{N}\phi_{ab}\Big(\frac{ \mathbf{v}_{b}}{T_b}-\frac{\mathbf{v}_{a}}{T_a}\Big),\\
		\displaystyle \frac{\mathrm{d}}{\mathrm{d} t}\left(\frac{2\chi +1}{2}T_a+\frac{1}{2}{|\mathbf{v}_a|^2}\right)
		=\frac{1}{N}\sum_{b=1}^{N}\zeta_{ab}\Big(\frac{1 }{T_a}-\frac{1 }{T_b}\Big).
	\end{cases}
\end{equation}
Using the same analogy between mixtures and Cucker–Smale models, Ha, Kim, and Ruggeri \cite{Ha-Kim-Ruggeri-ARMA-2020} were able to construct a relativistic counterpart of the classical TCS \eqref{TCS-complete}. The corresponding RTCS model then reads as follows:
\begin{equation} \label{A-2}
	\begin{cases}
		\displaystyle \frac{\mathrm{d}\mathbf{x}_a}{\mathrm{d} t}=\mathbf{v}_a, \quad  a \in [N], \\[1ex]
		\displaystyle \frac{\mathrm{d}(\rho_a \Gamma_a)}{\mathrm{d} t} =0, \\[1ex]
		\displaystyle \frac{\mathrm{d}}{\mathrm{d} t} \left [ \rho_a \Gamma_a^2\mathbf{v}_a\left(1+\frac{p_a+\rho_a\varepsilon_a}{\rho_a c^2}\right)\right ]
		=\frac{1}{N}\sum_{b=1}^{N}\phi_{ab}\Big(\frac{\Gamma_b \mathbf{v}_{b}}{T_b}-\frac{\Gamma_a \mathbf{v}_{a}}{T_a}\Big), \\[1ex]
		\displaystyle \frac{\mathrm{d}}{\mathrm{d} t} \Big[ p_a(\Gamma_a^2-1)+\rho_a\big(\Gamma_a^2\varepsilon_a+c^2\Gamma_a(\Gamma_a-1)\big) \Big ] 
		=\frac{1}{N}\sum_{b=1}^{N}\zeta_{ab}\Big(\frac{\Gamma_a}{T_a}-\frac{\Gamma_b}{T_b}\Big),
	\end{cases}
\end{equation}
where $\rho_a = n_a m_a$ is the density, $n_a$ is the particle number, and $m_a$ is the rest mass. Here $p_a$ is the pressure, $\Gamma_a$  and $e_a$ are the Lorentz factor for the $a$-th particle and energy, respectively defined by the following explicit forms:
\begin{equation} \label{suminterres}
	\Gamma_a := \frac{1}{\sqrt{1-\frac{v^2_a}{c^2}}}, \qquad  e_a := \rho_a(c^2+\varepsilon_a), \quad  \forall~a\in[N].
\end{equation}
Note that energy is the sum of the rest-frame energy and the internal energy $\varepsilon_a$.  To close the system \eqref{A-2}, we also need additional constitutive equations. System \eqref{A-2} is very complex, in particular for the presence of a pressure. For simplicity, it was assumed that $p_a$ is negligible  in \cite{Ha-Kim-Ruggeri-ARMA-2020}. Moreover, as in the case of a single fluid, the expression of $\varepsilon_a$ was assumed to be the same as in the classical case \eqref{cvac} with $c_{V_a}=1$. From the second equation in \eqref{A-2}, we also set the constant $\rho_a \Gamma_a=1$:
\begin{equation}
	\begin{cases} \label{A-2s}
		\displaystyle \frac{\mathrm{d}\mathbf{x}_a}{\mathrm{d} t}=\mathbf{v}_a,\quad a\in[N],\\
		\displaystyle \frac{\mathrm{d}}{\mathrm{d} t} \left [ \Gamma_a \mathbf{v}_a \Big(1 + \frac{T_a}{c^2} \Big) \right ]
		=\frac{1}{N}\sum_{b=1}^{N}\phi_{ab}\Big(\frac{\Gamma_b \mathbf{v}_{b}}{T_b}-\frac{\Gamma_a \mathbf{v}_{a}}{T_a}\Big),\\
		\displaystyle \frac{\mathrm{d}}{\mathrm{d} t} \left( \Gamma_a  T_a  + c^2   (\Gamma_a -1)  \right)  =\frac{1}{N}\sum_{b=1}^{N}\zeta_{ab}\Big(\frac{\Gamma_a }{T_a}-\frac{\Gamma_b }{T_b}\Big).    
	\end{cases}
\end{equation}
Even with these simplifications, the RTCS model remains to be difficult to handle. Therefore, the relativistic CS mechanical case was considered by using the theory of principal subsystems due to Boillat and Ruggeri \cite{Boillat-Ruggeri}.  
This simplified case is obtained by neglecting the last energy equation in \eqref{A-2} and taking $T_a = T^* \Gamma_a$ with $T^* = \text{const.} = 1$ in the remaining equations,  see \cite{Ha-Kim-Ruggeri-ARMA-2020} for more details.  
This approximate relativistic CS model \eqref{A-2} reads as follows.
\begin{equation}
	\begin{cases} \label{A-2ss}
		\displaystyle \frac{\mathrm{d}\mathbf{x}_a}{\mathrm{d} t}=\mathbf{v}_a,\quad a\in[N],\\
		\displaystyle \frac{\mathrm{d}}{\mathrm{d} t} \left [ \Gamma_a \mathbf{v}_a \left(1 + \frac{\Gamma_a}{c^2} \right) \right ]
		=\frac{1}{N}\sum_{b=1}^{N}\phi_{ab} ( \mathbf{v}_{b}- \mathbf{v}_{a}).
	\end{cases}
\end{equation}
System  \eqref{A-2s} and in particular the simplified one \eqref{A-2ss} have been extensively investigated from different perspectives, e.g., derivation and asymptotic flocking of the relativistic CS model \cite{Ha-Kim-Ruggeri-ARMA-2020}, non-relativistic limit \cite{A-H-K1,H-R-X}, uniform stability and mean-field limit \cite{A-H-K2}, kinetic and hydrodynamic counterparts \cite{H-K-R2}, etc.

So far, the model has been constructed by analogy with phenomenological mixture theory.  
However, for rarefied gases, an alternative formulation can be obtained from kinetic theory, which provides a counterpart of mixture theory.  
In the classical setting, this leads to models with distinct production terms appearing on the right-hand side of the system.  
A detailed comparison between the phenomenological and kinetic approaches was presented in a recent paper \cite{PhysicaD}.  

In the relativistic framework, an even subtler issue arises, which has already been seen in the case of a single fluid.  
A crucial aspect of relativistic fluid dynamics is the closure relation, namely the constitutive equation that links pressure, energy density, and particle number density.  
In contrast to classical fluid mechanics, where the equation of state is well established (for instance, the ideal gas law), in the relativistic context the formulation of a general and accurate equation of state is a highly nontrivial task.  
Most relations currently in use are valid only in asymptotic regimes: either in the classical (non-relativistic) limit or in the ultrarelativistic regime, where the rest mass of the particles becomes negligible or the temperature extremely high.  

As emphasized by Ruggeri, Xiao and Zhao \cite{Ruggeri-Xiao-Zhao-ARMA-2021}, this limitation represents one of the weaknesses of the present theory of the relativistic fluids.  
To address this challenge, insights from the mesoscopic scale—particularly from relativistic kinetic theory—prove to be essential.  
For rarefied gases, Synge \cite{Synge} derived a constitutive equation based on the equilibrium Jüttner distribution, yielding a complex yet physically consistent expression involving ratios of modified Bessel functions of the second kind.  
Synge’s equation has the important property of reduction to the classical ideal gas law in the low-temperature limit and of approaching simplified forms in the ultrarelativistic regime.  
This framework was later extended to polyatomic gases by Pennisi and Ruggeri \cite{Pennisi_Ruggeri}, who obtained Synge-type energies with explicit formulas in the diatomic case.  
Although they are physically sound, these models are mathematically intricate and only  the Riemann problem has been rigorously solved \cite{Ruggeri-Xiao-Zhao-ARMA-2021}.

In this case, it follows from the general system \eqref{A-2} that we have (see Appendix \ref{App-C}):
\begin{equation}
	\begin{cases} \label{A-4}
		\displaystyle \frac{\mathrm{d}\mathbf{x}_a}{\mathrm{d} t}=\mathbf{v}_a,\quad t > 0, \quad a\in[N],\\
		\displaystyle \frac{\mathrm{d}}{\mathrm{d} t}\left(\Gamma_a\mathbf{v}_aH_a\right)
		=\frac{1}{N}\sum_{b=1}^{N}\phi_{ab}\Big(\frac{\Gamma_b \mathbf{v}_{b}}{T_b}-\frac{\Gamma_a \mathbf{v}_{a}}{T_a}\Big),\\
		\displaystyle \frac{\mathrm{d}}{\mathrm{d} t}\left[c^2\left(\Gamma_aH_a
		-1-\frac{1}{\gamma_a \Gamma_a}\right)\right]=\frac{1}{N}\sum_{b=1}^{N}\zeta_{ab}\Big(\frac{\Gamma_a }{T_a}-\frac{\Gamma_b }{T_b}\Big),   
	\end{cases}
\end{equation}
with
\begin{equation} \label{B-1-1}
	H_a:=
	\begin{cases}
		\displaystyle  \frac{K_1(\gamma_a)}{K_2(\gamma_a)}+\frac{4}{\gamma_a},\quad \qquad \qquad \qquad \qquad \mbox{monatomic~ gases}, \\
		\displaystyle \frac{K_0(\gamma_a)}{K_1(\gamma_a)}+\frac{4}{\gamma_a},\quad \qquad \qquad \qquad \qquad \mbox{diatomic~ gases}, \\
		\displaystyle \frac{K_1(\gamma_a)}{\gamma_a\int_{\gamma_a}^{\infty}\frac{K_1(y)}{y}\mathrm{d}y}+\frac{3}{\gamma_a}, \qquad \qquad \qquad \mbox{triatomic~ gases}, \\
		\displaystyle \frac{K_0(\gamma_a)}{\gamma_aK_0(\gamma_a)-\gamma_a^2 \int_{\gamma_a}^{\infty}\frac{K_1(y)}{y}\mathrm{d}y}+\frac{3}{\gamma_a},\quad \mbox{tetratomic~ gases},
	\end{cases}
\end{equation}
where $K_j(\cdot)$ denotes the modified Bessel functions of the second kind, defined in (\ref{defini}).  
For the simplicity of presentation, we normalize constants $m_a$ and $k_B$ to unity, so that  
\[
\gamma_a := \frac{m_a c^2}{k_B T_a} = \frac{c^2}{T_a}.
\]  
In the relativistic regime, we have $\gamma_a \ll 1$.  Moreover, in the limits $\gamma_a \to \infty$ and $\gamma_a \to 0$, we recover, respectively, the classical and ultrarelativistic regimes.






\section{The classical polyatomic TCS model} \label{sec:2.1}
\setcounter{equation}{0}
In this section, we provide basic properties of the TCS model \eqref{A-3}.
For notational simplicity, we use the following handy notations:
\begin{align*}
	\begin{aligned}
		& \max_{a} := \max_{a \in [N]}, \quad \max_{a, b} := \max_{a,b \in [N]}, \quad  \min_{a} := \min_{a \in [N]}, \quad \min_{a, b} := \min_{a,b \in [N]}, \\
		& \overline{T}(t):=\max_{a\in [N]}T_a(t),  \quad \underline{T}(t):=\min_{a\in [N]}T_a(t), \quad t \geq 0. \\
	\end{aligned}
\end{align*}
Now we recall the concept of asymptotic themo-mechanical flocking in the following definition.
\begin{definition} \label{D2.1}
	The configuration $\{ (\mathbf{x}_a, \mathbf{v}_a, T_a) \}$ exhibits asymptotic thermo-mechanical flocking if and only if the following relations hold:
	\begin{enumerate}
		\item
		Spatial cohesion holds:
		\[ \sup_{0\leq t<\infty} \max_{a, b} |\mathbf{x}_a(t)-\mathbf{x}_b(t)| < \infty. \]
		\item
		Velocity and temperature  alignment holds asymptotically: 
		\[
		\lim_{t \to \infty} \max_{a, b} |\mathbf{v}_a(t) -\mathbf{v}_b(t)| = 0, \quad \lim_{t \to \infty} \max_{a, b} |T_a(t)-T_b(t)| =0. 
		\]
	\end{enumerate}
\end{definition}

\vspace{0.2cm}

In the sequel, we introduce total momentum, entropy, and energy as follows:
\begin{equation} \label{B-0}
	\mathbf{M}:=\sum_{a=1}^{N}\mathbf{v}_a,\quad S:=\sum_{a=1}^{N}\ln{T_a}, \quad 
	E_\chi:=\sum_{a=1}^{N}\left(\frac{2\chi+1}{2}{T_a} +\frac{1}{2}|{\mathbf{v}_a}|^2\right), \quad \chi \in [4].
\end{equation}
Next, we provide several preparatory estimates to be used in later analysis. 
\begin{lemma}\label{L2.1}
	Let $\left\{\left(\mathbf{x}_a, \mathbf{v}_a, T_a\right)\right\}$ be a global solution to \eqref{A-3} satisfying 
	\begin{equation} \label{B-0-0}
		|\mathbf{M}(0)|<\infty, \quad E_\chi(0)<\infty, \quad \min_{a} T_a(0) > 0.
	\end{equation}
	Then, the following assertions hold: 
	\begin{enumerate}
		\item
		Total momentum and energy are conserved:
		\[
		\mathbf{M}(t)=\mathbf{M}(0), \quad E_\chi(t)=E_\chi(0), \quad t \geq 0.
		\]
		\item
		Entropy is monotonically increasing:
		\[
		\frac{\mathrm{d}S(t)}{\mathrm{d}t}\geq0, \quad t > 0.
		\]
		\item
		Temperatures are bounded below away from zero and bounded above: 
		\[ 0<\frac{\prod_{a=1}^{N}T_a(0)}{\left[2E_\chi(0)/(2\chi+1)\right]^{N-1}} < T_a(t)< \frac{2 E_\chi(0)}{2\chi+1}<\infty. \]
	\end{enumerate}
\end{lemma}
\begin{proof} Since the proofs in the first two assertions have already been verified in \cite{Ha-Ruggeri-ARMA-2017}, we do not repeat them here.  Thus, we focus on the proof of the third assertion. \newline
	
	\noindent (3) We define a set and its supremum:
	\[
	{\mathcal T} := \{ \tau \in (0, \infty):~\min_{a} T_a(t) > 0, \quad t \in [0, \tau) \}, \quad \tau^{\infty}: = \sup {\mathcal T}.
	\]
	It follows from \eqref{B-0-0} and the continuity of $T_a$ that there exists a positive constant $\delta$ such that 
	\[ \min_{a} T_a(t) > 0, \quad t \in [0, \delta). \]
	Hence, we have
	\[ \delta \in {\mathcal T} \quad \mbox{and}  \quad \tau^{\infty} \in [0, \infty]. \]
	We claim that 
	\[ \tau^{\infty} = \infty. \]
	Suppose the contrary holds, i.e., $\tau^{\infty} < \infty$ and we have
	\begin{equation} \label{B-0-1}
		\min_{a} T_a(t)  > 0, \quad t \in [0, \tau^{\infty}) \quad \mbox{and} \quad \lim_{t \to \tau^{\infty}} \min_{a} T_a(t) = 0.
	\end{equation}
	In what follows, we will show that the relations in \eqref{B-0-1} lead to a contradiction and conclude that  
	\[ \tau^{\infty} = \infty \quad \mbox{and} \quad  \min_{a} T_a(t) > 0, \quad \forall~t \in [0, \infty). \]
	For this, we derive the desired upper and lower bounds for temperatures. \newline
	
	\noindent $\bullet$~(Upper bound estimate):~We use the positivity of the temperature and the first assertion to obtain
	\begin{align}\label{up-T}
		\frac{2\chi+1}{2}T_a(t)\leq \sum_{a=1}^N \frac{2\chi+1}{2}T_a(t)\leq \sum_{a=1}^N\frac{2\chi+1}{2}T_a(0)+\frac{1}{2}\sum_{a=1}^N|\mathbf{v}_a(0)|^2= E_\chi(0), \quad \forall~t \in [0, \tau^{\infty}).
	\end{align}

	\noindent $\bullet$~(Lower bound estimate):~We use the positivity of temperature and the second assertion to obtain
	\begin{align*}
		\frac{\mathrm{d}}{\mathrm{d}t}\sum_{a=1}^{N}\ln{T_a}\ge 0, \quad \mbox{i.e.,} \quad  \sum_{a=1}^{N}\ln{T_a(t)}\geq&\sum_{a=1}^{N}\ln{T_a(0)},  \quad \forall~t \in [0, \tau^{\infty}).
	\end{align*}
	This implies
	\begin{align} \label{low-T}
		\begin{aligned}
			\left(\max_{a}T_a(t)\right)^{N-1}\min_{a}T_a(t)\geq&\prod_{a=1}^{N}T_a(0),  \quad \forall~t \in [0, \tau^{\infty}).
		\end{aligned}
	\end{align}
	On the other hand, it follows from \eqref{up-T} that 
	\begin{equation} \label{B-0-2}
		\max_{a} T_a(t) \leq \frac{2 E_\chi(0)}{2\chi+1}.
	\end{equation}
	We combine \eqref{low-T} and \eqref{B-0-2} to find 
	\[
	\left(\frac{2 E_\chi(0)}{2\chi+1}\right)^{N-1}\min_{a}T_a(t) \geq \prod_{a=1}^{N}T_a(0),  \quad \forall~t \in [0, \tau^{\infty}).
	\]
	This implies 
	\[ \min_{a}T_a(t) \geq \frac{\prod_{a=1}^{N}T_a(0)}{\left[2E_\chi(0)/(2\chi+1)\right]^{N-1}}, \quad  \forall~t \in [0, \tau^{\infty}), \]
	which is contradictory to \eqref{B-0-1}. Hence, we have $\tau^{\infty}$ and derive the desired estimates.
\end{proof}
\noindent Throughout the paper, for the the TCS model \eqref{A-3}, we set 
\begin{equation} \label{B-0-1}
	\underline{T}: = \frac{\prod_{a=1}^{N}T_a(0)}{\left[2E_\chi(0)/(2\chi+1)\right]^{N-1}} \quad \mbox{and} \quad  \overline{T} := \frac{2 E_\chi(0)}{2\chi+1}.
\end{equation}
Note that $\underline{T}$ and $\overline{T}$ are independent of $t$.  As the direct application of Lemma \ref{L2.1}, we derive a characterization of asymptotic states for \eqref{A-3}.
\begin{corollary} \label{C2.1}
	Let $\left\{\left(\mathbf{x}_a, \mathbf{v}_a, T_a\right)\right\}$ be a global solution to \eqref{A-3}.  If asymptotic velocity and temperature alignment occurs, i.e., 
	\[ \lim_{t \to \infty} (\mathbf{v}_a, T_a) = (\mathbf{v}^\infty, T^\infty), \quad \forall~a \in [N], \]
	then the asymptotic limits  are explicitly given in terms of initial data:
	\[
	\mathbf{v}^\infty =  \frac{ \mathbf{M}(0)}{N} \quad \mbox{and} \quad  T^\infty  = \Big( \chi + \frac{1}{2}  \Big)^{-1} \Big[  \frac{E_\chi(0)}{N} - \frac{1}{2N^2} |  \mathbf{M}(0)|^2  \Big ].
	\]
\end{corollary}
\begin{proof}
	It follows from the momentum conservation law in Lemma \ref{L2.1} that 
	\[  N \textbf{v}^\infty = \lim_{t \to \infty} \mathbf{M}(t)=\mathbf{M}(0) \quad \mbox{i.e.,} \quad   \textbf{v}^\infty = \frac{1}{N} \mathbf{M}(0). \]
	On the other hand, the conservation of energy implies 
	\begin{equation} \label{B-1}
		\sum_{a=1}^{N}\left(\frac{2\chi +1}{2}{T_a}(t)+\frac{1}{2}|{\mathbf{v}_a}(t)|^2\right) = E_\chi(0).
	\end{equation}
	Letting $t \to \infty$ to \eqref{B-1}, we get the desired analytical expression for $T^{\infty}$:
	\[
	N \Big[ \Big( \chi + \frac{1}{2} \Big) T^\infty +  \frac{1}{2N^2} | \mathbf{M}(0)|^2 \Big] = E_\chi(0). 
	\]
\end{proof}
\section{The relativistic RTCS model with Synge energy} \label{sec:2.2}
In this section, we study basic properties of system \eqref{A-4} with \eqref{B-1-1}.
As in \cite{Ha-Ruggeri-ARMA-2017}, we also introduce auxiliary observables:
\begin{align*}
	F_a:=\Gamma_aH_a, \quad    \mathbf{w}_a := F_a \mathbf{v}_a, \quad \forall~a\in[N].
\end{align*}
From \eqref{K0/K1}, \eqref{H_a(D=7)} and \eqref{H_a(D=9)}, we can rewrite $H_a$ as follows:
\begin{align}\label{H_a(same num)}
	H_a=1+\frac{2\chi +3}{2\gamma_a}+\mathcal{O}(\gamma_a^{-2}).
\end{align}
According to \eqref{energy-m(D=3)}, \eqref{energy-m(D=5)}, \eqref{energy-m(D=7)} and \eqref{energy-m(D=9)}, this implies
\begin{align}\label{energy-m(same num)}
	c^2\left(\Gamma_aH_a
	-1-\frac{1}{\gamma_a \Gamma_a}\right)
	=\frac{2 \chi +1 }{2}T_a+\frac{|\mathbf{w}_a|^2}{2}+\mathcal{O}(c^{-2}).
\end{align}
To connect the RTCS model \eqref{A-4} with the classical TCS model \eqref{A-3}, we introduce an error function  ${\mathcal F}_a$  as follows
\begin{align}\label{def-fa}
	{\mathcal F}_a := {\mathcal F}_a(T_a,|\mathbf{w}_a|^2):=c^2\left[c^2\left(\Gamma_aH_a
	-1-\frac{1}{\gamma_a \Gamma_a}\right)-\frac{2\chi +1}{2}{T_a}-\frac{1}{2}|\mathbf{w}_a|^2\right].
\end{align}
Now, we use  \eqref{H_a(same num)}, \eqref{energy-m(same num)} and \eqref{def-fa}  to rewrite \eqref{A-4} as
\begin{equation} 
	\begin{cases} \label{RTCSmodel(same num)-1}
		\displaystyle \frac{\mathrm{d}\mathbf{x}_a}{\mathrm{d} t}=\frac{\mathbf{w}_a}{F_a}, \quad t > 0, \quad  a \in [N], \\
		\displaystyle \frac{\mathrm{d} \mathbf{w}_a}{\mathrm{d} t}
		=\frac{1}{N}\sum_{b=1}^{N}\phi_{ab}\Big(\frac{ \Gamma_b\mathbf{w}_{b}}{F_bT_b}-\frac{\Gamma_a\mathbf{w}_{a}}{F_aT_a}\Big),\\
		\displaystyle \frac{\mathrm{d}}{\mathrm{d} t}\left(\frac{2\chi +1}{2}{T_a}+\frac{1}{2} |\mathbf{w}_a|^2+\frac{{\mathcal F}_a}{c^2}\right)=\frac{1}{N}\sum_{b=1}^{N}\zeta_{ab}\Big(\frac{\Gamma_a }{T_a}-\frac{\Gamma_b }{T_b}\Big).
	\end{cases}
\end{equation}
Similar to \eqref{B-0}, we set 
\begin{align}
	\begin{aligned}  \label{B-2}
		\mathbb{M} :=\sum_{a=1}^{N}\mathbf{w}_a,\quad \mathbb{E}_\chi:=\sum_{a=1}^N\left[c^2\left(\Gamma_aH_a
		-1-\frac{1}{\gamma_a \Gamma_a}\right)\right], \quad \chi \in [4]. \\
	\end{aligned}
\end{align}
To define an entropy, we recall Gibbs' equation:
\begin{align*}
	T_a \mathrm{d} {\mathbb S}_a=\mathrm{d}\varepsilon_a-\frac{p_a}{\rho_a^2}\mathrm{d}\rho_a,
\end{align*}
and observables
\begin{align*}
	p_a=\frac{\rho_a c^2}{\gamma_a},\quad
	e_a=\rho_a c^2\left(H_a-\frac{1}{\gamma_a}\right),\quad
	\varepsilon_a=c^2\left(H_a-\frac{1}{\gamma_a}-1\right)
\end{align*}
from \eqref{p0}, \eqref{pei-5d}, \eqref{pei-7d} and \eqref{pei-9d}. \newline

Note that
$$\rho_a\Gamma_a=1,\quad -\frac{p_a}{\rho_a^2}\mathrm{d}\rho_a=\frac{1 }{\gamma_a\Gamma_a}\mathrm{d}\Gamma_a,$$
and we can also define the entropy $\mathbb{S}_a$ and the global entropy ${\mathbb S} = \sum_{a=1}^{N} {\mathbb S}_a$ as follows:
\begin{align*}
	&\mathrm{d} {\mathbb S}_a=\frac{c^2}{T_a}\left(\mathrm{d}\left(H_a-\frac{1}{\gamma_a}\right)+\frac{1 }{\gamma_a\Gamma_a}\mathrm{d}\Gamma_a\right),\nonumber
	\\ &\mathrm{d} {\mathbb S}=\sum_{a=1}^{N}\mathrm{d}S_a=\sum_{a=1}^{N}\frac{c^2}{T_a}\left(\mathrm{d}\left(H_a-\frac{1}{\gamma_a}\right)+\frac{1 }{\gamma_a\Gamma_a}\mathrm{d}\Gamma_a\right).
\end{align*}

\begin{definition} \label{D2.2}
	Time-dependent configuration $\{ (\mathbf{x}_a, \mathbf{w}_a, T_a) \}$ exhibits asymptotic thermo-mechanical flocking if and only if the following relations hold: 
	\begin{enumerate}
		\item
		Spatial cohesion holds:
		\[ \sup_{0\leq t <\infty} \max_{a, b} |\mathbf{x}_a(t)-\mathbf{x}_b(t)| < \infty. \]
		\item
		Velocity and temperature  alignment holds asymptotically: 
		\[
		\lim_{t \to \infty} \max_{a, b} |\mathbf{w}_a(t) -\mathbf{w}_b(t)| = 0, \quad \lim_{t \to \infty} \max_{a, b} |T_a(t)-T_b(t)| =0. 
		\]
	\end{enumerate}
\end{definition}
\vspace{0.2cm}

\noindent In the following lemma, we list several analytical properties of the relativistic model \eqref{RTCSmodel(same num)-1}.
\begin{lemma}\label{L2.2}
	For $\chi \in [4]$, let $\left\{\left(\mathbf{x}_a, \mathbf{v}_a, T_a\right)\right\}$ be a global solution to \eqref{RTCSmodel(same num)-1} satisfying 
	\[
	|\mathbb{M}(0)|<\infty, \quad {\mathbb E}_\chi(0)<\infty, \quad  \min_{a} T_a(0) > 0.
	\]
	Then, the following assertions hold.
	\begin{enumerate}
		\item
		The total momentum and energy are conserved:
		\[
		\mathbb{M}(t) = \mathbb{M}(0), \quad {\mathbb E}_\chi(t) =  {\mathbb E}_\chi(0), \quad t \geq 0.
		\]
		\item
		Then entropy is monotonically increasing:
		\[
		\frac{\mathrm{d} {\mathbb S}(t)}{\mathrm{d}t} = \frac{1}{2N}\sum_{a,b=1}^{N}\phi_{ab}\left|\frac{ \Gamma_{b}\mathbf{v}_{b}}{T_b}-\frac{\Gamma_{a}\mathbf{v}_{a}}{T_a}\right|^2
		+\frac{1}{2N}\sum_{a,b=1}^{N}\zeta_{ab}\Big| \frac{\Gamma_{b}}{T_b}-\frac{\Gamma_{a}}{T_a}\Big |^2 \geq 0, \quad t > 0.
		\]
		\item When the light speed $c$ is sufficiently large, temperatures are bounded below away from zero and bounded above: there exists a positive constant $K\geq1$ such that 
		\begin{equation*} \label{B-2-1}
			\frac{\prod_{a=1}^NT_a(0)}{K |\mathbb{E}_\chi(0)|^{N-1}} \leq T_a(t) \leq  K\mathbb{E}_\chi(0), \quad t \geq 0, \quad a \in [N].
		\end{equation*}
	\end{enumerate}
\end{lemma}
\begin{proof}
	\noindent (i)~We use the definition of ${\mathbb M}$ and the symmetry $\phi_{ab} = \phi_{ba}$ to find 
	\[
	\frac{\mathrm{d}}{\mathrm{d}t} {\mathbb M}  =\frac{1}{N}\sum_{a, b=1}^{N} \phi_{ab}\Big(\frac{ \Gamma_b\mathbf{w}_{b}}{F_bT_b}-\frac{\Gamma_a\mathbf{w}_{a}}{F_aT_a}\Big) =-\frac{1}{N}\sum_{a, b=1}^{N } \phi_{ab}\Big(\frac{ \Gamma_b\mathbf{w}_{b}}{F_bT_b}-\frac{\Gamma_a\mathbf{w}_{a}}{F_aT_a}\Big) = - \frac{\mathrm{d}}{\mathrm{d}t} {\mathbb M}.
	\]
	This implies 
	\[ \frac{\mathrm{d}}{\mathrm{d}t} {\mathbb M} = 0, \quad t > 0, \quad  \mbox{i.e.,} \quad {\mathbb M}(t) = {\mathbb M}(0), \quad t \geq 0. \]
	On the other hand, we sum up  $\eqref{A-4}$ over $a \in [N]$ and $\zeta_{ab} =\zeta_{ba}$ to find 
	\begin{align*}
		\begin{aligned}
			\frac{\mathrm{d}}{\mathrm{d}t} {\mathbb E}_\chi &= \frac{\mathrm{d}}{\mathrm{d}t} \sum_{a=1}^N \left[c^2\left(\Gamma_aH_a
			-1-\frac{1}{\gamma_a \Gamma_a}\right)\right] =\frac{1}{N}\sum_{a, b=1}^{N}\zeta_{ab}\Big(\frac{\Gamma_a }{T_a}-\frac{\Gamma_b }{T_b}\Big)  \\
			&= -\frac{1}{N}\sum_{a, b=1}^{N}\zeta_{ab}\Big(\frac{\Gamma_a }{T_a}-\frac{\Gamma_b }{T_b}\Big )= -\frac{\mathrm{d}}{\mathrm{d}t} {\mathbb E}_\chi.
		\end{aligned}
	\end{align*}
	This implies 
	\[ \frac{\mathrm{d}}{\mathrm{d}t} {\mathbb E}_\chi = 0, \quad \mbox{i.e.,} \quad  {\mathbb E}_\chi(t) =  {\mathbb E}_\chi(0), \quad t \geq 0. \]

	\noindent (ii)~We split its proof into two steps.\newline
	
	\noindent $\bullet$~Step A:~Recall the defining relation for ${\mathbb S}_a$:
	\[ \mathrm{d} {\mathbb S}_a=\frac{c^2}{T_a}\left(\mathrm{d}\left(H_a-\frac{1}{\gamma_a}\right)+\frac{1 }{\gamma_a\Gamma_a}\mathrm{d}\Gamma_a\right). \]
	Then, we claim that 
	\begin{equation} \label{B-3}
		\frac{\mathrm{d}{\mathbb S}_a}{\mathrm{d} t} =-\frac{1}{N}\sum_{b=1}^{N}\phi_{ab}\Big(\frac{ \Gamma_{b}\mathbf{v}_{b}}{T_b}-\frac{\Gamma_{a}\mathbf{v}_{a}}{T_a}\Big)\cdot \frac{\Gamma_{a}\mathbf{v}_{a}}{T_a}+\frac{1}{N}
		\sum_{b=1}^{N}\zeta_{ab}\Big(\frac{\Gamma_b }{T_b}-\frac{\Gamma_a }{T_a}\Big)\frac{\Gamma_a}{T_a}.
	\end{equation}
	{\it Proof of \eqref{B-3}}: ~We take the inner product of \eqref{A-4}$_2$ with $\displaystyle \frac{\Gamma_a \mathbf{v}_{a}}{T_a}$ to find  
	\begin{align*}
		& \Big ( \frac{\mathrm{d}}{\mathrm{d} t}\left(\Gamma_a\mathbf{v}_aH_a\right) \Big)  \cdot \frac{\Gamma_a \mathbf{v}_{a}}{T_a}   \\
		& \hspace{1cm} = \frac{1}{N}\sum_{b=1}^{N}\phi_{ab}\Big(\frac{\Gamma_b \mathbf{v}_{b}}{T_b}-\frac{\Gamma_a \mathbf{v}_{a}}{T_a}\Big)\cdot\frac{\Gamma_a \mathbf{v}_{a}}{T_a}  \nonumber\\
		& \hspace{1cm} =\frac{\Gamma_a^2 |\mathbf{v}_a|^2}{T_a}\frac{\mathrm{d}H_a}{\mathrm{d}\gamma_a}
		\frac{\mathrm{d}\gamma_a}{\mathrm{d}t}\nonumber+\frac{\Gamma_a }{T_a}H_a
		\left(|\mathbf{v}_a|^2\frac{\mathrm{d}\Gamma_a}{\mathrm{d}t}+\Gamma_a \mathbf{v}_{a}\cdot\frac{\mathrm{d}\mathbf{v}_a}{\mathrm{d}t}\right)
		\nonumber\\
		& \hspace{1cm} =\frac{\Gamma_a^2 |\mathbf{v}_a|^2}{T_a}\frac{\mathrm{d}H_a}{\mathrm{d}\gamma_a}\frac{\mathrm{d}\gamma_a}{\mathrm{d}t}
		+\frac{\Gamma_a^2 }{T_a}H_a
		\left(\frac{\Gamma_a^2|\mathbf{v}_a|^2}{c^2}+1 \right)\mathbf{v}_{a}\cdot\frac{\mathrm{d}\mathbf{v}_a}{\mathrm{d}t}
		\\
		&\hspace{1cm} =\frac{\Gamma_a^2 |\mathbf{v}_a|^2}{T_a}\frac{\mathrm{d}H_a}{\mathrm{d}\gamma_a}
		\frac{\mathrm{d}\gamma_a}{\mathrm{d}t}+\frac{\Gamma_a^4 }{T_a}H_a
		\mathbf{v}_{a}\cdot\frac{\mathrm{d}\mathbf{v}_a}{\mathrm{d}t}.\nonumber
	\end{align*}
	On the other hand, we multiply \eqref{A-4}$_3$ by $\frac{\Gamma_a }{T_a}$ to obtain
	\begin{align*}
		& \Big(  \frac{\mathrm{d}}{\mathrm{d} t}\left[c^2\left(\Gamma_aH_a
		-1-\frac{1}{\gamma_a \Gamma_a}\right)\right] \Big) \times \frac{\Gamma_a }{T_a} \nonumber\\
		& \hspace{1cm} = \frac{1}{N}\sum_{b=1}^{N}\zeta_{ab}\Big(\frac{\Gamma_a }{T_a}-\frac{\Gamma_b }{T_b}\Big)\frac{\Gamma_a }{T_a} \nonumber\\
		&\hspace{1cm}=\frac{\Gamma_a }{T_a}\frac{\mathrm{d}}{\mathrm{d} t}\left[c^2\left(\Gamma_aH_a
		-1-\frac{1}{\gamma_a \Gamma_a}\right)\right]  \\
		& \hspace{1cm} =\frac{c^2 }{T_a}\left(\Gamma_a^2\frac{\mathrm{d}H_a}{\mathrm{d}\gamma_a}
		+\frac{1}{\gamma^2_a}\right)\frac{\mathrm{d} \gamma_a}{\mathrm{d} t}+\frac{c^2 \Gamma_a }{T_a}\left(H_a+\frac{1}{\Gamma_a^2\gamma_a}\right)
		\frac{\mathrm{d} \Gamma_a}{\mathrm{d} t} \nonumber\\
		& \hspace{1cm} =\frac{c^2 }{T_a}\left(\Gamma_a^2\frac{\mathrm{d}H_a}{\mathrm{d}\gamma_a}
		+\frac{1}{\gamma^2_a}\right)\frac{\mathrm{d} \gamma_a}{\mathrm{d} t}+\frac{ \Gamma^4_a }{T_a}\left(H_a+\frac{1}{\Gamma_a^2\gamma_a}\right)
		\mathbf{v}_{a}\cdot\frac{\mathrm{d}\mathbf{v}_a}{\mathrm{d}t}.\nonumber
	\end{align*}
	Note that
	\begin{align*}
		\frac{c^2 }{T_a}\left(\Gamma_a^2\frac{\mathrm{d}H_a}{\mathrm{d}\gamma_a}
		+\frac{1}{\gamma^2_a}\right)-\frac{\Gamma_a^2}{T_a}|\mathbf{v}_a|^2
		\frac{\mathrm{d}H_a}{\mathrm{d}\gamma_a}&=
		\frac{c^2}{T_a}\left(\frac{\mathrm{d}H_a}{\mathrm{d}\gamma_a}+\frac{1}{\gamma^2_a}\right),
	\end{align*}
	
	\noindent $\bullet$~Step B: We sum up  \eqref{B-3} over all $a \in [N]$ to find 
	\begin{align}\label{entropygrow-mix}
		\begin{aligned}
			\frac{\mathrm{d}S}{\mathrm{d}t}&= \sum_{a} \frac{\mathrm{d}{\mathbb S}_a}{\mathrm{d}t} \\
			&=-\frac{1}{N}\sum_{a, b=1}^{N}\phi_{ab}\Big(\frac{ \Gamma_{b}\mathbf{v}_{b}}{T_b}-\frac{\Gamma_{a}\mathbf{v}_{a}}{T_a}\Big)\cdot\frac{\Gamma_{a}\mathbf{v}_{a}}{T_a} +\frac{1}{N} \sum_{a, b=1}^{N}\zeta_{ab}\Big(\frac{\Gamma_{b} }{T_b}-\frac{\Gamma_{a} }{T_a}\Big)\frac{\Gamma_{a}}{T_a} \\
			&=: {\mathcal I}_{21} + {\mathcal I}_{22}.
		\end{aligned}
	\end{align}
	Below, we estimate the terms ${\mathcal I}_{2i}$ one by one. \newline
	
	\noindent $\diamond$~Case B.1 (Estimate of ${\mathcal I}_{21}$):~By direct calculation, one use $\phi_{ab} = \phi_{ba}$ to find 
	\begin{align}
		\begin{aligned} \label{B-6}
			{\mathcal I}_{21} &:= -\frac{1}{N}\sum_{a, b=1}^{N}\phi_{ab}\Big(\frac{ \Gamma_{b}\mathbf{v}_{b}}{T_b}-\frac{\Gamma_{a}\mathbf{v}_{a}}{T_a}\Big)\cdot\frac{\Gamma_{a}\mathbf{v}_{a}}{T_a} = -\frac{1}{N}\sum_{a, b=1}^{N}\phi_{ba}\Big(\frac{ \Gamma_{a}\mathbf{v}_{a}}{T_a}-\frac{\Gamma_{b}\mathbf{v}_{b}}{T_b}\Big)\cdot\frac{\Gamma_{b}\mathbf{v}_{b}}{T_b} \\
			&= -\frac{1}{2N}\sum_{a,b=1}^{N}\phi_{ab}\left|\frac{ \Gamma_{b}\mathbf{v}_{b}}{T_b}-\frac{\Gamma_{a}\mathbf{v}_{a}}{T_a}\right|^2.
		\end{aligned}
	\end{align}
	
	\noindent $\diamond$~Case B.2 (Estimate of ${\mathcal I}_{22}$):~We use $\zeta_{ab} = \zeta_{ba}$ to get 
	\begin{align}
		\begin{aligned} \label{B-7}
			{\mathcal I}_{21} &:= \frac{1}{N} \sum_{a, b=1}^{N}\zeta_{ab}\Big(\frac{\Gamma_{b} }{T_b}-\frac{\Gamma_{a} }{T_a}\Big)\frac{\Gamma_{a}}{T_a} =  \frac{1}{N} \sum_{a, b=1}^{N}\zeta_{ba}\Big(\frac{\Gamma_{a} }{T_a}- 
			\frac{\Gamma_{b} }{T_b}\Big)\frac{\Gamma_{b}}{T_b}  \\
			&= -\frac{1}{2N}\sum_{a,b=1}^{N}\zeta_{ab}\Big| \frac{\Gamma_{b}}{T_b}-\frac{\Gamma_{a}}{T_a} \Big |^2.
		\end{aligned}
	\end{align}  
	In \eqref{entropygrow-mix}, we combine \eqref{B-6} and \eqref{B-7} to obtain the desired estimate. 
	
	\vspace{0.2cm}
	
	\noindent (iii)~ We consider upper and lower bound estimates separately. \newline
	
	\noindent$\bullet$ Case C.1 (Upper bound of the temperature):~We use \eqref{energy-m(same num)}, \eqref{B-2}, and Lemma \ref{L2.2} for $c$ sufficiently large to derive a desired estimate.  \newline
	\vspace{0.2cm}
	
	\noindent$\bullet$ Case C.2 (Lower bound of the temperature):~Note that
	\begin{align*}
		\begin{aligned}
			& \sum_{a=1}^N\frac{c^2}{T_a}\left[\frac{\mathrm{d}}{\mathrm{d}t}\left(H_a-\frac{1}{\gamma_a}\right)+\frac{1 }{\gamma_a}\frac{\mathrm{d}}{\mathrm{d}t}\ln{\Gamma_a}\right] \\
			& \hspace{1cm} =\sum_{a=1}^N\left[\frac{c^2}{T_a}\frac{\mathrm{d}}{\mathrm{d}t}\left(1
			+\frac{2\chi+1}{2\gamma_a}+\ldots\right)+\ln \Gamma_a\right] \\
			&\hspace{1cm} =\frac{1}{2}\sum_{a=1}^N\left[(2\chi+1)\frac{ \mathrm{d}}{\mathrm{d}t}\big[\ln{T_a}+c^2\mathcal{O}(\gamma_a^{-2})\big]-\frac{ \mathrm{d}}{\mathrm{d}t}\ln{\left(1-\frac{\mathbf{v}_a^2}{c^2}\right)}     \right]
			\\ 
			& \hspace{1cm} =\frac{1}{2}\frac{ \mathrm{d}}{\mathrm{d}t}\sum_{a=1}^N\left[
			\ln \Big\{\frac{T_a^{2\chi+1}}{1-\frac{\mathbf{v}_a^2}{c^2}}\Big\}
			+\mathcal{O}(c^{-2})\right].
		\end{aligned}
	\end{align*}
	Then, the relation \eqref{entropygrow-mix}$_2$ implies
	\begin{align*}
		\sum_{a=1}^N\left[
		\ln \Big\{\frac{T_a^{2\chi+1}(t)}{1-\frac{\mathbf{v}_a^2(t)}{c^2}}\Big\}
		+\mathcal{O}(c^{-2})\right]\geq \sum_{a=1}^N\left[
		\ln \Big\{\frac{T_a^{2\chi+1}(0)}{1-\frac{\mathbf{v}_a^2(0)}{c^2}}\Big\}
		+\mathcal{O}(c^{-2})\right],
	\end{align*}
	and it also leads to
	\begin{align*}
		\prod_{a=1}^N\left[
		\frac{T_a^{2\chi+1}(t)}{1-\frac{\mathbf{v}_a^2(t)}{c^2}}
		\exp\left\{\mathcal{O}(c^{-2})\right\}\right]\geq \prod_{a=1}^N\left[
		\frac{T_a^{2\chi+1}(0)}{1-\frac{\mathbf{v}_a^2(0)}{c^2}}
		\exp\left\{\mathcal{O}(c^{-2})\right\}\right].
	\end{align*}
	For a sufficiently large $c$, this further implies that there exists a positive constant $K\geq1$ independent of $c$ such that
	\begin{align*}
		\prod_{a=1}^NT_a(t)    \geq \frac{1}{K} \prod_{a=1}^NT_a(0).
	\end{align*}
\end{proof}
\begin{remark}
	For \eqref{RTCSmodel(same num)-1}, as in \cite{Ha-Kim-Ruggeri-ARMA-2020}, we can also define global quantities $\rho, v^j, \varepsilon$, and $p$  such that the sum of densities is the same as a single fluid. 
	Moreover, since total momentum and total energy are conserved, these global quantities satisfy the following set of relations: 
	\begin{align*} 
		\begin{aligned} \label{ME-RTCS}
			&\rho\Gamma=N,\qquad\Gamma \mathbf{v}H
			=\sum_{a=1}^N\left(\Gamma_a\mathbf{v}_aH_a\right)=\mathbf{0},\\
			&c^2\left(\Gamma H
			-1-\frac{1}{\gamma \Gamma}\right)=\sum_{a=1}^N\left[c^2\left(\Gamma_aH_a
			-1-\frac{1}{\gamma_a \Gamma_a}\right)\right]=\mbox{constant}.\nonumber
		\end{aligned}
	\end{align*} 
\end{remark}
As a direct application of Lemma \ref{L2.2}, we can also characterize asymptotic limits for \eqref{A-3}.
\begin{corollary} \label{C2.2}
	Let $\left\{\left(\mathbf{x}_a, \mathbf{w}_a, T_a\right)\right\}$ be a global solution to \eqref{RTCSmodel(same num)-1}. If asymptotic flocking occurs, i.e., 
	\[ \lim_{t \to \infty} \left(\mathbf{w}_a, T_a\right) = (\mathbf{w}^\infty,T^\infty), \]
	then, asymptotic limits are explicitly given by the following relations:
	\[
	\mathbf{w}^\infty =  \frac{\mathbb{M}(0)}{N} \quad \mbox{and} \quad \frac{(2\chi +1)}{2} T^{\infty} + \frac{1}{2N^2} |\mathbb{M}(0) |^2 + \frac{{\mathcal F}^{\infty}}{c^2} = \frac{{\mathbb E}_\chi(0)}{N}.
	\]
\end{corollary}
\begin{proof}
	\noindent It follows from Lemma \ref{L2.2} that 
	\begin{equation} \label{B-8}
		\frac{1}{N}\sum_{a=1}^{N}\mathbf{w}_a(t)=\frac{1}{N}\sum_{a=1}^{N}\mathbf{w}_a(0).
	\end{equation}
	We take $t \to \infty$ in \eqref{B-8} to get 
	\[ 
	\mathbf{w}^\infty = \frac{1}{N}\sum_{a=1}^{N}\mathbf{w}_a(0) = \frac{{\mathbb M}(0)}{N}. 
	\]
	By the conservation of energy, we have
	\begin{equation} \label{B-9}
		{\mathbb E}_\chi =  \sum_{a=1}^{N}\left(\frac{2\chi +1 }{2}{T_a(t)}+\frac{1}{2}|{\mathbf{w}_a}(t)|^2+\frac{{\mathcal F}_a(t)}{c^2}\right) =  {\mathbb E}_\chi(0). 
	\end{equation}
	Since we assume that $T_a$ and $\mathbf{w}_a$ have the same asymptotic limits independent of $a$, it follows from \eqref{def-fa} that ${\mathcal F}_a = {\mathcal F}_a(T_a,|\mathbf{w}_a|^2)$ has an asymptotic limit:  there exists ${\mathcal F}^{\infty}$ such that 
	\begin{equation*} \label{B-10}
		\lim_{t \to \infty} {\mathcal F}_a(t) = {\mathcal F}^{\infty}, \quad \forall~a \in [N]. 
	\end{equation*}
	In \eqref{B-9}, letting $t \to \infty$ in \eqref{B-9} and we use \eqref{B-8} to find 
	\[
	\frac{(2\chi +1)}{2} T^{\infty} + \frac{1}{2} |\mathbf{w}^{\infty} |^2 + \frac{{\mathcal F}^{\infty}}{c^2} = \frac{\mathbb{E}_\chi(0)}{N},
	\]
	or equivalently, 
	\[
	\frac{(2\chi +1)}{2} T^{\infty} + \frac{1}{2N^2} |\mathbb{M}(0) |^2 + \frac{{\mathcal F}^{\infty}}{c^2} = \frac{\mathbb{E}_\chi(0)}{N}.
	\]
\end{proof}

\section{Communication weights and relativistic corrections}\label{sec:3}
\setcounter{equation}{0}
In this section, we discuss three types of communication weight functions to be used in later sections, and present several preparatory estimates on the relativistic corrections ${\mathcal F}_a$ to be used in the classical limit $c \to \infty$.
\subsection{Communication weight functions} \label{sec:3.1}
In this subsection, we consider three different ansatzs for communication weight functions $\phi_{ab}$ and $\zeta_{ab}$.
\subsubsection{Uniform constants} 
Consider the constant communication functions which are independent on spatial distance:
\[ \phi_{ab} \equiv 1, \quad \zeta_{ab}\equiv1 \quad \forall~ a, b\in[N]. \]
In this case, the classical and relativistic TCS models read as follows: 
\begin{equation}\label{TCSmodel-3}
	\begin{cases}
		\displaystyle \frac{\mathrm{d}\mathbf{x}_a}{\mathrm{d} t}=\mathbf{v}_a, \quad a \in [N],\\
		\displaystyle \frac{\mathrm{d} \mathbf{v}_a}{\mathrm{d} t}
		=\frac{1}{N}\sum_{b=1}^{N}\Big(\frac{ \mathbf{v}_{b}}{T_b}-\frac{\mathbf{v}_{a}}{T_a}\Big),\\
		\displaystyle \frac{\mathrm{d}}{\mathrm{d} t}\left({\frac{2 \chi +1}{2} T_a}+\frac{1}{2}{|\mathbf{v}_a|^2}\right)=\frac{1}{N}\sum_{b=1}^{N}\Big(\frac{1 }{T_a}-\frac{1 }{T_b}\Big),
	\end{cases}
\end{equation}
and 
\begin{equation}\label{RTCSmodel(same num)-3}
	\begin{cases}
		\displaystyle \frac{\mathrm{d}\mathbf{x}_a}{\mathrm{d} t}=\frac{\mathbf{w}_a}{F_a},\quad a \in [N],\\
		\displaystyle \frac{\mathrm{d} \mathbf{w}_a}{\mathrm{d} t}
		=\frac{1}{N}\sum_{b=1}^{N}\Big(\frac{ \Gamma_b\mathbf{w}_{b}}{F_bT_b}-\frac{\Gamma_a\mathbf{w}_{a}}{F_aT_a}\Big),\\
		\displaystyle \frac{\mathrm{d}}{\mathrm{d} t}\left(\frac{2 \chi +1}{2}{T_a}+\frac{1}{2}|\mathbf{w}_a|^2+\frac{{\mathcal F}_a}{c^2}\right)=\frac{1}{N}\sum_{b=1}^{N}\Big(\frac{\Gamma_a }{T_a}-\frac{\Gamma_b }{T_b}\Big).
	\end{cases}
\end{equation}

\subsubsection{Small perturbations of uniform constant}
Consider constant communication weight matrix $(\phi_{ab)}$  which is a small perturbation around some given constant. More precisely, we set 
\[
\overline{\phi}:=\max_{a,b}{\phi_{ab}}, \quad \underline{\phi} :=\min_{a,b}{\phi_{ab}}, \quad \mbox{and}\quad \overline{\zeta}:=\max_{a,b}{\zeta_{ab}},\quad \underline{\zeta}:=\min_{a,b}{\zeta_{ab}}.
\]
Then we assume that $(\phi_{ab})$ satisfies 
\begin{align} \label{phiMm}
	| \overline{\phi} - \underline{\phi} |<\varepsilon\ll 1.
\end{align}
In this case, the classical and relativistic TCS models can be rewritten as 
\begin{equation}\label{TCSmodel-2}
	\begin{cases}
		\displaystyle \frac{\mathrm{d}\mathbf{x}_a}{\mathrm{d} t}=\mathbf{v}_a, \quad a \in [N],\\
		\displaystyle \frac{\mathrm{d} \mathbf{v}_a}{\mathrm{d} t}
		=\frac{\underline{\phi}}{N}\sum_{b=1}^{N}\Big(\frac{ \mathbf{v}_{b}}{T_b}-\frac{\mathbf{v}_{a}}{T_a}\Big)+\frac{1}{N}\sum_{b=1}^{N}(\phi_{ab}- \underline{\phi})\Big(\frac{ \mathbf{v}_{b}}{T_b}-\frac{\mathbf{v}_{a}}{T_a}\Big),\\
		\displaystyle \frac{\mathrm{d}}{\mathrm{d} t}\left(\frac{2 \chi +1}{2}{T_a}+\frac{1}{2}{|\mathbf{v}_a|^2}\right)=\frac{1}{N}\sum_{b=1}^{N} \zeta_{ab} \Big(\frac{1}{T_a}-\frac{1 }{T_b}\Big),
	\end{cases}
\end{equation}
and 
\begin{equation} \label{RTCSmodel(same num)-2}
	\begin{cases}
		\displaystyle \frac{\mathrm{d}\mathbf{x}_a}{\mathrm{d} t}=\frac{\mathbf{w}_a}{F_a},\quad a \in [N],\\
		\displaystyle \frac{\mathrm{d}\mathbf{w}_a}{\mathrm{d} t}
		=\frac{\underline{\phi}}{N}\sum_{b=1}^{N}\Big(\frac{\Gamma_b \mathbf{w}_{b}}{T_b}-\frac{\Gamma_a \mathbf{w}_{a}}{F_aT_a}\Big)+\frac{1}{N}\sum_{b=1}^{N}(\phi_{ab}- \underline{\phi})\Big(\frac{\Gamma_b \mathbf{w}_{b}}{F_bT_b}-\frac{\Gamma_a \mathbf{w}_{a}}{F_aT_a}\Big),\\
		\displaystyle \frac{\mathrm{d}}{\mathrm{d} t}\left(\frac{2\chi+1 }{2}{T_a}+\frac{1}{2}|\mathbf{w}|^2+\frac{{\mathcal F}_a}{c^2}\right)=\frac{1}{N}\sum_{b=1}^{N}\zeta_{ab}\Big(\frac{\Gamma_a }{T_a}-\frac{\Gamma_b }{T_b}\Big).
	\end{cases}
\end{equation}

\subsubsection{Bounded and monotone communication weight functions}
Two mother functions $\phi$ and $\zeta$ are assumed to be Lipschitz continuous and monotonic:
\begin{align}\label{properties of phi zeta}
	\begin{aligned}
		& \phi\in \mathrm{Lip} (\mathbb{R_+};\mathbb{R_+}),\quad (\phi(r_2)-\phi(r_1))(r_2-r_1)\le 0,\quad r_1,r_2\ge 0,
		\\&\zeta\in \mathrm{Lip} (\mathbb{R_+};\mathbb{R_+}),\quad (\zeta(r_2)-\zeta(r_1))(r_2-r_1)\le 0,\quad r_1,r_2\ge 0.
	\end{aligned}
\end{align}

\subsection{Estimates on relativistic corrections} \label{sec:3.2}
In this subsection, we study relativistic perturbations, when $c$ is sufficiently large.  We first study the temporal variations of $|\mathbf{v}_a|^2$ and $T_a$ in the following lemma.
\begin{lemma}\label{L3.1} Let $\{(\mathbf{x}_a, \mathbf{v}_a, T_a)\}$ be a global solution to \eqref{A-4}. Then, the following estimates hold. 
	\[ \Big|\frac{\mathrm{d} |\mathbf{v}_a|^2}{\mathrm{d} t}\Big|
	\leq\frac{C\left( \overline{\phi} + \overline{\zeta} \right)}{N \underline{T}^3}\sum_{a=1}^N|\mathbf{w}_a|^2, \quad \Big|\frac{\mathrm{d} T_a}{\mathrm{d} t}\Big|\leq \frac{C( \overline{\phi} + \overline{\zeta})}{N \underline{T}^4}
	\sum_{a=1}^N\left(|\mathbf{w}_a|^2+|T_a-T^{\infty}|\right).
	\]
\end{lemma}
\begin{proof}
	Since the proof are very lengthy, we split them into several steps. \newline
	
	\noindent $\bullet$~Step A (Derivation of equations for $\frac{\mathrm{d} |\mathbf{v}_a|^2}{\mathrm{d}t}$ and $\frac{\mathrm{d} T_a}{\mathrm{d} t}$): Recall that 
	\begin{equation} \label{C-0}
		\Gamma_a =\frac{1}{\sqrt{1-\frac{|\mathbf{v}_a|^2}{c^2}}},  \quad  \gamma_a :=\frac{c^2}{T_a}, \quad \frac{\mathrm{d}}{\mathrm{d} t}\left(\Gamma_a\mathbf{v}_aH_a\right) =\frac{1}{N}\sum_{b=1}^{N}\phi_{ab}\Big(\frac{\Gamma_b \mathbf{v}_{b}}{T_b}-\frac{\Gamma_a \mathbf{v}_{a}}{T_a}\Big).
	\end{equation}
	It follows from $\eqref{C-0}_1$ that  
	\begin{equation}\label{D Gamma}
		\frac{\mathrm{d}\Gamma_a}{\mathrm{d} t}=\frac{\Gamma_a^3\mathbf{v}_a}{c^2} \cdot \frac{\mathrm{d}\mathbf{v}_a}{\mathrm{d} t}  = \frac{\Gamma_a^3}{2c^2}  \frac{\mathrm{d}|\mathbf{v}_a|^2}{\mathrm{d} t} \quad \mbox{and} \quad  c^2\Gamma_a^2-\Gamma_a^2|\mathbf{v}_a|^2=c^2.
	\end{equation}
	Then, we use the relations \eqref{D Gamma} to get
	\[
	\frac{\mathrm{d}}{\mathrm{d} t}\left(\Gamma_a\mathbf{v}_aH_a\right) = \frac{\mathrm{d}
		\left(\Gamma_a\mathbf{v}_a\right)}{\mathrm{d} t} H_a +\Gamma_a\mathbf{v}_a\frac{\mathrm{d}H_a}{\mathrm{d} t} =H_a \Gamma_a \left(\frac{\mathrm{d}\mathbf{v}_a}{\mathrm{d} t }+ \frac{\Gamma_a^2}{2c^2} \frac{\mathrm{d}|\mathbf{v}_a|^2}{\mathrm{d} t}\mathbf{v}_a\right)+\Gamma_a\mathbf{v}_a\frac{\mathrm{d}H_a}{\mathrm{d}\gamma_a}
	\frac{\mathrm{d}\gamma_a}{\mathrm{d}t}.
	\]
	Thus, the equation $\eqref{C-0}_2$ can be rewritten as 
	\begin{equation} \label{C-1}
		\frac{1}{N}\sum_{b=1}^{N}\phi_{ab}\Big(\frac{\Gamma_b \mathbf{v}_{b}}{T_b}-\frac{\Gamma_a \mathbf{v}_{a}}{T_a}\Big) = H_a \Gamma_a \left(1+\Gamma_a^2 \frac{|\mathbf{v}_a|^2}{c^2}\right)\frac{\mathrm{d}\mathbf{v}_a}{\mathrm{d} t}+\Gamma_a\mathbf{v}_a\frac{\mathrm{d}H_a}{\mathrm{d}\gamma_a}
		\frac{\mathrm{d}\gamma_a}{\mathrm{d}t}.
	\end{equation}
	Now, we take an inner product of \eqref{C-1} with $\mathbf{v}_{a}$ and use $\eqref{D Gamma}_2$  to obtain
	\begin{align} 
		\begin{aligned} \label{vv-D3}
			&\frac{1}{N}\sum_{b=1}^{N}\phi_{ab}\Big(\frac{\Gamma_b \mathbf{v}_{b}}{T_b}-\frac{\Gamma_a \mathbf{v}_{a}}{T_a}\Big)\cdot\mathbf{v}_a  \\
			& \hspace{1cm} =\frac{1}{2}H_a\Gamma_a\left(1+\frac{\Gamma_a^2|\mathbf{v}_a|^2}{c^2}\right)\frac{\mathrm{d}|\mathbf{v}_a|^2}{\mathrm{d} t}+\Gamma_a|\mathbf{v}_a|^2\frac{\mathrm{d}H_a}{\mathrm{d}\gamma_a}\frac{\mathrm{d}\gamma_a}{\mathrm{d}t} \\
			& \hspace{1cm} =\frac{1}{2}H_a\Gamma^3_a\frac{\mathrm{d}|\mathbf{v}_a|^2}{\mathrm{d} t}
			+\Gamma_a|\mathbf{v}_a|^2\frac{\mathrm{d}H_a}{\mathrm{d}\gamma_a}
			\frac{\mathrm{d}\gamma_a}{\mathrm{d}t}.
		\end{aligned}
	\end{align}
	On the other hand, we recall the energy equation:
	\begin{equation} \label{C-2}
		\frac{\mathrm{d}}{\mathrm{d} t}\left[c^2\left(\Gamma_aH_a
		-1-\frac{1}{\gamma_a \Gamma_a}\right)\right]=\frac{1}{N}\sum_{b=1}^{N}\zeta_{ab}\Big(\frac{\Gamma_a }{T_a}-\frac{\Gamma_b }{T_b}\Big).    
	\end{equation}
	We rewrite the left-hand side of \eqref{C-2} using \eqref{D Gamma} as follows: \newline
	\begin{align}
		\begin{aligned}
			& \frac{\mathrm{d}}{\mathrm{d} t}\left[c^2\left(\Gamma_aH_a
			-1-\frac{1}{\gamma_a \Gamma_a}\right)\right] = c^2  \frac{\mathrm{d}}{\mathrm{d} t} \left(\Gamma_aH_a + \frac{1}{\gamma_a \Gamma_a}\right) \\
			& \hspace{1cm} =  c^2  \left( \frac{\mathrm{d}\Gamma_a}{\mathrm{d}t} H_a + \Gamma_a \frac{dH_a}{\mathrm{d}t}  + \frac{1}{\gamma^2_a \Gamma_a} \frac{d\gamma_a}{\mathrm{d}t} + \frac{1}{\gamma_a \Gamma^2_a} \frac{\mathrm{d}\Gamma_a}{\mathrm{d}t} \right) \\
			& \hspace{1cm} =c^2 \left(\Gamma_a\frac{\mathrm{d} H_a}{\mathrm{d} \gamma_a}
			+\frac{1}{\Gamma_a\gamma^2_a}\right)\frac{\mathrm{d} \gamma_a}{\mathrm{d} t}+c^2\left(H_a+\frac{1}{\Gamma_a^2\gamma_a}\right)
			\frac{\mathrm{d} \Gamma_a}{\mathrm{d} t}\\
			& \hspace{1cm} =c^2\left(\Gamma_a\frac{\mathrm{d} H_a}{\mathrm{d} \gamma_a}
			+\frac{1}{\Gamma_a\gamma^2_a}\right)\frac{\mathrm{d} \gamma_a}{\mathrm{d} t}+\frac{1}{2}\Gamma_a \left(\Gamma_a^2H_a+\frac{1}{\gamma_a}\right)
			\frac{\mathrm{d}|\mathbf{v}_a|^2}{\mathrm{d}t}.\nonumber
		\end{aligned}
	\end{align}
	Thus, the relation \eqref{C-2} can be rewritten as 
	\begin{equation}  \label{C-3}
		\frac{1}{N}\sum_{b=1}^{N}\zeta_{ab}\Big(\frac{\Gamma_a }{T_a}-\frac{\Gamma_b }{T_b}\Big) =c^2\left(\Gamma_a\frac{\mathrm{d} H_a}{\mathrm{d} \gamma_a}
		+\frac{1}{\Gamma_a\gamma^2_a}\right)\frac{\mathrm{d} \gamma_a}{\mathrm{d} t}+\frac{1}{2}\Gamma_a \left(\Gamma_a^2H_a+\frac{1}{\gamma_a}\right)
		\frac{\mathrm{d}|\mathbf{v}_a|^2}{\mathrm{d}t}.
	\end{equation}
	We combine \eqref{vv-D3} and \eqref{C-3} as 
	\begin{equation} \label{C-4}
		\begin{cases}
			\displaystyle a_1\frac{\mathrm{d} \gamma_a}{\mathrm{d} t}+b_1\frac{\mathrm{d}|\mathbf{v}_a|^2}{\mathrm{d}t}=c_1,\\
			\displaystyle  a_2\frac{\mathrm{d} \gamma_a}{\mathrm{d} t}+b_2\frac{\mathrm{d}|\mathbf{v}_a|^2}{\mathrm{d}t}=c_2,
		\end{cases}
	\end{equation}
	where $a_i, b_i$ and $c_i$ are given as follows.  
	\begin{align*}
		a_1 &=\Gamma_a|\mathbf{v}_a|^2\frac{\mathrm{d} H_a}{\mathrm{d} \gamma_a}, \quad a_2=c^2\left(\Gamma_a\frac{\mathrm{d} H_a}{\mathrm{d} \gamma_a}
		+\frac{1}{\Gamma_a\gamma^2_a}\right), \\
		b_1 &=\frac{1}{2}H_a\Gamma^3_a, \quad b_2 =\frac{1}{2}\Gamma_a \left(\Gamma_a^2H_a+\frac{1}{\gamma_a}\right), \\
		c_1 &=\frac{1}{N}\sum_{b=1}^{N}\phi_{ab}\Big(\frac{\Gamma_b \mathbf{v}_{b}}{T_b}-\frac{\Gamma_a \mathbf{v}_{a}}{T_a}\Big)\cdot\mathbf{v}_a, \quad c_2=\frac{1}{N}\sum_{b=1}^{N}\zeta_{ab}
		\Big(\frac{\Gamma_a }{T_a}-\frac{\Gamma_b }{T_b}\Big).
	\end{align*}
	It follows from \eqref{C-4} that 
	\begin{equation*}
		\frac{\mathrm{d} |\mathbf{v}_a|^2}{\mathrm{d} t}=\frac{a_2c_1-a_1c_2}{a_2b_1-a_1b_2} \quad \mbox{and} \quad 
		\frac{\mathrm{d} \gamma_a}{\mathrm{d} t}= \frac{c_1b_2-c_2b_1}{a_1b_2-a_2b_1}.
	\end{equation*}

	\vspace{0.2cm}
	
	\noindent $\bullet$~Step B: we claim that 
	\begin{align}
		\begin{aligned} \label{C-5}
			& (i)~a_1b_2-a_2b_1 = \frac{2\chi +1}{4c^2}T_a^2+\mathcal{O}(c^{-4}).  \\
			& (ii)~c_1b_2-c_2b_1 \leq \frac{C( \overline{\phi} + \overline{\zeta})}{2N \underline{T}^2}\sum_{b=1}^N\left(|\mathbf{v}_b|^2+|T_b-T_{\infty}|\right).  \\
			& (iii)~a_1c_2-a_2c_1 \leq \frac{C\left(\overline{\phi} +\overline{\zeta} \right)}{N \underline{T} c^2}\sum_{b=1}^N|\mathbf{v}_b|^2.
		\end{aligned}
	\end{align}
	Below, we provide estimates for \eqref{C-5} one by one. \newline
	
	\noindent $\diamond$~Case D.1 (Estimate of $\eqref{C-1}_1$):~By \eqref{H_a(same num)}, we obtain
	\begin{align*}
		\frac{\mathrm{d}H_a}{\mathrm{d}\gamma_a}=\frac{\mathrm{d}}{\mathrm{d}\gamma_a}
		\left(1+\frac{2 \chi +3}{2\gamma_a}+\mathcal{O}(\gamma^{-2}_a)\right)
		=-\frac{2 \chi +3}{2\gamma_a^2}+\mathcal{O}(\gamma_a^{-3}).
	\end{align*}
	We use the above relation to find 
	\begin{align} 
		\begin{aligned} \label{a1b2-a2b1}
			&2(a_1b_2-a_2b_1) \\
			& \hspace{0.5cm} =-c^2\Gamma_a^2\left(\frac{\mathrm{d} H_a}{\mathrm{d} \gamma_a}+\frac{1}{\gamma_a^2}\right)H_a+\Gamma_a^2|\mathbf{v}_a|^2\frac{1}{\gamma_a}\frac{\mathrm{d} H_a}{\mathrm{d} \gamma_a} \\
			& \hspace{0.5cm}  =\Gamma_a^2\left[-c^2\left(-\frac{2 \chi +1}{2\gamma_a^2}\right)
			\left(1+\frac{2 \chi +3}{2\gamma_a}\right)+\frac{1}{\gamma_a}|\mathbf{v}_a|^2
			\left(-\frac{2 \chi +3}{2\gamma_a^2}\right)\right]+c^2\mathcal{O}(\gamma_a^{-3}) \\
			& \hspace{0.5cm} = \frac{2 \chi +1}{2\gamma_a^2}c^2\Gamma_a^2+c^2\mathcal{O}(\gamma_a^{-3}) =\frac{2 \chi +1}{2c^2}T_a^2+\mathcal{O}(c^{-4}).
		\end{aligned}
	\end{align}
	\noindent $\diamond$~Case D.2 (Estimate of $\eqref{C-1}_2$):~Similarly, we have
	\begin{align}
		\label{c1b2-c2b1}
		&2(c_1b_2-c_2b_1) \nonumber \\
		&  \hspace{0.5cm}  =\Gamma_a\left[\Gamma_a^2\left(1+\frac{2 \chi +3}{2\gamma_a}\right)+\frac{1}{\gamma_a}
		+\mathcal{O}(\gamma_a^{-2})\right]\frac{1}{N}\sum_{b=1}^{N}\phi_{ab}\Big(\frac{\Gamma_b \mathbf{v}_{b}}{T_b}-\frac{\Gamma_a \mathbf{v}_{a}}{T_a}\Big)\cdot\mathbf{v}_a \nonumber\\
		& \hspace{0.5cm} -\Gamma^3_a\left[1+\frac{2 \chi +3}{2\gamma_a}+\mathcal{O}(\gamma_a^{-2})\right]
		\frac{1}{N}\sum_{b=1}^{N}\zeta_{ab}\Big(\frac{\Gamma_a }{T_a}-\frac{\Gamma_b }{T_b}\Big)
		\\ 
		& \hspace{0.5cm} \le \frac{C \overline{\phi}}{N \underline{T} }\sum_{b=1}^N|\mathbf{v}_b|^2
		+\frac{C \overline{\zeta}}{N}\sum_{b=1}^N\left|\frac{\Gamma_a-\Gamma_b }{T_a}-\frac{(T_b-T^{\infty})-(T_a-T^{\infty}) }{T_aT_b}\right|\nonumber\\
		& \hspace{0.5cm} \le \frac{C( \overline{\phi} + \overline{\zeta})}{N \underline{T}^2}\sum_{b=1}^N\left(|\mathbf{v}_b|^2+|T_b-T^{\infty}|\right).\nonumber
	\end{align}
	\noindent $\diamond$~Case D.3 (Estimate of $\eqref{C-1}_3$): By direct calculations, we have
	\begin{align*} 
		\begin{aligned}
			\label{a2c1-a1c2}
			&a_1c_2-a_2c_1 \nonumber\\
			& \hspace{0.5cm} = -c^2\left[-\frac{2 \chi +3}{2\gamma_a^2}\Gamma_a+\frac{1}{\Gamma_a\gamma^2_a}+\mathcal{O}(\gamma_a^{-3})\right]\frac{1}{N}\sum_{b=1}^{N}\phi_{ab}\Big(\frac{\Gamma_b \mathbf{v}_{b}}{T_b}-\frac{\Gamma_a \mathbf{v}_{a}}{T_a}\Big)\cdot\mathbf{v}_a
			\\
			& \hspace{0.5cm} +\left[\frac{2 \chi +3}{2\gamma_a^2}+\mathcal{O}(\gamma_a^{-3})\right]\Gamma_a|\mathbf{v}_a|^2\frac{1}{N}\sum_{b=1}^{N}\zeta_{ab}\Big(\frac{\Gamma_a }{T_a}-\frac{\Gamma_b }{T_b}\Big) \nonumber
			\\ 
			& \hspace{0.5cm} \leq \frac{C\left( \overline{\phi} + \overline{\zeta} \right)}{N \underline{T} c^2}\sum_{b=1}^N|\mathbf{v}_b|^2.\nonumber
		\end{aligned}
	\end{align*}
	
	\noindent $\bullet$~Step C:  It follows from $\eqref{C-0}_2$ that 
	\[  \frac{\mathrm{d}\gamma_a}{\mathrm{d}t} = -\frac{c^2}{T^2_a} \frac{\mathrm{d}T_a}{\mathrm{d}t}, \quad \mbox{i.e.} \quad  \frac{\mathrm{d}T_a}{\mathrm{d}t} = -\frac{T_a^2}{c^2}  \frac{\mathrm{d}\gamma_a}{\mathrm{d}t} =  -\frac{T_a^2}{c^2}   \frac{c_1b_2-c_2b_1}{a_1b_2-a_2b_1}.  \]
	Now, we use the relations \eqref{a1b2-a2b1} and  \eqref{c1b2-c2b1} to obtain estimates:
	\begin{align*}
		\begin{aligned}
			& \Big|\frac{\mathrm{d} |\mathbf{v}_a|^2}{\mathrm{d} t}\Big|\leq \frac{\frac{C\left(\overline{\phi}+\overline{\zeta}\right)}{N\underline{T}c^2}
				\sum_{b=1}^N|\mathbf{v}_b|^2}{\frac{2\chi+1}{2c^2}T_a^2+\mathcal{O}(c^{-4})}
			\leq\frac{C\left(\overline{\phi}+\overline{\zeta}\right)}{N\underline{T}^3}\sum_{b=1}^N|\mathbf{w}_b|^2  ,
			\\ 
			& \Big|\frac{\mathrm{d} T_a}{\mathrm{d} t}\Big|\leq \frac{\frac{C(\overline{\phi}+\overline{\zeta})}{N\underline{T}^2}\sum_{b=1}^N\left(|\mathbf{v}_b|^2+|T_b-T_{\infty}|\right)}
			{\frac{2\chi+1}{2}T_a^2+\mathcal{O}(c^{-2})} \leq\frac{C(\overline{\phi}+\overline{\zeta})}{N\underline{T}^4}
			\sum_{b=1}^N\Big (|\mathbf{w}_b|^2+|T_b-T_{\infty}| \Big).
		\end{aligned}
	\end{align*}
\end{proof}
\noindent In the next lemma, we estimate on the relativistic fluctuations and their temporal derivatives. 
\begin{lemma}
\label{L3.2}
Let $\{(\mathbf{x}_a, \mathbf{w}_a, T_a)\}$ be a global solution to \eqref{RTCSmodel(same num)-1}. Then, for light speed $c$ large enough, there exist positive constants $C_1$ and $C_2$ independent of $c$ such that 
\[ \sum_{a=1}^N| {\mathcal F}_a|\le C_1\mathbb{E}^2_{\chi}, \quad  \sum_{a=1}^N\left|\frac{\mathrm{d} {\mathcal F}_a}{\mathrm{d}t}\right|\leq C_2\sum_{a=1}^N\Big(|\mathbf{w}_a|^2+|T_a-T_{\infty}| \Big), \quad 
\frac{\Gamma_a}{F_a}-1=-\frac{(2\chi +3)T_a}{2c^2}+\mathcal{O}(c^{-4}).
\]
\end{lemma}
\begin{proof} 
\noindent (i)~Note that
$$\Gamma_a=1+\frac{|\mathbf{v}_a|^2}{2c^2}+\frac{3|\mathbf{v}_a|^2}{8c^4}+\mathcal{O}(c^{-6}),\quad
\frac{1}{\Gamma_a}=1-\frac{|\mathbf{v}_a|^2}{2c^2}+\mathcal{O}(c^{-4}).$$ 
Thus, it follows from \eqref{H_a(same num)} and the definition of ${\mathcal F}_a$ in \eqref{def-fa} that 
\begin{align*}
\sum_{a=1}^N\frac{{\mathcal F}_a}{c^2} &= \sum_{a=1}^N \left[c^2\left(\Gamma_aH_a
-1-\frac{1}{\gamma_a \Gamma_a}\right)-\frac{2\chi +1}{2}{T_a}-\frac{1}{2}|\mathbf{w}_a|^2\right] \\
&=\sum_{a=1}^Nc^2\left[\left(1+\frac{|\mathbf{v}_a|^2}{2c^2}+\frac{3|\mathbf{v}_a|^4}{8c^4}\right)
\left(1+\frac{(2 \chi +3)T_a}{2c^2}+\frac{\mathcal{O}(T_a^2)}{c^{4}}\right)
-1-\left(1-\frac{|\mathbf{v}_a|^2}{2c^2}\right)\frac{T_a}{c^2}\right] \\
&-\sum_{a=1}^{N}\left[\frac{2\chi+1}{2}T_a+\frac{|\mathbf{v}_a|^2}{2}
\left(1+\frac{|\mathbf{v}_a|^2}{2c^2}\right)^2\left(1+\frac{(2 \chi +3)T_a}{2c^2}\right)^2
\right]+\mathcal{O}(c^{-4})
\\
&= \sum_{a=1}^{N} \frac{\mathcal{O}(T_a^2)}{c^{2}} +\sum_{a=1}^{N}\frac{\mathcal{O}(|\mathbf{v}_a|^4)}{c^{2}}+\mathcal{O}(c^{-4}).
\end{align*}
This implies the first estimate. On the other hand, we use \eqref{energy-m(same num)} and the conservation of energy $\mathbb{E}_\chi(t)$ to find 
\[ |\mathbf{v}_a|^2\le 3\mathbb{E}_\chi(0), \]
when $c$ is large enough. Then, the expression of  ${\mathcal F}_a$ implies
\begin{align*}
\sum_{a=1}^{N}\left|\frac{\mathrm{d} {\mathcal F}_a}{\mathrm{d}t} \right| 
&\leq \sum_{a=1}^{N}\left|\left(\mathcal{O}(T_a)-\frac{(2 \chi +1)|\mathbf{v}_a|^2}{4}+\mathcal{O}(c^{-2})\right)\frac{\mathrm{d} T_a}{\mathrm{d} t}\right| \\
& +\sum_{a=1}^{N}\left|\left(-\frac{(2\chi+1)T_a}{4}-\frac{|\mathbf{v}_a|^2}{8}+\mathcal{O}(c^{-2})\right)\frac{\mathrm{d} |\mathbf{v}_a|^2}{\mathrm{d} t}\right|
\\
& \le C \sum_{a=1}^{N}\left(T_a+|\mathbf{v}_a|^2\right)\left(\left|\frac{\mathrm{d}T_a}{\mathrm{d}t} \right| +\left|\frac{\mathrm{d}|\mathbf{v}_a|^2}{\mathrm{d}t} \right| \right) \leq  C_2 \mathbb{E}_\chi \left(\left|\frac{\mathrm{d}T_a}{\mathrm{d}t} \right| +\left|\frac{\mathrm{d}|\mathbf{v}_a|^2}{\mathrm{d}t} \right| \right).\nonumber
\end{align*}
This inequality and Lemma \ref{L3.1} yield the second estimate.  \newline

\noindent (ii)~We use \eqref{H_a(same num)} to have
\begin{align*}
\frac{\Gamma_a}{F_a}-1=\frac{1-H_a}{H_a}
=-\frac{\frac{(2 \chi +3)T_a}{2c^2}+\mathcal{O}(c^{-4})}{1+\frac{(2 \chi +3)T_a}{2c^2}
	+\mathcal{O}(c^{-4})}=-\frac{(2 \chi +3)T_a}{2c^2}+\mathcal{O}(c^{-4}).
	\end{align*}
\end{proof}
\noindent We set the ensemble averages of velocity and temperature configurations as follows.
\[  \langle T \rangle := \frac{1}{N} \sum_{a=1}^{N} \mathbf{T}_a, \quad   \langle |{\mathbf v}|^2 \rangle := \frac{1}{N} \sum_{a=1}^{N} |\mathbf{v}_a|^2, \quad \langle |{\mathbf w}|^2 \rangle := \frac{1}{N} \sum_{a=1}^{N} |\mathbf{w}_a|^2. \]
In the next lemma, we continue to estimate $|T_\infty- \langle T \rangle|$.
\begin{lemma}\label{L3.3}
Let $\{(\mathbf{x}_a, \mathbf{w}_a, T_a)\}$ be a global solution to \eqref{RTCSmodel(same num)-1} satisfying 
\[ \langle \mathbf{w}(t) \rangle =  0, \quad \forall~t \geq 0. \]
Then we have
\begin{align*}
\begin{aligned}
	& (i)~|T^{\infty} -\langle T \rangle|\le \frac{1}{(2\chi+1)N} \sum_{a=1}^N|\mathbf{w}_a|^2 +\frac{C}{Nc^2}\sum_{a=1}^N\left(|\mathbf{w}_a|^2+|T_a-T_{\infty}|\right). \\
	& (ii)~   \sum_{a,b=1}^{N} (T_a-T_b)^2\le2N\left[1+\mathcal{O}(c^{-2})\right]\sum_{a=1}^{N} \left(|T_a-T^\infty|^2-\frac{1}{(2 \chi +1)^2N^2}\sum_{a=1}^N|\mathbf{w}_a|^2 \right).
\end{aligned}
\end{align*}
\end{lemma}
\begin{proof}
\noindent (i)~ We use the conservation of energy, \eqref{energy-m(same num)} and \eqref{def-fa}, and we further apply the mean value theorem for the error function ${\mathcal F}_a(T_a,|\mathbf{w}_a|^2)$ with respect to  its variables to get 
\begin{align}\label{Tinfty-varT b}
\begin{aligned}
	\frac{2\chi +1}{2}|T^\infty- \langle T \rangle | &= \frac{1}{N}\left| \frac{1}{2}\sum_{a=1}^N|\mathbf{w}_a|^2  +
	\sum_{a=1}^{N}\frac{{\mathcal F}_a(T_a,|\mathbf{w}_a|^2)- {\mathcal F}_a(T_{\infty},0)}{c^2}\right|\\
	\leq& \frac{1}{2N} \sum_{a=1}^N|\mathbf{w}_a|^2 +\frac{C}{Nc^2}\sum_{a=1}^N \Big (|\mathbf{w}_a|^2+|T_a-T_{\infty}| \Big),
\end{aligned}
\end{align}
where $C$ is a uniform constant independent of $c$.

\noindent (ii)~We first use $\displaystyle \sum_{a=1}^{N} (T_a- \langle T \rangle ) = 0$ to see
\[ \sum_{a=1}^{N} (T_a-T^\infty)(T^\infty- \langle T \rangle) = - |T^\infty- \langle T \rangle|^2.
\]
This yields
\begin{align} 
\begin{aligned} \label{Ta-Tb(TCS)}
	\sum_{a,b=1}^{N} |T_a-T_b|^2 &= \sum_{a,b=1}^{N} \left|(T_a- \langle T \rangle )-(T_b- \langle T \rangle )\right|^2\\
	&= 2N \sum_{a,b=1}^{N} \left| T_a- \langle T \rangle \right|^2 -2\sum_{a,b=1}^{N} (T_a- \langle T \rangle)(T_b- \langle T \rangle)           \\
	&= 2N\sum_{a=1}^{N} |(T_a-T^\infty)+(T^\infty- \langle T \rangle )|^2-2\sum_{a,b=1}^{N} (T_a- \langle T \rangle )(T_b- \langle T \rangle ) 
	\\
	&= 2N\sum_{a=1}^{N}\Big ( |T_a-T^\infty|^2+2(T_a-T_\infty)(T_\infty- \langle T \rangle)+(T^\infty- \langle T \rangle)^2 \Big)
	\\
	&= 2N\sum_{a=1}^{N} \Big(|T_a-T^\infty|^2-|T^\infty- \langle T \rangle |^2 \Big).
\end{aligned} 
\end{align}
Then,  \eqref{Ta-Tb(TCS)} and \eqref{Tinfty-varT b} imply
\begin{align}\label{Ta-Tb(RTCS)}
\sum_{a,b=1}^{N} (T_a-T_b)^2\le2N \Big ( 1+\mathcal{O}(c^{-2}) \Big) \sum_{a=1}^{N} \left(|T_a-T^\infty|^2-\frac{1}{(2\chi+1)^2N^2} \sum_{a=1}^N|\mathbf{w}_a|^2\right).
\end{align}
\end{proof}

%
%

\section{Asymptotic flocking of classical TCS model}\label{sec:4}
\setcounter{equation}{0}
In this section, we investigate the asymptotic flocking dynamics of the TCS model \eqref{A-3} depending on the different types of communication weight functions described in Section \ref{sec:3}. First, we set
\begin{align*}
\begin{aligned}
\mathbf{X} &:= (\mathbf{x}_1, \cdots, \mathbf{x}_N), \quad \mathbf{V} := (\mathbf{v}_1, \cdots, \mathbf{v}_N), \quad \mathbf{\hat T} := \mathbf{T} - T^\infty {\mathbf 1} := (T_1  - T^\infty,  \cdots, T_N  - T^\infty), \\
\| \mathbf{X} \| &:=\left(\sum_{a=1}^N|\mathbf{x}_a|^2\right)^{\frac{1}{2}},\quad \| \mathbf{V} \|:=\left( \sum_{a=1}^N|\mathbf{v}_a|^2\right)^{\frac{1}{2}}, \quad 
\| \mathbf{\hat T} \| :=  \Big( \frac{2\chi+1}{2}\sum_{a=1}^N|T_a -T^\infty|^2 \Big)^{\frac{1}{2}}.
\end{aligned}
\end{align*}
\subsection{Flocking dynamics I} \label{sec:4.1}
In this subsection, we show the asymptotic flocking dynamics of the simplified TCS model \eqref{TCSmodel-3} for arbitrary initial data. First, we derive a system of dissipative inequalities for position and velocity.
\begin{lemma} \label{L4.1}
Let $\{(\mathbf{x}_a, \mathbf{v}_a, T_a) \}$ be a global solution to \eqref{TCSmodel-3} with zero initial momentum:
\[ \langle \mathbf{v}(0) \rangle= 0. \]
Then, $\| \mathbf{X} \|$ and $\| \mathbf{V} \|$ satisfy
\begin{equation} \label{D-1}
\frac{\mathrm{d} \| \mathbf{X} \|}{\mathrm{d}t} \leq \| \mathbf{V} \|, \qquad \frac{\mathrm{d}| \mathbf{V} \|}{\mathrm{d}t}\le -\frac{1}{\overline{T}} \| \mathbf{V} \|, \quad \mbox{a.e.}~ t > 0.
\end{equation}
\end{lemma}
\begin{proof} 
By Lemma \ref{L2.1}, the total momentum is zero:
\begin{equation} \label{D-1-1}
\langle {\mathbf v}(t) \rangle = 0, \quad \forall~t \geq 0.
\end{equation}
\noindent (i)~The first differential inequality follows directly from the Cauchy-Schwarz inequality.  \newline

\noindent (ii)~We take the inner product $\mathbf{v}_a$ with \eqref{TCSmodel-3}$_2$ to obtain
\begin{align}\label{dv(TCS 3)}
\begin{aligned}
	\frac{\mathrm{d}}{\mathrm{d}t} \sum_{a=1}^{N}|\mathbf{v}_a|^2 = \frac{2}{N}\sum_{a,b=1}^{N}\left(\frac{\mathbf{v}_a\cdot\mathbf{v}_b}{T_b}-\frac{|\mathbf{v}_a|^2}{T_a}\right)
	\le  -\frac{2}{\overline{T}}\sum_{a=1}^{N} |\mathbf{v}_a|^2,
\end{aligned}
\end{align}
where we used \eqref{D-1-1} and Lemma \ref{L4.1}:
\[
\sum_{a,b=1}^{N} \frac{\mathbf{v}_a\cdot\mathbf{v}_b}{T_b} = \Big( \sum_{a=1}^{N} \mathbf{v}_a \Big) \cdot \Big(  \sum_{b=1}^{N} \frac{\mathbf{v}_b}{T_b}\Big) = 0, \quad \sup_{0 \leq t < \infty} \max_{a} T_a(t) \leq \overline{T}.
\]
\end{proof}
As a direct application of Lemma \ref{L4.1}, we obtain the exponential decay of $\| {\mathbf V} \|$ and  the uniform boundedness of $\| {\mathbf X} \|$. 
\begin{corollary}
\label{C4.1}
Let $\left\{\left(\mathbf{x}_a, \mathbf{v}_a, T_a\right)\right\}$ be a global solution to  \eqref{TCSmodel-3} with zero initial momentum:
\[ \langle \mathbf{v}(0) \rangle= 0. \]
Then, one has 
\[
\| \mathbf{V}(t) \|  \le  \| \mathbf{V}(0) \|  e^{-\frac{t}{\overline{T}}}, \quad \forall~t \geq 0, \quad \sup_{0 \leq t < \infty}  \| \mathbf{X}(t) \| \leq \| \mathbf{X}(0) \|  +  \| {\mathbf V}(0) \| \overline{T}.
\]
\end{corollary}
\begin{proof} (i)~ It follows from $\eqref{D-1}_2$ that 
\[  \| \mathbf{V}(t) \| \le  \| \mathbf{V}(0) \|e^{-\frac{t}{\overline{T}}}. \]
Again, we use $\eqref{D-1}_1$  and the above estimate to find 
\begin{equation} \label{D-1-2}
\frac{\mathrm{d} \| \mathbf{X}(t) \|}{\mathrm{d}t} \leq \| \mathbf{V}(t) \| \leq  \| \mathbf{V}(0) \|e^{-\frac{t}{\overline{T}}}.
\end{equation}
We integrate \eqref{D-1-2} to find 
\[
\| \mathbf{X}(t) \| \leq \| \mathbf{X}(0) \| +  \| {\mathbf V}(0) \| \overline{T}(1-e^{-\frac{t}{\overline{T}}}) \leq \| \mathbf{X}(0) \| +  \| {\mathbf V}(0) \| \overline{T}.
\]
Now, we take a supremum the above relation over all $t$ to get the desired estimate. 
\end{proof}
Next, we derive an exponential decay for temperatures.  For this, we choose $A_1$ and $\lambda_1$ to satisfy 
\begin{align}
\begin{aligned} \label{D-2}
A_1 &>\frac{ \| \mathbf{V}(0) \|^2 }{(2 \chi +1)^2N^2 \overline{T}}+ \overline{T}
+\frac{(N-1)^2(\overline{T} -\underline{T})^2 \overline{T}}{4N^2 \underline{T} T^\infty}, \\
\lambda_1 &= \min\left\{\frac{4}{(2\chi +1) \overline{T}^2}, ~~\frac{2}{\overline{T}}-\frac{1}{A_1}\left(\frac{ 2\| \mathbf{V}(0) \|^2 }{(2 \chi +1)^2N^2 \overline{T}^2}+2
+\frac{(N-1)^2(\overline{T} -\underline{T})^2}{2N^2 \underline{T} T^\infty}\right)\right\},
\end{aligned}
\end{align}
where 
\[ \underline{T}: = \frac{\prod_{a=1}^{N}T_a(0)}{\left[2E_\chi(0)/(2\chi+1)\right]^{N-1}}, \quad  \overline{T} := \frac{2E_\chi(0)}{2\chi+1}, \quad  T^\infty  = \Big( \chi + \frac{1}{2}  \Big)^{-1} \Big(  \frac{E_\chi(0)}{N} - \frac{1}{2N^2} |  \mathbf{M}(0)|^2  \Big). \]

\begin{lemma}\label{L4.2}
Let $\left\{\left(\mathbf{x}_a, \mathbf{v}_a, T_a\right)\right\}$ be a global solution to \eqref{TCSmodel-3}. 
For positive constants $A_1$ and $\lambda_1$ satisfying \eqref{D-2}, one has 
\begin{align}\label{diss-vT(TCS-3)}
\| \hat {\mathbf T}(t) \|^2 +A_1 \| {\mathbf V}(t) \| ^2 \le \Big(  \| \hat{\mathbf T}(0) \|^2 +A_1 \| \mathbf{V}(0) \|^2 \Big)e^{-\lambda_1t}, \quad t \geq 0.
\end{align}
\end{lemma}
\begin{proof}
It follows from \eqref{TCSmodel-3}$_2$ and \eqref{TCSmodel-3}$_3$ that 
\begin{align}\label{dT(TCS-3)}
\begin{aligned}    
	\frac{2 \chi +1}{2}\frac{\mathrm{d}T_a}{\mathrm{d}t} &= \frac{\mathrm{d}}{\mathrm{d}t}\left(\frac{2\chi +1}{2}{T_a}+\frac{1}{2}|{v_a}|^2\right)
	-\mathbf{v}_a\cdot\frac{\mathrm{d}\mathbf{v}_a}{\mathrm{d}t}  \\
	&= \frac{1}{N}\sum_{b=1}^{N}\left(\frac{1}{T_a}-\frac{1}{T_b}\right)
	-\frac{1}{N} \sum_{b=1}^{N} \mathbf{v}_a\cdot 
	\left(\frac{\mathbf{v}_b}{T_b}-\frac{\mathbf{v}_a}{T_a}\right).
\end{aligned}
\end{align}
Consider the linear combination of $|T_a-T^\infty|^2$ and $|\mathbf{v}_a|^2$:
\[ 
\sum_{a=1}^N\left[\frac{2\chi +1}{2}|T_a-T^\infty|^2+A_1|\mathbf{v}_a|^2\right] 
\]
and we choose $A_1$ such that the above quadratic quantity satisfies Gronwall's differential inequality. By direct calculation, one has
\begin{align}
\begin{aligned} \label{cT11}
	&\frac{\mathrm{d}}{\mathrm{d}t}\sum_{a=1}^N\left[\frac{2\chi +1}{2}|T_a-T^\infty|^2+A_1|\mathbf{v}_a|^2\right] \\
	& \hspace{0.5cm} = \frac{2}{N}\sum_{a,b=1}^{N}(T_a-T^\infty) \left(\frac{1}{T_a}-\frac{1}{T_b}\right)
	+\frac{2}{N}\sum_{a,b=1}^{N}
	(T^\infty-T_a+A_1) \mathbf{v}_a\cdot \left(\frac{\mathbf{v}_b}{T_b}-\frac{\mathbf{v}_a}{T_a}\right) \\
	& \hspace{0.5cm} =:\mathcal{I}_{31}+\mathcal{I}_{32}.
\end{aligned}
\end{align}
Below, we estimate the term ${\mathcal I}_{3i},~i = 1,2$ one by one. \newline   

\noindent $\bullet$ Case E.1 (Estimate of $\mathcal{I}_{31}$): Note that
\begin{align}
\begin{aligned} \label{D-3}
	&\sum_{a,b=1}^{N}(T_a-T^\infty) \left(\frac{1}{T_a}-\frac{1}{T_b}\right)=\sum_{a,b=1}^{N}(T_b-T^\infty) \left(\frac{1}{T_b}-\frac{1}{T_a}\right)
	\\
	& \hspace{1cm} =-\sum_{a,b=1}^{N}(T_a-T^\infty) \left(\frac{1}{T_b}-\frac{1}{T_a}\right)
	=\frac{1}{2}\sum_{a,b=1}^{N}(T_b-T_a) \left(\frac{1}{T_b}-\frac{1}{T_a}\right).
\end{aligned}
\end{align}
We use  \eqref{Ta-Tb(TCS)} and \eqref{D-3} to find 
\begin{align*}
\begin{aligned}
	& \mathcal{I}_{31} = \frac{1}{N}\sum_{a,b=1}^{N}(T_b-T_a) \left(\frac{1}{T_b}-\frac{1}{T_a}\right)
	= -\frac{1}{N}\sum_{a,b=1}^{N} \frac{|T_a-T_b|^2}{T_aT_b} \\
	& \hspace{0.5cm} \leq -\frac{1}{N \overline{T}^2}\sum_{a,b=1}^{N} |T_a-T_b|^2 = -\frac{2}{\overline{T}^2}\sum_{a=1}^N|T_a-T^\infty|^2+\frac{2}{(2 \chi+1)^2N^2 \overline{T}^2}\| {\mathbf V}(t) \| ^4, \\
	& \sum_{a,b=1}^{N} |T_a-T_b|^2 = 2N\sum_{a=1}^{N} \Big(|T_a-T^\infty|^2-|T^\infty- \langle T \rangle |^2 \Big).
\end{aligned} 
\end{align*} 

Then we use  Corollary \ref{C4.1} to obtain
\begin{align*}
\begin{aligned}
	\mathcal{I}_{31}\le -\frac{2}{\overline{T}^2}\sum_{a=1}^N|T_a-T^\infty|^2+\frac{2\| {\mathbf V}(0) \| ^2}{(2 \chi+1)^2N^2 \overline{T}^2}\| {\mathbf V}(t) \| ^2.
\end{aligned}
\end{align*}
\vspace{0.2cm}

\noindent $\bullet$ Case E.2 (Estimate of $\mathcal{I}_{32}$): Note that
$$|T_a- \langle T \rangle |=\Big|\frac{1}{N}\sum_{b=1}^{N}\left(T_a-T_b\right)\Big|\leq \frac{N-1}{N}(\overline{T} - \underline{T}). $$
Then, we can obtain
\begin{align*}
\mathcal{I}_{32}=& -\frac{2}{N}\sum_{a,b=1}^{N}(T_a-\overline{T})\frac{\mathbf{v}_a\cdot \mathbf{v}_b}{T_b}+2\sum_{a=1}^{N}|\mathbf{v}_a|^2-2\sum_{a=1}^N(A_1+T^\infty)\frac{|\mathbf{v}_a|^2}{T_a}\\
\leq& \frac{2}{N}\sum_{a,b=1}^{N}\left[T^\infty\frac{ |\mathbf{v}_b|^2}{T_b}+\frac{(T_a-\overline{T})^2}{4T^\infty}\frac{ |\mathbf{v}_b|^2}{\underline{T}}\right]+2\sum_{a=1}^{N}|\mathbf{v}_a|^2-2\sum_{a=1}^N(A_1+T^\infty)\frac{|\mathbf{v}_a|^2}{T_a}\\
\leq& \left(2+\frac{(N-1)^2( \overline{T} - \underline{T})^2}{2N^2 \underline{T} T^\infty}\right)\sum_{a=1}^N|\mathbf{v}_a|^2
-2A_1\sum_{a=1}^N\frac{|\mathbf{v}_a|^2}{\overline{T}}.
\end{align*}
We collect all the estimates of $ \mathcal{I}_{11}$ and  $\mathcal{I}_{12}$  in \eqref{cT11}, and choose the constant $A_1$ such that
\begin{align*}
A_1>\frac{ \| \mathbf{V}(0) \|^2 }{(2 \chi +1)^2N^2 \overline{T}}+ \overline{T}
+\frac{(N-1)^2(\overline{T} -\underline{T})^2 \overline{T}}{4N^2 \underline{T} T^\infty}
\end{align*}
to obtain 
\begin{align*}
&\frac{\mathrm{d}}{\mathrm{d}t}\Big( \| \hat {\mathbf T}(t) \|^2 +A_1 \| {\mathbf V}(t) \|^2 \Big)+\frac{4\| \hat {\mathbf T}(t) \|^2}{(2\chi +1) \overline{T}^2}\\
&\qquad-\left(\frac{ 2\| \mathbf{V}(0) \|^2 }{(2 \chi +1)^2N^2 \overline{T}^2}+2
+\frac{(N-1)^2(\overline{T} - \underline{T})^2}{2N^2 \underline{T} T^\infty}-\frac{2A_1}{\overline{T}}\right)  \| {\mathbf V}(t) \|^2 \leq  0.
\end{align*}
This implies \eqref{diss-vT(TCS-3)}.
\end{proof}

\vspace{0.2cm}

\noindent Next, we are ready to present our first result as follows. 
\begin{theorem}\label{T4.1}
Let $\left\{\left(\mathbf{x}_a, \mathbf{v}_a, T_a\right)\right\}$ be a global solution to \eqref{TCSmodel-3}. 
Then, for positive constants $A_1$ and $\lambda_1$ in \eqref{D-2}, we have
\begin{align*}
\begin{aligned}
	& \| \mathbf{X}(t) \| \leq \| \mathbf{X}(0) \|  + \| {\mathbf V}(0) \| \overline{T}, \quad \| \mathbf{V}(t) \|  \le  \| \mathbf{V}(0) \|  e^{-\frac{t}{\overline{T}}}, \\
	& \| \hat {\mathbf T}(t) \|^2 +A_1 \| {\mathbf V}(t) \| ^2 \le \Big(  \| \hat{\mathbf T}(0) \|^2 +A_1 \| \mathbf{V}(0) \|^2 \Big)e^{-\lambda_1t},
\end{aligned}
\end{align*}
for $t \geq 0$. 
\end{theorem}
\begin{proof}
We combine Lemma \ref{L4.1} and Lemma \ref{L4.2} to derive the desired estimates. 
\end{proof}

\subsection{Flocking dynamics II} \label{sec:4.2}
In this subsection, we present the asymptotic flocking dynamics of  \eqref{TCSmodel-2} in which the communication weight functions $\phi_{ab}$ satisfy \eqref{phiMm} for $a, b\in[N].$ In what follows, we show that for $\varepsilon$ sufficiently small, flocking estimates hold without additional conditions. \newline

In the next lemma, we derive a system of dissipative inequalities for position and velocity.
\begin{lemma}\label{L4.3}
Let $\left\{\left(\mathbf{x}_a, \mathbf{v}_a, T_a\right)\right\}$ be a global solution to \eqref{TCSmodel-2} with zero initial momentum:
\[ \langle \mathbf{v}(0) \rangle= 0. \]
Then, we have a system of dissipative inequalities: 
\begin{equation} \label{D-4}
\frac{\mathrm{d} \| \mathbf{X}\|}{\mathrm{d}t}  \le \| \mathbf{V} \|, \qquad
\frac{\mathrm{d} \| \mathbf{V} \|}{\mathrm{d}t}  \le -\left(\frac{\underline{\phi}}{\overline{T}} -\frac{\varepsilon}{\underline{T}}\right) \| \mathbf{V} \|, \quad \mbox{a.e.}~t > 0.
\end{equation}
\end{lemma}
\begin{proof}
Since the derivation of the first differential inequality in \eqref{D-4} is the same as in Lemma \ref{L4.1}. Thus, we focus on the second differential inequality. First, we take an inner product of $\mathbf{v}_a$ with \eqref{TCSmodel-2}$_2$ and take sum the resulting relation over  $a\in[N]$ to get the desired estimate:
\begin{align*} 
\begin{aligned} \label{D-5}
	\frac{\mathrm{d}}{\mathrm{d}t}\left(\sum_{a=1}^{N}|\mathbf{v}_a|^2\right) &=\frac{2 \underline{\phi}}{N}\sum_{a,b=1}^{N}\left(\frac{\mathbf{v}_a\cdot\mathbf{v}_b}{T_b}-\frac{|\mathbf{v}_a|^2}{T_a}\right) +\frac{2}{N}\sum_{a,b=1}^{N}\left(\phi_{ab} - \underline{\phi} \right)\left(\frac{\mathbf{v}_a\cdot\mathbf{v}_b}{T_b}-\frac{|\mathbf{v}_a|^2}{T_a}\right)
	\\ 
	&=: {\mathcal I}_{41} + {\mathcal I}_{42}.
\end{aligned}
\end{align*}
Below, we estimate the terms ${\mathcal I}_{4i}$ one by one. \newline

\noindent $\bullet$~(Estimate of ${\mathcal I}_{41}$): We use  
\begin{align*}
\sum_{a,b=1}^{N} \frac{\mathbf{v}_a\cdot\mathbf{v}_b}{T_a}  =  \Big( \sum_{a=1}^{N} \frac{\mathbf{v}_a}{T_a} \Big) \cdot \Big( \sum_{b=1}^{N} \mathbf{v}_b \Big) = 0 \quad  \mbox{and} \quad 
\sum_{a,b=1}^{N} \frac{|\mathbf{v}_a|^2}{T_a}  \geq \frac{1}{\overline{T}} \sum_{a,b=1}^{N} |\mathbf{v}_a|^2 =  \frac{N}{\overline{T}}  \sum_{a=1}^{N} |\mathbf{v}_a|^2
\end{align*}
to find 
\begin{equation} \label{D-6}
{\mathcal I}_{41} = \frac{2 \underline{\phi}}{N}\sum_{a,b=1}^{N}\left(\frac{\mathbf{v}_a\cdot\mathbf{v}_b}{T_b}-\frac{|\mathbf{v}_a|^2}{T_a}\right) = \frac{2 \underline{\phi}}{N}\sum_{a,b=1}^{N}\frac{\mathbf{v}_a\cdot\mathbf{v}_b}{T_b}-  \frac{2 \underline{\phi}}{N}\sum_{a,b=1}^{N} \frac{|\mathbf{v}_a|^2}{T_a} \leq  -\frac{2 \underline{\phi} }{ \overline{T}}\sum_{a=1}^{N} |\mathbf{v}_a|^2.
\end{equation}
\noindent $\bullet$~(Estimate of ${\mathcal I}_{42}$): We use 
\begin{align*}
\begin{aligned}
	& \Big| \frac{2}{N}\sum_{a,b=1}^{N}\left(\phi_{ab} - \underline{\phi} \right) \frac{\mathbf{v}_a\cdot\mathbf{v}_b}{T_b} \Big| \leq  \frac{2 \varepsilon}{N \underline{T}}  \sum_{a,b=1}^{N} |\mathbf{v}_a| |\mathbf{v}_b| \\ 
	&\hspace{1.5cm}\leq \frac{2 \varepsilon}{N \underline{T}}  \Big( \sum_{a,b=1}^{N} |\mathbf{v}_a|^2 \Big)^{\frac{1}{2}}   \Big( \sum_{a,b=1}^{N} |\mathbf{v}_b|^2 \Big)^{\frac{1}{2}} 
	\leq \frac{2\varepsilon}{\underline{T}}\sum_{a=1}^{N} |\mathbf{v}_a|^2 
\end{aligned}
\end{align*}
to get 

\begin{equation} \label{D-7}
{\mathcal I}_{42}\leq \frac{2}{N}\sum_{a,b=1}^{N}\left(\phi_{ab} - \underline{\phi} \right) \frac{\mathbf{v}_a\cdot\mathbf{v}_b}{T_b} \leq \frac{2\varepsilon}{\underline{T}}\sum_{a=1}^{N} |\mathbf{v}_a|^2.
\end{equation}
We combine \eqref{D-6} and \eqref{D-7} to get the desired estimate. 

\end{proof}
As a direct application of Lemma \ref{L4.3}, we have the following estimates. 
\begin{corollary}
\label{C4.2}
Suppose the communication weight functions and $\varepsilon$ satisfy 
\begin{equation} \label{D-8}
\underline{\phi} > 0, \quad \underline{T}  > 0, \quad   0<  \varepsilon \leq \frac{ \underline{\phi} \underline{T}}{2 \overline{T}},
\end{equation}
and let $\left\{\left(\mathbf{x}_a, \mathbf{v}_a, T_a\right)\right\}$ be a global solution to \eqref{TCSmodel-2}. Then, one has
\begin{align*}
\| \mathbf{X}(t) \| \le \| \mathbf{X}(0) \|+\frac{2\overline{T} }{ \underline{\phi}} \| \mathbf{V}(0)\|, \qquad \| \mathbf{V}(t) \| \le  \| \mathbf{V}(0) \| e^{-\frac{\underline{\phi}}{2 \overline{T}}t}, \quad t \geq 0.
\end{align*}
\end{corollary} 
\begin{proof} By $\eqref{D-8}_3$, one has 
\[
-\left(\frac{2\underline{\phi}}{\overline{T}} -\frac{2\varepsilon}{\underline{T}}\right) <  - \frac{\underline{\phi}}{ \overline{T}}.
\]
Thus, we have
\[   \frac{d \| \mathbf{V} \|}{\mathrm{d}t}  \le  - \frac{\underline{\phi}}{2\overline{T}} \| \mathbf{V} \| \quad \mbox{a.e.}~t > 0, \quad  \mbox{i.e.,} \quad \| \mathbf{V}(t) \| \leq   \| \mathbf{V}(0) \| e^{-\frac{\underline{\phi}}{2 \overline{T}}t}, \quad t \geq 0.
\]
Again, we integrate the estimate for $\| {\mathbf V} \|$ to find 
\[
\| \mathbf{X}(t) \| \leq  \| \mathbf{X}(0) \| +\frac{2\|\overline{T} }{ \underline{\phi}} \| \mathbf{V}(0) \Big(1 -   e^{-\frac{\underline{\phi}}{ \overline{T}}t} \Big) \leq  \| \mathbf{X}(0) \| +\frac{2\overline{T} }{ \underline{\phi}} \| \mathbf{V}(0) \|.
\]
\end{proof}

\vspace{0.2cm}

Next, we derive exponential time-decay estimates for temperature and velocity configurations. We choose positive constants  $A_2$ and $\lambda_2$ such that 
\begin{align}
\label{AL2}
A_2 &> \frac{4}{7}\left[\frac{2 \underline{\zeta} \| {\mathbf V}(0) \| ^2}{(2\chi +1)^2 N^2\underline{\phi} \overline{T} }+\left(2
+\frac{(N-1)^2( \overline{T} - \underline{T})^2}{2N^2 \underline{T} T^\infty}\right) \overline{T} +\left(2+\frac{\| {\mathbf V}(0) \| ^2}{2(2\chi +1) \underline{T}}\right)\frac{\varepsilon \overline{T}}{ \underline{\phi} }\right],\nonumber \\
\lambda_2 &:=\min\left\{\frac{4 \underline{\zeta}}{(2 \chi +1) \overline{T}^2},~\frac{7 \underline{\phi}}{4 \overline{T}}-\frac{1}{A_2}\Big[ \frac{2 \underline{\zeta} \| {\mathbf V}(0) \| ^2}{(2\chi +1)^2 N^2\overline{T}^2}
+\left(2+\frac{ \| {\mathbf V}(0) \| ^2}{2(2\chi +1) \underline{T} }\right)\varepsilon \right.\\
&\hspace{1.5cm}\left.+\left( 2+\frac{(N-1)^2(\overline{T} - \underline{T})^2}{2N^2 \underline{T} T^\infty}\right) \underline{\phi}  \Big ]\right\}.\nonumber
\end{align}
\begin{lemma}\label{L4.4}
Let $\left\{\left(\mathbf{x}_a, \mathbf{v}_a, T_a\right)\right\}$ be a global solution to \eqref{TCSmodel-2}. For positive constants $A_2$ and $\lambda_2$ satisfying \eqref{AL2}, one has 
\begin{align*}
\| \hat {\mathbf T}(t) \|^2 +A_2 \| {\mathbf V}(t) \| ^2  \le \Big(  \| \hat {\mathbf T}(0) \|^2 +A_2 \| {\mathbf V}(0) \| ^2      \Big)e^{-\lambda_2t}, \quad t \geq 0.
\end{align*}
\end{lemma}
\begin{proof}
Similar to \eqref{dT(TCS-3)}, it follows from $\eqref{TCSmodel-2}_2$ and  $\eqref{TCSmodel-2}_3$ that 
\[
\frac{2 \chi +1}{2} \frac{\mathrm{d}T_a}{\mathrm{d}t}
=\frac{1}{N}\sum_{b=1}^{N}\zeta_{ab}\left(\frac{1}{T_a}-\frac{1}{T_b}\right)
-\frac{\underline{\phi}}{N} \sum_{b=1}^{N} \mathbf{v}_a \cdot
\left(\frac{\mathbf{v}_b}{T_b}-\frac{\mathbf{v}_a}{T_a}\right) -\frac{1}{N}  \sum_{b=1}^{N}\left(\phi_{ab}- \underline{\phi} \right) \mathbf{v}_a\cdot \left(\frac{\mathbf{v}_b}{T_b}-\frac{\mathbf{v}_a}{T_a}\right). 
\]
For the positive constant $A_2$  in \eqref{AL2}, we can obtain
\begin{align}\label{cT12}
&\frac{\mathrm{d}}{\mathrm{d}t}\sum_{a=1}^N\left( \frac{2 \chi +1}{2}|T_a-T^\infty|^2+A_2|\mathbf{v}_a|^2\right) \nonumber\\
& \hspace{0.5cm} =\frac{2}{N}\sum_{a,b=1}^{N}\zeta_{ab}(T_a-T^\infty) \left(\frac{1}{T_a}-\frac{1}{T_b}\right)  +\frac{2 \phi_{ab}}{N}\sum_{a,b=1}^{N}
(T^\infty-T_a) \mathbf{v}_a\cdot  \left(\frac{\mathbf{v}_b}{T_b}-\frac{\mathbf{v}_a}{T_a}\right) 
\nonumber\\
& \hspace{0.5cm} +\frac{2\phi_{ab}}{N}\sum_{a,b=1}^{N} A_2 \mathbf{v}_a\cdot \left(\frac{\mathbf{v}_b}{T_b}-\frac{\mathbf{v}_a}{T_a}\right)
\\
& \hspace{0.5cm} =\frac{2}{N}\sum_{a,b=1}^{N}\zeta_{ab}(T_a-T^\infty) \left(\frac{1}{T_a}-\frac{1}{T_b}\right)  +\frac{2 \underline{\phi}}{N}\sum_{a,b=1}^{N}
(T^\infty-T_a+A_2) \mathbf{v}_a\cdot  \left(\frac{\mathbf{v}_b}{T_b}-\frac{\mathbf{v}_a}{T_a}\right) 
\nonumber\\
& \hspace{0.5cm} +\frac{2}{N}\sum_{a,b=1}^{N}\left(\phi_{ab}- \underline{\phi} \right)(T^\infty-T_a+A_2)   \mathbf{v}_a\cdot \left(\frac{\mathbf{v}_b}{T_b}-\frac{\mathbf{v}_a}{T_a}\right)
\nonumber\\
& \hspace{0.5cm} =:\mathcal{I}_{41}+\mathcal{I}_{42}+\mathcal{I}_{43}. \nonumber 
\end{align}
\vspace{0.2cm}

In what follows, we estimate the term ${\mathcal I}_{4i},~i = 1,2,3$ one by one. \newline

\noindent $\bullet$~Case F.1 (Estimate of $\mathcal{I}_{41}$): We use the index exchange transform $ (a, b)~ \leftrightarrow~(b,a)$ and $\zeta_{ba} = \zeta_{ab}$ to find 
\begin{align}
\begin{aligned} \label{D-9}
	&\sum_{a,b=1}^{N}\zeta_{ab}(T_a-T_\infty) \left(\frac{1}{T_a}-\frac{1}{T_b}\right)   \\
	& \hspace{0.5cm} =  \sum_{a,b=1}^{N}\zeta_{ba}(T_b-T_\infty) \left(\frac{1}{T_b}-\frac{1}{T_a}\right) =  -\sum_{a,b=1}^{N}\zeta_{ab}(T_b-T_\infty) \left(\frac{1}{T_a}-\frac{1}{T_b}\right)  \\
	& \hspace{0.5cm} =\frac{1}{2}\sum_{a,b=1}^{N}\zeta_{ab}(T_a-T_b) \left(\frac{1}{T_a}-\frac{1}{T_b}\right) = -\frac{1}{2}\sum_{a,b=1}^{N}\zeta_{ab} \frac{|T_a-T_b|^2} {T_aT_b}.
\end{aligned}
\end{align}
Then, we use \eqref{D-9}, the relations:
\[
\zeta_{ab} \geq \underline{\zeta} \quad \mbox{and} \quad  T_a T_b \leq \overline{T}^2 
\]
to find 
\begin{align*}
\mathcal{I}_{41} \leq -\frac{\underline{\zeta}}{N \overline{T}^2}\sum_{a,b=1}^{N} |T_a-T_b |^2      \leq -\frac{2 \underline{\zeta}}{\overline{T}^2}\sum_{a=1}^{N} |T_a-T^\infty|^2+\frac{2 \underline{\zeta} \|{\mathbf V}(0) \|^2}{(2\chi+1)^2N^2 \overline{T}^2} \|{\mathbf V}(t) \|^2.
\end{align*}

\noindent $\bullet$~Case F.2  (Estimation of $\mathcal{I}_{42}$): By the same computation as in the Case E.2, one has 
\begin{align*}
\begin{aligned}
	\mathcal{I}_{42}  &\leq  \left(2+\frac{(N-1)^2(\overline{T} - \underline{T})^2}{2N^2 \underline{T} T^\infty}\right) \underline{\phi} \sum_{a=1}^{N}|\mathbf{v}_a|^2
	-2A_2 \underline{\phi} \sum_{a=1}^N\frac{|\mathbf{v}_a|^2}{\overline{T}}.
\end{aligned}
\end{align*}

\noindent $\bullet$~Case F.3 (Estimation of $\mathcal{I}_{43}$): According to Lemma \ref{C4.2}, we use
\[
\frac{1}{N}\sum_{a=1}^{N}\frac{(T^{\infty}-T_a)^2}{T^\infty} =T^{\infty}- \frac{1}{N}\sum_{a=1}^{N}T_a- \frac{1}{N}\sum_{a=1}^{N}\left(T_a-\frac{T_a^2}{T^\infty}\right)\leq  \frac{\|{\mathbf V}(t) \|^2}{(2\chi+1)N} \| \leq  \frac{\|{\mathbf V}(0) \|^2}{(2\chi+1)N}
\]
to find 
\begin{align*}
\mathcal{I}_{43} &= \frac{2}{N}\sum_{a,b=1}^{N}(T^\infty-T_a)\left(\phi_{ab}- \underline{\phi} \right)\frac{\mathbf{v}_a\cdot \mathbf{v}_b}{T_b}+2\sum_{a=1}^{N}\left(\phi_{ab}- \underline{\phi} \right)|\mathbf{v}_a|^2
\\
&+\frac{2}{N}\sum_{a,b=1}^{N} A_2\left(\phi_{ab} - \underline{\phi} \right)\frac{\mathbf{v}_a\cdot\mathbf{v}_b}{T_a}-\frac{2}{N}\sum_{a,b=1}^N
\left(\phi_{ab}-\overline{\phi}\right)(A_2+T_\infty)\frac{|\mathbf{v}_a|^2}{T_a}
\\
& \le \frac{2}{N}\sum_{a,b=1}^{N}(\phi_{ab}- \underline{\phi} )\left[T^\infty\frac{ |\mathbf{v}_b|^2}{T_b}+\frac{(T^{\infty}-T_a)^2}{4T^\infty}\frac{ |\mathbf{v}_b|^2}{\underline{T} }\right]+2\varepsilon\sum_{a=1}^{N}|\mathbf{v}_a|^2
\\
&+\frac{2A_2}{N}\sum_{a,b=1}^N(\phi_{ab}- \underline{\phi})\left(\frac{|\mathbf{v}_a|^2}{T_a}+\frac{|\mathbf{v}_b|^2}{4\underline{T} }\right)
-\frac{2}{N}\sum_{a,b=1}^N\left(\phi_{ab}- \underline{\phi} \right)(A_2+T_\infty)\frac{|\mathbf{v}_a|^2}{T_a}
\\ 
& \leq \varepsilon\left(2+\frac{\|{\mathbf V}(0) \|^2}{2(2\chi+1) \underline{T} }+\frac{A_2}{2\underline{T}}\right)\sum_{a=1}^{N}|\mathbf{v}_a|^2.
\end{align*}
In \eqref{cT12}, we collect all the estimates in Case F.1 - Case F.3 to obtain the desired estimate:
\begin{align*}
&\frac{\mathrm{d}}{\mathrm{d}t}\Big(  \| \hat {\mathbf T} \|^2 +A_2 \| {\mathbf V} \|^2              \Big)
+\frac{4 \underline{\zeta}}{(2\chi +1) \overline{T}^2} \| \hat{\mathbf T} \|^2  \\
& \hspace{1cm} -\Big[ \Big (2+\frac{(N-1)^2( \overline{T}  - \underline{T})^2}{2N^2 \underline{T} T^\infty} \Big) \underline{\phi}  +\frac{2\underline{\zeta}}{(2\chi +1)^2 N^2\overline{T}^2}  \| {\mathbf V}(0) \|^2  \\
&  \hspace{1cm} +\left(2+\frac{\|{\mathbf V}(0) \|^2}{2(2\chi+1) \underline{T} }\right) \varepsilon -A_2\Big(\frac{2 \underline{\phi}}{\overline{T}}-\frac{\varepsilon}{2 \underline{T}} \Big)\Big]  \| {\mathbf V}(t) \|^2  \leq  0.
\end{align*}
\end{proof}
Next, we provide our second main result by combining the result of Lemma \ref{L4.3}, Corollary \ref{C4.2}, and Lemma \ref{L4.4} to derive the asymptotic flocking dynamics. 
\begin{theorem}\label{T4.2}
Let $\left\{\left(\mathbf{x}_a, \mathbf{v}_a, T_a\right)\right\}$ be a global solution to \eqref{TCSmodel-2}. For positive constants $A_2, \lambda_2$ in \eqref{AL2}, we have
\begin{align*}
\begin{aligned}
	&  \| \mathbf{X}\| \le \| \mathbf{X}(0) \|+\frac{2\overline{T} }{ \underline{\phi}} \| \mathbf{V}(0)\|, \quad \| \mathbf{V}(t) \| \le  \| \mathbf{V}(0) \| e^{-\frac{\underline{\phi}}{ 2\overline{T}}t}, \\
	&  \| \hat {\mathbf T}(t) \|^2 +A_2 \| {\mathbf V}(t) \| ^2  \le \Big(  \| \hat {\mathbf T}(0) \|^2 +A_2 \| {\mathbf V}(0) \| ^2      \Big)e^{-\lambda_2t}, \quad t \geq 0.
\end{aligned}
\end{align*}
\end{theorem}
\subsection{Flocking dynamics III} \label{sec:4.3}
In this subsection,  we show that for well-prepared initial data, flocking estimates hold for the TCS model \eqref{A-3}. \newline

It follows from \eqref{properties of phi zeta} and the relation
\begin{align*}
|\mathbf{x}_a-\mathbf{x}_b|\le \sqrt{2(|\mathbf{x}_a|^2+|\mathbf{x}_b|^2)}\le \sqrt{2}\| \mathbf{X}\|,\quad a, b \in [N],
\end{align*}
that 
\begin{align*}
\phi( \sqrt{2}\| \mathbf{X} \| )\le\phi_{ab}\le\phi(0),\quad \zeta(\sqrt{2}\| \mathbf{X} \|)\le\zeta_{ab}\le\zeta(0).
\end{align*}
Next, we derive a system of dissipative inequalities for position and velocity.
\begin{lemma}\label{L4.5}
Let $\left\{\left(\mathbf{x}_a, \mathbf{v}_a, T_a\right)\right\}$ be a global solution to \eqref{A-3}. Then, one has the following estimates:
\begin{align*}
\frac{\mathrm{d} \| \mathbf{X} \|}{\mathrm{d}t}\le  \| \mathbf{V} \|, \quad \frac{\mathrm{d}\| \mathbf{V} \|}{\mathrm{d}t}\le -\frac{\underline{T} \phi^2(0)- \overline{T} (\phi(0)-\phi( \| \mathbf{X} \|))^2}{2 \overline{T} \underline{T} \phi(0)}  \| \mathbf{V} \|, \quad \mbox{a.e.}~t > 0.
\end{align*}
\end{lemma}
\begin{proof}
The first estimate is the same as in Lemma \ref{L4.1}. Thus, we focus on the derivation of the second differential inequality. We take an inner product of $\mathbf{v}_a$ with \eqref{A-3}$_2$, and sum up the resulting relation over  $a\in[N]$ to get
\begin{align} 
\begin{aligned} \label{dv(TCS 1)}
	\frac{\mathrm{d}}{\mathrm{d}t}\left(\sum_{a=1}^{N}|\mathbf{v}_a|^2\right) &=-\frac{2}{N}\sum_{a,b=1}^{N}\phi_{ab}\frac{|\mathbf{v}_a|^2}{T_a}
	+\frac{2}{N}\sum_{a,b=1}^{N}\left(\phi_{ab}- \underline{\phi} \right)
	\left(\frac{\mathbf{v}_a\cdot\mathbf{v}_b}{T_b}\right)
	\\  
	&\le -\frac{2}{N}\sum_{a,b=1}^{N}\phi_{ab}\frac{|\mathbf{v}_a|^2}{T_a}
	+\frac{1}{N}\sum_{a,b=1}^{N}\phi_{ab}\frac{|\mathbf{v}_b|^2}{T_b}
	+\frac{1}{N}\sum_{a,b=1}^{N}\frac{(\phi_{ab}- \underline{\phi} )^2|\mathbf{v}_a|^2}{\phi_{ab}T_b}
	\\
	&\le -\frac{1}{N}\sum_{a,b=1}^{N}\left[\frac{\phi_{ab}}{T_a}-
	\frac{(\phi_{ab}- \underline{\phi})^2}{\phi_{ab} \underline{T}}\right]|\mathbf{v}_a|^2 \\
	& \leq-\left(\frac{\phi(0)}{\overline{T}}
	-\frac{(\phi(0)- \underline{\phi})^2}{\phi(0) \underline{T}}\right)\sum_{a=1}^{N}|\mathbf{v}_a|^2.
\end{aligned}
\end{align}
This yields the desired differential inequality.   
\end{proof}
\noindent To estimate the position and the velocity, we define 
\begin{align*}
\chi_0(\sqrt{2} \| \mathbf{X} \|) :=\frac{\underline{T} \phi^2(0)- \overline{T} (\phi(0)-\phi(\sqrt{2} \| \mathbf{X} \|))^2}{2 \overline{T} \underline{T} \phi(0)}.
\end{align*}
As a direct application of Lemma \ref{L4.5}, we have the following estimate. 
\begin{corollary} \label{C4.3}
For given initial data, we assume that there exists a positive constant $U_0$ such that
\begin{align}\label{phiX phi0}
\phi(U_0)>\left(1-\sqrt{\frac{\underline{T}}{\overline{T}}}\right)\phi(0),\quad
U_0\le \sqrt{2}\| \mathbf{X}(0) \| +\frac{\sqrt{2}}{\chi_0(U_0)} \| \mathbf{V}(0) \|,
\end{align}
and let $\left\{\left(\mathbf{x}_a, \mathbf{v}_a, T_a\right)\right\}$ be a global solution to \eqref{A-3}. Then, one has 
\begin{equation}\label{XU3v-bound}
\sup_{0 \leq t < \infty} \sqrt{2}\| \mathbf{X}(t) \| \le U_0,\quad  \| \mathbf{V} (t) \| \le  \| \mathbf{V}(0) \| e^{-\chi_0(U_0) t}, \quad t \geq 0.
\end{equation}
\end{corollary}
\begin{proof}
Note that $\chi_0(U_0)>0$ under the assumption  \eqref{phiX phi0}.
Then,  we have \eqref{XU3v-bound} from Lemma \ref{L4.5}.
\end{proof}

\vspace{0.2cm}

Next, we derive exponential time-decay estimates for temperature and velocity. For this, we set positive constants $A_3$ and $\lambda_3$ such that 
\begin{align}
\begin{aligned} \label{AL3}
A_3 &> \frac{1}{2\chi_0(U_0)}\Big\{  \frac{2\| \mathbf{V}(0) \|^2 \zeta\left(U_0\right)}{(2\chi +1)^2N^2 \overline{T}^2}
+\left(2+\frac{\| \mathbf{V}(0) \|^2}{2(2\chi +1) \underline{T}}\right)\phi(0)  \Big \}, 
\\\lambda_3 :=&\min\left\{\frac{4\zeta\left(U_0\right)}{(2 \chi +1) \overline{T}^2},
2\chi_0\left(U_0\right)-\frac{1}{A_3}\left(\frac{2 \| \mathbf{V}(0) \|^2 \zeta\left(U_0\right)}{(2 \chi +1)^2N^2 \overline{T}^2}+\left(2
+\frac{\| \mathbf{V}(0) \|^2}{2(2 \chi +1)\underline{T}}\right)\phi(0)\right)\right\}.
\end{aligned}
\end{align}
\begin{lemma}\label{L4.6}
Suppose the conditions \eqref{phiX phi0} hold, and  let $\{(\mathbf{x}_a, \mathbf{v}_a, T_a)\}$ be a global solution to \eqref{A-3}. Then one has 
\begin{equation*}\label{diss-vT(TCS-1)}
\| \hat {\mathbf T}(t) \|^2 +A_3 \| {\mathbf V}(t) \|^2  \le \Big(  \| \hat {\mathbf T}(0) \|^2 +A_3 \| {\mathbf V}(0) \| ^2      \Big)e^{-\lambda_3 t}, \quad t \geq 0.
\end{equation*}
\end{lemma}
\begin{proof}
Similar to \eqref{dT(TCS-3)}, it follows from \eqref{A-3} that 
\begin{align*}
\frac{2\chi +1}{2}\frac{\mathrm{d}T_a}{\mathrm{d}t}
=&\frac{1}{N}\sum_{b=1}^{N}\zeta_{ab}\left(\frac{1}{T_a}-\frac{1}{T_b}\right)-\frac{1}{N}\mathbf{v}_a\cdot\sum_{b=1}^{N}\phi_{ab}\left(\frac{\mathbf{v}_b}{T_b}-\frac{\mathbf{v}_a}{T_a}\right) .
\end{align*}
For the positive constant $A_3$ given in \eqref{AL3}, we obtain
\begin{align}\label{cT13}
\begin{aligned}
	&\frac{\mathrm{d}}{\mathrm{d}t}\sum_{a=1}^N
	\left[\frac{2\chi +1}{2}|T_a-T^\infty|^2+A_3|\mathbf{v}_a|^2\right]\\
	&=\frac{2}{N}\sum_{a,b=1}^{N}\zeta_{ab}(T_a-T^\infty) \left(\frac{1}{T_a}-\frac{1}{T_b}\right)
	+ \frac{2}{N}\sum_{a,b=1}^{N}\mathbf{v}_a\cdot\phi_{ab}(T^\infty-T_a+A_3)\left(\frac{\mathbf{v}_b}{T_b}-\frac{\mathbf{v}_a}{T_a}\right)
	\\
	&=:\mathcal{I}_{51}+\mathcal{I}_{52}.
\end{aligned}
\end{align}
In what follows, we estimate the terms $\mathcal{I}_{5i}$ one by one. \newline

\noindent $\bullet$~Case G.1 (Estimate of $\mathcal{I}_{51}$): By direct calculation, we have
\begin{align*}
\mathcal{I}_{51} &= -\frac{1}{N}\sum_{a,b=1}^{N} \zeta_{ab}\frac{(T_a-T_b)^2}{T_aT_b}
\leq -\frac{\zeta\left(U_0\right)}{N \overline{T}^2}\sum_{a,b=1}^{N} (T_a-T_b)^2
\\ 
&=  -\frac{4\zeta\left(U_0\right)}{(2 \chi+1)\overline{T}^2}   \| \hat {\mathbf T}(t) \|^2 +\frac{2\zeta\left(U_0\right)\| {\mathbf V}(0) \|^2}{(2 \chi+1)^2N^2 \overline{T}^2} \| {\mathbf V} (t)\|^2.
\end{align*}

\noindent $\bullet$~Case G.2 (Estimate of $\mathcal{I}_{52}$):~We use $\langle \mathbf{v}(t)\rangle=0$ to obtain
\begin{align*}
\mathcal{I}_{52} \le& \frac{2}{N}\sum_{a,b=1}^{N}\phi_{ab}\left[T^\infty\frac{ |\mathbf{v}_b|^2}{T_b}+\frac{(T^\infty-T_a)^2}{4T^\infty}\frac{ |\mathbf{v}_b|^2}{\underline{T}}\right]+\frac{2}{N}\sum_{a,b=1}^{N}\phi_{ab}|\mathbf{v}_a|^2
\\ &+\frac{2A_3}{N}\sum_{a,b=1}^N\left(\phi_{ab}\frac{|\mathbf{v}_a|^2}{2T_a}+\frac{|\phi_{ab}-\phi\left(U_0\right)|^2}{\phi_{ab}}\frac{|\mathbf{v}_b|^2}{2T_a}\right)-\frac{2(A_3+T^\infty)}{N}\sum_{a,b=1}^N\phi_{ab}\frac{|\mathbf{v}_a|^2}{T_a}
\\ \le& \left(2+\frac{\| \mathbf{V}(0) \|^2}{2(2 \chi +1)\underline{T}}\right) \overline{\phi} \sum_{a=1}^{N}|\mathbf{v}_a|^2
-\frac{A_3}{N}\sum_{a,b=1}^{N}\left(\frac{\phi_{ab}}{T_a}
-\frac{(\phi_{ab}-\phi\left(U_0\right))^2}{\phi_{ab}T_b}\right)|\mathbf{v}_a|^2.
\end{align*}
In \eqref{cT13}, we collect all the estimates in Case G.1 and Case G.2 to obtain  the desired estimate:
\begin{align*}
\begin{aligned}
	& \frac{\mathrm{d}}{\mathrm{d}t}  \Big (\| \hat {\mathbf T}(t) \|^2 +A_3 \| {\mathbf V}(t) \|^2 \Big ) + \frac{4\zeta\left(U_0\right) \| \hat {\mathbf T}(t) \|^2}{(2\chi +1) \overline{T}^2}\\
	&\qquad-\Big[\frac{2 \| \mathbf{V}(0) \|^2 \zeta\left(U_0\right)}{(2 \chi +1)^2N^2 \overline{T}^2}+\left(2
	+\frac{\| \mathbf{V}(0) \|^2}{2(2 \chi +1)\underline{T}}\right)\phi(0)-2A_3\chi_0\left(U_0\right)\Big]  \| {\mathbf V}(t) \|^2
	\le  0.
\end{aligned}
\end{align*}
\end{proof}
Next, we verify the emergence of asymptotic flocking by the estimation of position, velocity, and temperature under the well-prepared condition using Corollary \ref{C4.3} and Lemma \ref{L4.6}.
\begin{theorem}\label{T4.3}
Suppose the conditions \eqref{phiX phi0} hold, and  let $\{(\mathbf{x}_a, \mathbf{v}_a, T_a)\}$ be a global solution to \eqref{A-3}. For positive constants $A_3, \lambda_3$ in \eqref{AL3}, we have
\begin{align*}
\begin{aligned}
	& \sup_{0 \leq t < \infty} \sqrt{2}\| \mathbf{X}(t) \| \le U_0, \quad \| \mathbf{V}(t) \| \le  \| \mathbf{V}(0) \| e^{-\chi_0(U_0) t}, \\
	&  \| \hat {\mathbf T}(t) \|^2 +A_3 \| {\mathbf V}(t) \|^2  \le \Big(  \| \hat {\mathbf T}(0) \|^2 +A_3 \| {\mathbf V}(0) \|^2 \Big)e^{-\lambda_3 t}, \quad t \geq 0.
\end{aligned}
\end{align*}
\end{theorem}

%
%
\section{Asymptotic flocking of the RTCS model with Synge energy} \label{sec:5} 
\setcounter{equation}{0}
In this section, we study the asymptotic flocking dynamics of the RTCS model \eqref{RTCSmodel(same num)-1} with different types of communication weight functions. More explicitly, we study the flocking estimation of \eqref{RTCSmodel(same num)-3} with $\phi_{ab}=\zeta_{ab}\equiv1$,  \eqref{RTCSmodel(same num)-2} with communication weight functions satisfying \eqref{phiMm}, and the RTCS model \eqref{RTCSmodel(same num)-1} with general communication functions.

\subsection{Flocking dynamics I} \label{sec:5.1}
In this subsection, we study the asymptotic flocking dynamics of the simplified RTCS model \eqref{RTCSmodel(same num)-3}. We will show that asymptotic flocking estimates hold without any condition on initial data. For convenience, we set
\[ \mathbf{W} := (\mathbf{w}_1, \cdots, \mathbf{w}_N), \quad  \| \mathbf{W} \|:=\left( \sum_{a=1}^N|\mathbf{w}_a|^2\right)^{\frac{1}{2}}.\]
Next,  we derive a system of dissipative inequalities for position and velocity configurations.

\begin{lemma} \label{L5.1}
Let $\{(\mathbf{x}_a, \mathbf{w}_a, T_a) \}$ be a global solution to \eqref{RTCSmodel(same num)-3} with zero initial momentum:
\[ \langle \mathbf{w}(0) \rangle  = 0. \]
Then, for $c$ suitably large such that
\begin{align}\label{condtion-c}
& \frac{2 \chi +3}{2c^2}+\frac{\mathcal{O}(c^{-4})}{\underline{T}}\leq \frac{1}{2\overline{T}},
\end{align} $\| \mathbf{X} \|$ and $\| \mathbf{W} \|$ satisfy
\[ \Big| \frac{\mathrm{d} \| \mathbf{X} \|}{\mathrm{d}t} \Big| \leq \frac{c^2}{1+c^2} \| \mathbf{W} \|, \qquad \frac{d| \mathbf{W} \|}{\mathrm{d}t}\le -\frac{1}{\overline{T}} \| \mathbf{W} \|, \quad \mbox{a.e.}~t > 0.
\]
\end{lemma}
\begin{proof} 
\noindent (i)~The first inequality follows from the Cauchy-Schwarz inequality and the following relation in \cite{Ha-Kim-Ruggeri-ARMA-2020}:
\[  |\mathbf{v}_a|\le\frac{c^2|\mathbf{w}_a|}{1+c^2}. \]

\noindent (ii)~We take an inner product between $\mathbf{w}_a$ and \eqref{RTCSmodel(same num)-3}$_2$,  and sum up the resulting relation over  $a \in[N]$ to get
\begin{align}\label{w 0(RTCS-3)}
\frac{\mathrm{d}}{\mathrm{d}t}\left(\sum_{a=1}^{N}|\mathbf{w}_a|^2\right)
=2\sum_{a=1}^{N}\mathbf{w}_a\cdot\frac{\mathrm{d}\mathbf{w}_a}{\mathrm{d}t}
=\frac{2}{N}\sum_{a,b=1}^{N}\mathbf{w}_a\cdot\left(\frac{\mathbf{w}_b\Gamma_b}{F_bT_b}-\frac{\mathbf{w}_a\Gamma_a}{F_aT_a}\right)
=\frac{2}{N}\sum_{a,b=1}^{N}\frac{-|\mathbf{w}_a|^2\Gamma_a}{F_aT_a}
\end{align}
by $\langle \mathbf{w}(0) \rangle  = 0$. Similar to \eqref{dv(TCS 3)}, we use \eqref{condtion-c} 
to have
\begin{align} 
\begin{aligned} \label{dw(RTCS-3)}
	\frac{\mathrm{d}}{\mathrm{d}t}\left(\sum_{a=1}^{N}|\mathbf{w}_a|^2\right) &=-\frac{2}{\overline{T}}\sum_{a=1}^{N} |\mathbf{w}_a|^2 +\frac{2}{N}\sum_{a,b=1}^{N}
	\frac{\mathbf{w}_b\cdot\mathbf{w}_a-|\mathbf{w}_a|^2}{T_a}
	\left(\frac{\Gamma_a}{F_a}-1\right) \\
	&\le-\frac{2}{\overline{T}}\sum_{a=1}^{N} |\mathbf{w}_a|^2 + \sum_{a=1}^{N}\frac{2|\mathbf{w}_a|^2}{T_a}
	\left(\frac{(2 \chi +3)T_a}{2c^2}+\mathcal{O}(c^{-4})\right)\\
	&\leq-\frac{1}{\overline{T}}\sum_{a=1}^{N} |\mathbf{w}_a|^2.
\end{aligned}
\end{align}

\end{proof}
As a direct application of Lemma \ref{L5.1}, one has the following result. 

\begin{corollary} \label{C5.1}
Suppose the condition \eqref{condtion-c} holds, and let $\left\{\left(\mathbf{x}_a, \mathbf{w}_a, T_a\right)\right\}$ be a global solution to  \eqref{RTCSmodel(same num)-3}. Then, one has 
\[
\| \mathbf{X}(t) \| \leq \| \mathbf{X}(0) \|  +  \frac{2c^2\overline{T} }{1+c^2}\| \mathbf{W}(0) \|, \quad \| \mathbf{W}(t) \|  \le  \| \mathbf{W}(0) \|  e^{-\frac{t}{2\overline{T}}}, \quad t \geq 0.
\]
\end{corollary}
\noindent To derive exponential time-decay estimates for temperature and velocity configurations, we first choose two positive constants $\tilde{A}_1$ and $\tilde{\lambda}_1$ such that 
\begin{align}
\begin{aligned} \label{F-2}
\hspace{-0.8cm} \tilde{A}_1 &> \overline{T}  + \frac{(N-1)^2( \overline{T} - \underline{T} )^2 \overline{T}}{4N^2 \underline{T} T^\infty } +\frac{[1+\mathcal{O}(c^{-2})] \| {\mathbf W} (0)\|^2}{(2\chi +1)^2N^2\overline{T}}+\mathcal{O}(c^{-2}), \\
\hspace{-0.8cm} \tilde{\lambda}_1 & = \min \Big\{ \frac{4+\mathcal{O}(c^{-2})}{(2\chi +1) \overline{T}^2},~  \frac{2}{ \overline{T}} -\frac{1}{\tilde{A}_1} \left(2+\frac{(N-1)^2( \overline{T} - \underline{T})^2}{2N^2 \underline{T} T^\infty} \right.\\ 
&\hspace{2cm}  \left.+\frac{[2+\mathcal{O}(c^{-2})] \| {\mathbf W} (0)\|^2}{(2\chi +1)^2N^2\overline{T}^2}+\mathcal{O}(c^{-2})\right) \Big\}.
\end{aligned}
\end{align}

\begin{lemma}\label{L5.2}
Suppose the condition \eqref{condtion-c} holds, and let $\left\{\left(\mathbf{x}_a, \mathbf{w}_a, T_a\right)\right\}$ be a global solution to  \eqref{RTCSmodel(same num)-3}. For positive constants $\tilde{A}_1$ and $\tilde{\lambda}_1$ satisfying \eqref{F-2}, one has 
\begin{align*}
\| \hat {\mathbf T}(t) \|^2 + \tilde{A}_1 \| {\mathbf W}(t) \| ^2 \le \Big(  \| \hat{\mathbf T}(0) \|^2 + \tilde{A}_1 \| \mathbf{W}(0) \|^2 \Big)e^{-\tilde{\lambda}_1t}, \quad t \geq 0.
\end{align*}
\end{lemma}
\begin{proof}
It follows from $\langle \mathbf{W}(0) \rangle  = 0$ , \eqref{RTCSmodel(same num)-3}, and \eqref{w 0(RTCS-3)} that 
\begin{align}\label{dT(RTCS-3)}
\frac{2 \chi +1}{2}  \frac{\mathrm{d}T_a}{\mathrm{d}t}
=& \frac{\mathrm{d}}{\mathrm{d}t}\left(\frac{2 \chi +1}{2}{T_a}+\frac{1}{2}|{\mathbf{w}_a}|^2+\frac{1}{c^2} {\mathcal F}_a\right)-\mathbf{w}_a\cdot\frac{\mathrm{d}\mathbf{w}_a}{\mathrm{d}t}-\frac{1}{c^2} \frac{\mathrm{d}{\mathcal F}_a}{\mathrm{d}t}
\\=&\frac{1}{N}\sum_{b=1}^{N}\left(\frac{\Gamma_a}{T_a}-\frac{\Gamma_b}{T_b}\right)-\frac{1}{N}\mathbf{w}_a\cdot\sum_{b=1}^{N}\left(\frac{\mathbf{w}_b\Gamma_b}{F_bT_b}-\frac{\mathbf{w}_a\Gamma_a}{F_aT_a}\right)-\frac{1}{c^2} \frac{\mathrm{d} {\mathcal F}_a}{\mathrm{d}t}\nonumber\\
=&\frac{1}{N}\sum_{b=1}^{N}\left(\frac{\Gamma_a}{T_a}-\frac{\Gamma_b}{T_b}\right)
+\frac{1}{N}\sum_{b=1}^{N}\frac{|\mathbf{w}_a|^2\Gamma_a}{F_aT_a}-\frac{1}{c^2} \frac{\mathrm{d} {\mathcal F}_a}{\mathrm{d}t}.\nonumber
\end{align}
For the positive constant $\tilde{A}_1$ given in \eqref{F-2}, we  obtain 
\begin{align} 
\label{cT14}
&\frac{\mathrm{d}}{\mathrm{d}t}\sum_{a=1}^N
\left( \frac{2 \chi +1}{2}|T_a-T^\infty|^2+ \tilde{A}_1 |\mathbf{w}_a|^2\right ) \nonumber\\ 
& \hspace{1cm} =\frac{2}{N}\sum_{a,b=1}^{N}\Big[(T_\infty-T_a)\left(\frac{\mathbf{w}_a\cdot
	\mathbf{w}_b}{T_b}\left(\frac{\Gamma_b}{F_b}-1\right)-\frac{|\mathbf{w}_a|^2}{T_a}
\left(\frac{\Gamma_a}{F_a}-1\right)\right) \nonumber\\
& \hspace{1cm} +(T_a-T_\infty) \left(\frac{\Gamma_a}{T_a}-\frac{\Gamma_b}{T_b}\right) \Big]
+\frac{2}{N}\sum_{a,b=1}^{N}\Big[\frac{-\tilde{A}_1
	|\mathbf{w}_a|^2}{T_a}\left(\frac{\Gamma_a}{F_a}-1\right)
\\
&\hspace{1cm} +\left(T_\infty-T_a+\tilde{A}_1\right)\left(\frac{\mathbf{w}_a\cdot\mathbf{w}_b}{T_b}-
\frac{|\mathbf{w}_a|^2}{T_a}\right)\Big] +\frac{1}{c^2}\sum_{a=1}^N(T_\infty-T_a)\frac{\mathrm{d} {\mathcal F}_a}{\mathrm{d}t}\nonumber
\\  
&\hspace{1cm} =:\mathcal{I}_{61}+\mathcal{I}_{62}+\mathcal{I}_{63}.\nonumber
\end{align}
In what follows, we estimate the terms ${\mathcal I}_{6i}$ one by one. \newline

\noindent $\bullet$~Case H.1 (Estimate of $\mathcal{I}_{61}$):  As in the estimation of $\mathcal{I}_{31}$, we use Lemma \ref{L3.2}, \eqref{Tinfty-varT b} and \eqref{Ta-Tb(RTCS)} to get
\begin{align*}  
\mathcal{I}_{61} &= \frac{2}{N}\sum_{a,b=1}^{N}(T_a-T^\infty) \left(\frac{1}{T_a}-\frac{1}{T_b}\right)+\frac{2}{N}\sum_{a,b=1}^{N}(T_a-T^\infty) \left(\frac{\Gamma_a-1}{T_a}-\frac{\Gamma_b-1}{T_b}\right)
\\ 
&+\frac{2}{N}\sum_{a,b=1}^{N}\mathbf{w}_a\cdot(T^\infty-T_a)\left(\frac{\mathbf{w}_b}{T_b}\left(\frac{\Gamma_b}{F_b}-1\right)-\frac{\mathbf{w}_a}{T_a}\left(\frac{\Gamma_a}{F_a}-1\right)\right)
\\ 
&\le -\frac{2+\mathcal{O}(c^{-2})}{\overline{T}^2}\left(\sum_{a=1}^{N} |T_a-T^\infty|^2-\frac{\| {\mathbf W}(0) \|^2}{(2\chi +1)^2N^2}\| {\mathbf W}(t) \|^2\right) \\
& +\mathcal{O}(c^{-2})\sum_{a=1}^N\left(|\mathbf{w}_a|^2+|T_a-T^\infty|^2\right)  \\ 
&\le -\frac{4+\mathcal{O}(c^{-2})}{(2\chi +1)\overline{T}^2}  \| \hat {\mathbf T} \|^2 +\frac{[2+\mathcal{O}(c^{-2})] \| {\mathbf W} (0)\|^2}{(2\chi +1)^2N^2\overline{T}^2}  \| {\mathbf W} \|^2. 
\end{align*}
\vspace{0.2cm}


\noindent $\bullet$~Case H.2 (Estimate of $\mathcal{I}_{62}$): We use Lemma \ref{L3.2} to obtain
\begin{align*}
\begin{aligned}
	\mathcal{I}_{62} &\le \left(2+\frac{(N-1)^2(\overline{T} - \underline{T})^2}{2N^2 \underline{T}T^\infty}\right)\sum_{a=1}^{N}|\mathbf{w}_a|^2
	-2 \tilde{A}_1 \sum_{a=1}^N\frac{|\mathbf{w}_a|^2}{\overline{T}}
	\\
	& \hspace{1cm} +\frac{2 \tilde{A}_1}{N}\sum_{a,b=1}^{N}\frac{\mathbf{w}_a\cdot\mathbf{w}_b-|\mathbf{w}_a|^2}{T_a}
	\left[ -\frac{(2\chi +3)T_a}{2c^2}+\mathcal{O}(c^{-4})\right]
	\\ 
	&\leq  \left(2+\frac{(N-1)^2( \overline{T} - \underline{T})^2}{2N^2 \underline{T} T^\infty}-\frac{2 \tilde{A}_1}{\overline{T}}+\mathcal{O}(c^{-2})\right)\sum_{a=1}^{N}|\mathbf{w}_a|^2.
\end{aligned}
\end{align*}

\noindent $\bullet$~Case H.3 (Estimate of $\mathcal{I}_{63}$): We use  Lemma \ref{L3.2} again to get 
\begin{align*}
\mathcal{I}_{63}\le \mathcal{O}(c^{-2})\sum_{a=1}^N\left(|\mathbf{w}_a|^2+|T_a-T_{\infty}|^2\right).
\end{align*}
In \eqref{cT14}, we collect all the estimates in Case H1 - Case H.3 to obtain the desired estimate:
\begin{align*}
&\frac{\mathrm{d}}{\mathrm{d}t} \left(   \| \hat {\mathbf T} \|^2 + \tilde{A}_1 \| {\mathbf W} \|^2  \right)
+\frac{4+\mathcal{O}(c^{-2})}{(2 \chi +1) \overline{T}^2}  \| \hat {\mathbf T} \|^2 \\
&\hspace{1cm} -\left(2+\frac{(N-1)^2( \overline{T} - \underline{T})^2}{2N^2 \underline{T} T^\infty}
+\frac{[2+\mathcal{O}(c^{-2})] \| {\mathbf W} (0)\|^2}{(2\chi +1)^2N^2\overline{T}^2}-\frac{2 \tilde{A}_1}{\overline{T}}+\mathcal{O}(c^{-2})\right) \| {\mathbf W} \|^2  \leq 0.
\end{align*}
\end{proof}
\noindent Next, we are ready to present asymptotic flocking estimates as follows. 
\begin{theorem}\label{T5.1}
Suppose the condition \eqref{condtion-c} holds, and let $\left\{\left(\mathbf{x}_a, \mathbf{w}_a, T_a\right)\right\}$ be a global solution to \eqref{RTCSmodel(same num)-3}. Then for positive constants $\tilde{A}_1, {\tilde \lambda}_1$ in \eqref{F-2}, we have
\begin{align*}
\begin{aligned}
	& \sup_{0 \leq t < \infty} \| \mathbf{X}(t) \| \leq \| \mathbf{X}(0) \|  +  \frac{2c^2\overline{T} }{1+c^2}\| \mathbf{W}(0) \| , \quad \| \mathbf{W}(t) \|  \le  \| \mathbf{W}(0) \|  e^{-\frac{t}{2\overline{T}}}, \\
	& \| \hat {\mathbf T}(t) \|^2 + {\tilde A}_1 \| {\mathbf V}(t) \| ^2 \le \Big(  \| \hat{\mathbf T}(0) \|^2 + \tilde{A}_1 \| \mathbf{V}(0) \|^2 \Big)e^{- \tilde{\lambda}_1t}, \quad t \geq 0.
\end{aligned}
\end{align*}
\end{theorem}
\begin{proof}
We combine Lemma \ref{L5.1} and Lemma \ref{L5.2} to derive the desired estimates. 
\end{proof}

\subsection{Flocking dynamics II} \label{sec:5.2}
In this subsection, we present asymptotic flocking dynamics of  \eqref{RTCSmodel(same num)-2}, where the communication weight functions $\phi_{ab}$ with $a, b\in[N]$ satisfy \eqref{phiMm}. We show that for sufficiently small  $\varepsilon$ and sufficiently large   $c$, flocking estimates hold.

\begin{lemma} \label{L5.3}
Suppose that positive constants $\varepsilon$ and $c$ satisfy 
\begin{align}\label{c-varep}
\frac{2\varepsilon}{ \underline{T} }+\mathcal{O}(c^{-2})\leq \frac{\underline{\phi}} {\overline{T} },
\end{align}
and let $\{(\mathbf{x}_a, \mathbf{w}_a, T_a) \}$ be a global solution to \eqref{RTCSmodel(same num)-2} with zero initial momentum:
\[ \langle {\mathbf w}(0) \rangle  = 0. \]
Then, $\| \mathbf{X} \|$ and $\| \mathbf{W} \|$ satisfy
\[
\Big| \frac{\mathrm{d} \| \mathbf{X} \|}{\mathrm{d}t} \Big| \leq \frac{c^2}{1+c^2} \| \mathbf{W} \|, \quad \frac{\mathrm{d}| \mathbf{W} \|}{\mathrm{d}t}\le -\frac{\underline{\phi}}{2\overline{T}} \| \mathbf{W} \|, \quad \mbox{a.e.}~~t > 0.
\]
\end{lemma}
\begin{proof}
We only derive the second differential inequality. Similar to \eqref{dw(RTCS-3)}, we take an inner product of $\mathbf{w}_a$ with \eqref{RTCSmodel(same num)-2}$_2$, sum up the resulting relation over $a\in[N]$, and use \eqref{c-varep} to get 
\begin{align*}
\begin{aligned}
	\frac{\mathrm{d}}{\mathrm{d}t}\left(\sum_{a=1}^{N}|\mathbf{w}_a|^2\right) &\le-\frac{2\underline{\phi}}{ \overline{T} }\sum_{a=1}^{N} |\mathbf{w}_a|^2 +\frac{2\varepsilon}{\underline{T}}\sum_{a=1}^{N} |\mathbf{w}_a|^2 \\
	& \hspace{0.2cm} +\frac{2}{N}\sum_{a,b=1}^{N}\phi_{ab}
	\frac{\mathbf{w}_a\cdot\mathbf{w}_b-|\mathbf{w}_a|^2}{T_a}
	\left(-\frac{(2 \chi +3)T_a}{2c^2}+\mathcal{O}(c^{-4})\right)  \\ 
	&\le -\left(\frac{2 \underline{\phi}}{\overline{T}}-\frac{2\varepsilon}{\underline{T}}-\mathcal{O}(c^{-2})\right)\sum_{a=1}^{N} |\mathbf{w}_a|^2 \\
	&\leq -\frac{\underline{\phi}}{\overline{T}}\sum_{a=1}^{N} |\mathbf{w}_a|^2,
\end{aligned}
\end{align*}
for sufficiently small  $\varepsilon$ and sufficiently large   $c$.
\end{proof}
As a direct application of Lemma \ref{L5.3}, we obtain the following estimate. 

\begin{corollary}  \label{C5.2}
Suppose that the communication weight functions, initial data and $\varepsilon$ satisfy \eqref{c-varep}
and let $\left\{\left(\mathbf{x}_a, \mathbf{w}_a, T_a\right)\right\}$ be a global solution to \eqref{RTCSmodel(same num)-2}. Then, one has
\[
\sup_{0 \leq t < \infty} \| \mathbf{X}(t) \| \le \| \mathbf{X}(0) \|+\frac{2c^2\overline{T} }{ (1+c^2)\underline{\phi}} \| \mathbf{W}(0)\|, \quad \| \mathbf{W}(t) \| \le  \| \mathbf{W}(0) \| e^{-\frac{\underline{\phi}}{ 2\overline{T}}t}, \quad t \geq 0.
\]
\end{corollary}
Next, we derive the  exponential time-decay estimates for temperature and velocity. For this, we set positive constants $\tilde{A}_2$ and $\tilde{\lambda}_2$ such that 
\begin{align}
\begin{aligned} \label{New-E-0}
\tilde{A}_2  &> \frac{4}{7}\Big[\frac{8 \underline{\zeta} \| {\mathbf W}(0) \|^2}{(2\chi +1)^2 \underline{\phi} \overline{T}}+\left(2 +\frac{(N-1)^2( \overline{T} - \underline{T})^2}{2N^2 \underline{T} T^\infty}\right) \overline{T}  + \frac{\varepsilon \overline{T} }{\underline{\phi}}\left(2+\frac{(T^{\infty}- \underline{T})^2}{2 \underline{T} T^\infty} \right)
\Big], \\ 
\tilde{\lambda}_2 & := \min\Big\{\frac{4 \underline{\zeta} +\mathcal{O}(c^{-2})}{(2 \chi +1) \overline{T}^2}, ~~\frac{7 \underline{\phi}}{4 \overline{T}}
-\frac{1}{\tilde{A}_2}\Big[\left(2+\frac{(N-1)^2(\overline{T} - \underline{T})^2}{2N^2 \underline{T} T^\infty}\right) \underline{\phi}  \\
& \hspace{5.5cm} +\frac{8 \underline{\zeta} \| {\mathbf W}(0) \|^2 }{(2\chi+1)^2 \overline{T}^2}
+\varepsilon\left(2+\frac{(T^{\infty}- \underline{T})^2}{2 \underline{T} T^\infty}\right)\Big]\Big\}.
\end{aligned}
\end{align}
\begin{lemma}\label{L5.4} Suppose the condition \eqref{c-varep} holds, and let $\left\{ (\mathbf{x}_a, \mathbf{w}_a, T_a)\right \}$ be a global solution to \eqref{RTCSmodel(same num)-2}. Then, for positive constants ${\tilde A}_2$  and ${\tilde \lambda}_2$ satisfying  \eqref{New-E-0}, we have
\begin{equation*}\label{diss-vT(TCS-2)}
\| \hat {\mathbf T}(t) \|^2 + {\tilde A}_2 \| {\mathbf W}(t) \| ^2  \le \Big(  \| \hat {\mathbf T}(0) \|^2 + {\tilde A}_2 \| {\mathbf W}(0) \| ^2      \Big)e^{-{\tilde \lambda}_2t}, \quad t \geq 0.
\end{equation*}
\end{lemma}
\begin{proof}
Similar to \eqref{dT(RTCS-3)}, it follows from \eqref{RTCSmodel(same num)-2} that 
\begin{align*}
\frac{2\chi +1}{2}  \frac{\mathrm{d}T_a}{\mathrm{d}t}
&= \frac{1}{N}\sum_{b=1}^{N}\zeta_{ab}\left(\frac{\Gamma_a}{T_a}-\frac{\Gamma_b}{T_b}\right)-\frac{\underline{\phi}}{N}\mathbf{w}_a\cdot\sum_{b=1}^{N}\left(\frac{\mathbf{w}_b\Gamma_b}{F_bT_b}-\frac{\mathbf{w}_a\Gamma_a}{F_aT_a}\right)
\\
& \hspace{0.2cm} -\frac{1}{N}\mathbf{w}_a\cdot\sum_{b=1}^{N}(\phi_{ab}-\underline{\phi})\left(\frac{\mathbf{w}_b\Gamma_b}{F_bT_b}-\frac{\mathbf{w}_a\Gamma_a}{F_aT_a}\right)-\frac{1}{c^2}\frac{\mathrm{d} {\mathcal F}_a}{\mathrm{d}t}.
\end{align*}
For the positive constant ${\tilde A}_2$ in \eqref{New-E-0}, we obtain
\begin{align}
\begin{aligned} \label{cT15}
	& \frac{\mathrm{d}}{\mathrm{d}t}\sum_{a=1}^N\left(\frac{2\chi +1}{2}|T_a-T^\infty|^2
	+{\tilde A}_2|\mathbf{w}_a|^2\right) 
	\\  
	& \hspace{0.2cm} =\frac{2}{N}\sum_{a,b=1}^{N}\Big[\phi_{ab}(T^\infty-T_a)\left(\frac{\mathbf{w}_a\cdot\mathbf{w}_b}{T_b}\left(\frac{\Gamma_b}{F_b}-1\right)-\frac{|\mathbf{w}_a|^2}{T_a}\left(\frac{\Gamma_a}{F_a}-1\right)\right) 
	\\
	& \hspace{0.2cm}  +\zeta_{ab}(T_a-T^\infty) \left(\frac{\Gamma_a}{T_a}-\frac{\Gamma_b}{T_b}\right)\Big]
	+\frac{2}{N}\sum_{a,b=1}^{N}\Big[ \underline{\phi} \left(T^\infty-T_a+A_5\right)
	\\
	& \hspace{0.2cm}   \times \left(\frac{\mathbf{w}_a\cdot\mathbf{w}_b}{T_b}-\frac{|\mathbf{w}_a|^2}{T_a}\right)
	+{\tilde A}_2\phi_{ab}\frac{\mathbf{w}_a\cdot\mathbf{w}_b-|\mathbf{w}_a|^2}{T_a}\left(\frac{\Gamma_a}{F_a}-1\right)\Big] 
	\\
	& \hspace{0.2cm}  +\Big[\frac{2}{N}\sum_{a,b=1}^{N}(\phi_{ab} - \underline{\phi})\left(T^\infty-T_a+{\tilde A}_2\right) \left(\frac{\mathbf{w}_a\cdot\mathbf{w}_b}{T_b}-\frac{|\mathbf{w}_a|^2}{T_a}\right) 
	+\frac{1}{c^2}\sum_{a=1}^N(T_\infty-T_a)\frac{\mathrm{d} {\mathcal F}_a}{\mathrm{d}t}\Big]\\
	& \hspace{0.2cm}  =: \mathcal{I}_{71}+\mathcal{I}_{72}+\mathcal{I}_{73}.
\end{aligned}
\end{align}
Next, we estimate the terms ${\mathcal I}_{7i}$ one by one. \newline

\noindent $\bullet$~Case I.1 (Estimation of $\mathcal{I}_{71}$): By direct calculation, we have
\begin{align*}
\mathcal{I}_{71}  \le& -\frac{2 \underline{\zeta} +\mathcal{O}(c^{-2})}{(2\chi +1)\overline{T}^2}  \| \hat {\mathbf T}(t) \|^2 +\frac{[2 \underline{\zeta} +\mathcal{O}(c^{-2})]  \| \mathbf{W} (0)\|^2}
{(2\chi +1)^2N^2\overline{T}^2} \| \mathbf{W} \|^2.
\end{align*}
\noindent $\bullet$~Case I.2 (Estimation of $\mathcal{I}_{72}$): Similarly, we have
\begin{align*}
\mathcal{I}_{72}
\le \left(2+\frac{(N-1)^2( \overline{T} - \underline{T})^2}{2N^2 \underline{T} T^\infty}-\frac{2 \tilde{A}_2}{\overline{T}}
+\mathcal{O}(c^{-2})\right) \underline{\phi}  \| \mathbf{W} \|^2.
\end{align*}

\noindent $\bullet$~Case I.3 (Estimation of $\mathcal{I}_{73}$):~We use 
\[
\frac{1}{2N}\sum_{a,b=1}^{N}(\phi_{ab}- \underline{\phi})\frac{(T^{\infty}-T_a)^2}{T^{\infty}  \underline{T}} |\mathbf{w}_b|^2
\le \frac{\varepsilon(T^{\infty}- \underline{T})^2}{2 \underline{T}}   \| \mathbf{W} \|^2
\]
to find 
\begin{align*}
\begin{aligned}
	\mathcal{I}_{73}
	&\le \frac{2}{N}\sum_{a,b=1}^{N}(\phi_{ab}- \underline{\phi})\left[T^\infty\frac{ |\mathbf{w}_b|^2}{T_b}+\frac{(T^{\infty}-T_a)^2}{4T^{\infty}}\frac{ |\mathbf{w}_b|^2}{\underline{T}}\right]+2\varepsilon\sum_{a=1}^{N}|\mathbf{w}_a|^2
	\\
	&+\frac{2A_5}{N}\sum_{a,b=1}^N(\phi_{ab}-\underline{\phi})\left(\frac{|\mathbf{w}_a|^2}{T_a}
	+\frac{|\mathbf{w}_b|^2}{4 \underline{T}}\right)
	\\
	&-\frac{2}{N}\sum_{a,b=1}^N\left(\phi_{ab}- \underline{\phi} \right)(\tilde{A}_2 +T^\infty)
	\frac{|\mathbf{w}_a|^2}{T_a}+\mathcal{O}(c^{-2})\sum_{a=1}^N\left(|\mathbf{w}_a|^2+|T_a-T^{\infty}|^2\right)
	\\  
	&\le \left[\varepsilon\left(2+\frac{(T^{\infty}- \underline{T})^2}{2 \underline{T} T^\infty}
	+\frac{\tilde{A}_2}{2 \underline{T} }\right)
	+\mathcal{O}(c^{-2})\right]  \| \mathbf{W} \|^2 +\mathcal{O}(c^{-2})  \|\hat {\mathbf T}(t) \|^2.
\end{aligned}
\end{align*}
In \eqref{cT15}, we combine all the estimates in Case I.1 - Case I.3 to find 
\begin{align*}
\begin{aligned}
	&\frac{\mathrm{d}}{\mathrm{d}t}\Big(   \| \hat {\mathbf T} \|^2 + {\tilde A}_2 \| {\mathbf W} \| ^2 \Big)
	+\frac{4 \underline{\zeta} +\mathcal{O}(c^{-2})}{(2\chi +1) \overline{T}^2}  \| \hat {\mathbf T} \|^2   +\Big[-\left(2+\frac{(N-1)^2(\overline{T} - \underline{T})^2}{2N^2 \underline{T} T^\infty}
	\right) \underline{\phi} \\
	&\hspace{1cm} -\frac{[2 \underline{\zeta} +\mathcal{O}(c^{-2})]  \| \mathbf{W} (0)\|^2}
	{(2\chi +1)^2N^2\overline{T}^2}-\varepsilon\left(2+\frac{(T^{\infty}- \underline{T} )^2}{2 \underline{T} T^\infty}
	\right) +\left(\frac{2 \underline{\phi} }{\overline{T} }-\frac{\varepsilon}{2 \underline{T} }-\mathcal{O}(c^{-2})\right) \tilde{A}_2 \Big] \| {\mathbf W} \| ^2\leq  0.
\end{aligned}
\end{align*}
This implies the desired estimate.
\end{proof}
Now, we are ready to provide the emergence of asymptotic flocking by estimating position, velocity, and temperature under small perturbations.

\begin{theorem}\label{T5.2}
Suppose the condition \eqref{c-varep} holds, and let $\left\{\left(\mathbf{x}_a, \mathbf{w}_a, T_a\right)\right\}$ be a global solution to \eqref{RTCSmodel(same num)-2}. For positive constants $\tilde{A}_2, {\tilde \lambda}_2$ in \eqref{New-E-0}, we have
\begin{align*}
\begin{aligned}
	& \sup_{0 \leq t < \infty} \| \mathbf{X}(t) \| \leq \| \mathbf{X}(0) \|  +  \frac{2\overline{T}}{\underline{\phi}} {\mathbf W}(0), \quad \| \mathbf{W}(t) \|  \le  \| \mathbf{W}(0) \|  e^{-\frac{ \underline{\phi} t}{2\overline{T}}},  \\
	& \| \hat {\mathbf T}(t) \|^2 + {\tilde A}_2 \| {\mathbf W}(t) \| ^2 \le \Big(  \| \hat{\mathbf T}(0) \|^2 + \tilde{A}_2 \| \mathbf{W}(0) \|^2 \Big)e^{- \tilde{\lambda}_2 t},\quad t \geq 0.
\end{aligned}
\end{align*}
\end{theorem}
\begin{proof}
We combine Lemma \ref{L5.3} and Lemma \ref{L5.4} to derive the desired estimates. 
\end{proof}
\subsection{Flocking dynamics III} \label{sec:5.3}
In this subsection, we show the asymptotic flocking of the RTCS model \eqref{RTCSmodel(same num)-1} for the sufficiently large speed velocity $c$ and the well-prepared initial datum. Next, we derive a system of dissipative inequality for the position and the velocity.
\begin{lemma}\label{L5.5}
Let $\left\{\left(\mathbf{x}_a, \mathbf{w}_a, T_a\right)\right\}$ be a global solution to \eqref{RTCSmodel(same num)-1}. Then, we have
\[
\Big| \frac{\mathrm{d} \| \mathbf{X} \|}{\mathrm{d}t} \Big| \le \frac{c^2}{1+c^2} \| \mathbf{W} \|, \quad \frac{d\| \mathbf{W} \|}{\mathrm{d}t}\le -\frac{1}{2}\left(\frac{\phi(0)}{\overline{T}}
-\frac{(\phi(0)- \underline{\phi})^2}{\phi(0) \underline{T}}-\mathcal{O}(c^{-2})\right)  \| \mathbf{W} \|, \quad \mbox{a.e.}~ t > 0.
\]
\end{lemma}
\begin{proof}
We only derive the second differential inequality since the first one is obvious. 
Similar to \eqref{dv(TCS 1)} and \eqref{dw(RTCS-3)}, we take an inner product of $\mathbf{w}_a$ with \eqref{RTCSmodel(same num)-1}$_2$ and sum up the resulting relation over  $a\in[N]$ to get
\begin{align*}
\frac{\mathrm{d}}{\mathrm{d}t} \| \mathbf{W} \|^2 &\le
-\frac{1}{N}\sum_{a,b=1}^{N}\left[\frac{\phi_{ab}}{T_a}-
\frac{(\phi_{ab}- \underline{\phi})^2}{\phi_{ab} \underline{T}}\right]|\mathbf{w}_a|^2 
\\
&+\frac{2}{N}\sum_{a,b=1}^{N}\phi_{ab}
\frac{\mathbf{w}_b\cdot\mathbf{w}_a-|\mathbf{w}_a|^2}{T_a}
\left(-\frac{(2\chi +3)T_a}{2c^2}+\mathcal{O}(c^{-4})\right)
\\ 
&\le-\left(\frac{\phi(0)}{\overline{T}}
-\frac{(\phi(0)- \underline{\phi})^2}{\phi(0) \underline{T}}-\phi(0)\mathcal{O}(c^{-2})\right)  \| \mathbf{W} \|^2.
\end{align*}
\end{proof}
Next, we define 
\begin{align*}
2\chi_1(\| \mathbf{X} \|) :=\left[\frac{1}{\overline{T} }-\mathcal{O}(c^{-2})\right]\phi(0)
-\frac{(\phi(0)-\phi(\sqrt{2}\| \mathbf{X} \|))^2}{\phi(0) \underline{T}}.
\end{align*}
As a direct application of Lemma \ref{L5.5}, we have the following estimate. 
\begin{corollary} \label{C5.3}
Suppose that for given initial data, there exists a positive constant $U_1$ such that
\begin{align}\label{phiX phi0(RTCS)}
\phi(U_1)>\phi(0)\left(1-\sqrt{\frac{\underline{T}}{\overline{T}}- \underline{T} \mathcal{O}(c^{-2})}\right),\quad U_1\leq \sqrt{2} \| \mathbf{X}(0) \| +\frac{\sqrt{2}c^2}{(1+c^2)\chi_1(U_1)}  \| \mathbf{W}(0) \|,
\end{align}
and let $\left\{\left(\mathbf{x}_a, \mathbf{w}_a, T_a\right)\right\}$ be a global solution to \eqref{RTCSmodel(same num)-1}. Then, one has 
\[
\sup_{0 \leq t < \infty} \sqrt{2}\| \mathbf{X}(t) \| \le U_1,\quad  \| \mathbf{W} (t) \| \le  \| \mathbf{W}(0) \| e^{-\chi_1(U_1) t}, \quad t \geq 0.
\]
\end{corollary}
For given initial data, we set positive constants $\tilde{A}_3$ and $\tilde{\lambda}_3$ as follows.
\begin{align}\label{AL3-r}
{\tilde A}_3 &> \frac{1}{2\chi_1(U_1)}\left[\frac{[2\zeta(U_1)+\mathcal{O}(c^{-2})]  \| {\mathbf W}(0) \|^2}
{(2\chi +1)^2N^2\overline{T}^2} +\left(2+\frac{(\overline{T} - \underline{T})^2}{2 \underline{T} T^\infty}\right)\phi(0)\right], \nonumber
\\ 
{\tilde \lambda}_3 &:=\min\Big\{\frac{4\zeta(U_1)+\mathcal{O}(c^{-2})}{(2 \chi +1)\overline{T}^2},
2\chi_1(U_1)-\frac{1}{{\tilde A}_3}\Big[\frac{[2\zeta(U_1)+\mathcal{O}(c^{-2})]  \| {\mathbf W}(0) \|^2}
{(2\chi +1)^2N^2\overline{T}^2}\\
&\hspace{2cm}     +\left(2+\frac{(\overline{T} - \underline{T} )^2}{2 \underline{T} T^\infty}\right)\phi(0)\Big]\Big\}.\nonumber
\end{align}
\begin{lemma}\label{L5.6}
Suppose the conditions  \eqref{phiX phi0(RTCS)} hold, and  let $\{(\mathbf{x}_a, \mathbf{w}_a, T_a)\}$ be a global solution to \eqref{RTCSmodel(same num)-1}. Then one has 
\begin{equation*}
\| \hat {\mathbf T}(t) \|^2 + \tilde{A}_3 \| {\mathbf W}(t) \|^2  \le \Big(  \| \hat {\mathbf T}(0) \|^2 + {\tilde A}_3 \| {\mathbf W}(0) \| ^2      \Big)e^{-{\tilde \lambda}_3 t}, \quad t \geq 0.
\end{equation*}
\end{lemma}
\begin{proof}
Similar to \eqref{dT(RTCS-3)}, it follows from \eqref{RTCSmodel(same num)-1} that 
\begin{align*}
\frac{2\chi +1}{2} \frac{\mathrm{d}T_a}{\mathrm{d}t}
=\frac{1}{N}\sum_{b=1}^{N}\zeta_{ab}\left(\frac{\Gamma_a}{T_a}-\frac{\Gamma_b}{T_b}\right)-\frac{1}{N}\mathbf{w}_a\cdot\sum_{b=1}^{N}\phi_{ab}\left(\frac{\mathbf{w}_b\Gamma_b}{F_bT_b}-\frac{\mathbf{w}_a\Gamma_a}{F_aT_a}\right)-\frac{1}{c^2}\frac{\mathrm{d} {\mathcal F}_a}{\mathrm{d}t}.
\end{align*}
For the positive constant $\tilde{A}_3$ in \eqref{AL3-r}, we obtain 
\begin{align} 
\begin{aligned} \label{cT16}
	&\frac{\mathrm{d}}{\mathrm{d}t}\sum_{a=1}^N\left(\frac{2\chi+1}{2}|T_a-T^\infty|^2
	+A_6|\mathbf{w}_a|^2\right) 
	\\ 
	& \hspace{1cm} =\frac{2}{N}\sum_{a,b=1}^{N}\Big[\phi_{ab}(T^\infty-T_a)\left(\frac{\mathbf{w}_a\cdot\mathbf{w}_b}{T_b}\left(\frac{\Gamma_b}{F_b}-1\right)-\frac{|\mathbf{w}_a|^2}{T_a}\left(\frac{\Gamma_a}{F_a}-1\right)\right) 
	\\  
	& \hspace{1cm}  +\zeta_{ab}(T_a-T^\infty)\left(\frac{\Gamma_a}{T_a}-\frac{\Gamma_b}{T_b}\right)\Big]
	+\frac{2}{N}\sum_{a,b=1}^{N}\Big[ \tilde{A}_3 \phi_{ab}
	\frac{\mathbf{w}_a\cdot\mathbf{w}_b-|\mathbf{w}_a|^2}{T_a}
	\\
	& \hspace{1cm} \times \left(\frac{\Gamma_a}{F_a}-1\right)+\phi_{ab}\left(T^\infty-T_a+ \tilde{A}_3 \right)
	\left(\frac{\mathbf{w}_a\cdot\mathbf{w}_b}{T_b}-\frac{|\mathbf{w}_a|^2}{T_a}\right)\Big] 
	\\& \hspace{1cm}+\frac{1}{c^2} \sum_{a=1}^N(T^\infty-T_a)\frac{\mathrm{d} {\mathcal F}_a}{\mathrm{d}t} \\
	& \hspace{1cm} =:\mathcal{I}_{81}+\mathcal{I}_{82}+\mathcal{I}_{83}. 
\end{aligned}
\end{align}
Next, we estimate the term in ${\mathcal I}_{8i}$ one by one. \newline

\noindent $\bullet$~Case J.1 (Estimate of $\mathcal{I}_{81}$):~By direct estimate, we have
\begin{align*}
\mathcal{I}_{81}\le  -\frac{4\zeta(U_1)+\mathcal{O}(c^{-2})}{(2\chi +1)\overline{T}^2} \| \hat {\mathbf T}(t) \|^2 +\frac{[2\zeta(U_1)+\mathcal{O}(c^{-2})]  \| {\mathbf W}(0) \|^2}
{(2\chi +1)^2N^2\overline{T}^2}  \| {\mathbf W}(t) \|^2.
\end{align*}
\vspace{0.2cm}

\noindent $\bullet$~Case J.2 (Estimate of $\mathcal{I}_{82}$): As in $\mathcal{I}_{52}$ and $\mathcal{I}_{72}$, we can obtain
\begin{align*}
\begin{aligned}
	\mathcal{I}_{82}\le& \frac{2}{N}\sum_{a,b=1}^{N}\phi_{ab}\left[T^\infty\frac{ |\mathbf{w}_b|^2}{T_b}+\frac{(T^\infty-T_a)^2}{4T_\infty}\frac{ |\mathbf{w}_b|^2}{ \underline{T} }\right]
	\\ &+\frac{2 \tilde{A}_3}{N}\sum_{a,b=1}^N\left(\phi_{ab}\frac{|\mathbf{w}_a|^2}{2T_a}
	+\frac{|\phi_{ab}-\phi\left(U_1\right)|^2}{\phi_{ab}}\frac{|\mathbf{w}_b|^2}{2T_a}\right)\\
	&+\left[2\phi(0)+ \tilde{A}_3 \mathcal{O}(c^{-2})\right]\sum_{a=1}^{N}|\mathbf{w}_a|^2-\frac{2( \tilde{A}_3 +T^\infty)}{N}\sum_{a,b=1}^N\phi_{ab}\frac{|\mathbf{w}_a|^2}{T_a}\\
	&\le \left(2+\frac{(\overline{T} - \underline{T})^2}{2 \underline{T} T^\infty}\right)\phi(0)
	\sum_{a=1}^{N}|\mathbf{w}_a|^2\\
	&-\frac{\tilde{A}_3}{N}\sum_{a,b=1}^{N}\left[\frac{\phi_{ab}}{T_a}
	-\frac{(\phi_{ab}-\phi\left(U_1\right))^2}{\phi_{ab}T_b}-\phi(0)\mathcal{O}(c^{-2})\right]
	|\mathbf{w}_a|^2\\
	\leq& -\left[2\tilde{A}_3\chi_1(U_1)-\left(2+\frac{( \overline{T} - \underline{T})^2}{2 \underline{T} T^\infty}\right)\phi(0)\right] \| {\mathbf W}(t) \|^2.
\end{aligned}
\end{align*}
\vspace{0.2cm}

\noindent $\bullet$~Case J.3 (Estimate of $\mathcal{I}_{83}$): As in $\mathcal{I}_{63}$, we have
\begin{align*}
\mathcal{I}_{83}\leq \mathcal{O}(c^{-2}) \left( \| \hat {\mathbf T}(t) \|^2+\| {\mathbf W}(t) \|^2\right).
\end{align*}
In  \eqref{cT16}, we combine all the estimates of Case J.1 - Case J.3 to obtain
\begin{align*}
& \frac{\mathrm{d}}{\mathrm{d}t}\Big(   \| \hat {\mathbf T}(t) \|^2 + \tilde{A}_3 \| {\mathbf W}(t) \|^2 \Big)
+\frac{4\zeta(U_1)+\mathcal{O}(c^{-2})}{(2\chi+1)\overline{T}^2}  \| \hat {\mathbf T}(t) \|^2  +\Big[-\frac{[2\zeta(U_1)+\mathcal{O}(c^{-2})]  \| {\mathbf W}(0) \|^2}
{(2\chi +1)^2N^2\overline{T}^2} 
\\
&\qquad+ \tilde{A}_3\chi_1(U_1)-\left(2+\frac{(\overline{T} - \underline{T})^2}{2 \underline{T} T^\infty}\right)\phi(0)\Big]  \| {\mathbf W}(t) \|^2
\le  0.
\end{align*}
This implies the desired estimate.
\end{proof}

Next, we show the emergence of asymptotic flocking by the estimation of the position, the velocity, and the temperature under the well-prepared condition using Corollary \ref{C5.3} and Lemma \ref{L5.6}.
\begin{theorem}\label{T5.3}
Suppose the conditions \eqref{phiX phi0(RTCS)} hold, and  let $\{(\mathbf{x}_a, \mathbf{w}_a, T_a)\}$ be a global solution to \eqref{RTCSmodel(same num)-1}. For positive constants $\tilde{A}_3, \tilde{\lambda}_3$ in \eqref{AL3-r}, we have
\begin{align*}
\begin{aligned}
	& \sup_{0 \leq t < \infty} \sqrt{2}\| \mathbf{X}\| \le U_1, \quad \| \mathbf{W}(t) \| \le  \| \mathbf{W}(0) \| e^{-\chi_1(U_1) t}, \\
	&  \| \hat {\mathbf T}(t) \|^2 + \tilde{A}_3 \| {\mathbf W}(t) \|^2  \le \Big(  \| \hat {\mathbf T}(0) \|^2 + {\tilde A}_3 \| {\mathbf W}(0) \|^2 \Big)e^{-{\tilde \lambda}_3 t}, \quad t \geq 0.
\end{aligned}
\end{align*}
\end{theorem}

\section{Conclusion}\label{conclusion}  \label{sec:6}
\setcounter{equation}{0}
In this paper, we have studied the asymptotic flocking dynamics of the classical TCS model and the RTCS model with Synge energy corresponding to monatomic and polyatomic gases. Based on the entropy principle, we have succeeded to derive a uniform lower bound of the temperature for the TCS model and the RTCS model, which is related to the Synge energy, when the speed of light velocity $c$ is large. Thanks to the lower bound of temperature and upper bounds of velocities and temperatures derived by the conservation of energy, we can establish the asymptotic flocking for the the TCS model and the RTCS model for arbitrary initial data, when the communication weight functions $\phi$ and $\psi$ are constants or the function $\phi$ is close to a constant. For a general communication weight function, we also provided sufficient frameworks on the flocking dynamics of the TCS model and the RTCS model with Synge energy. Of course, there are several issues unexplored in this work. For instance, we have only provided sufficient framework for the flocking dynamics of the TCS model and the RTCS model. Thus, looking for necessary frameworks for asymptotic flocking will be also an interesting problem. We also have not considered the classical limit from the RTCS model with Synge energy to the classical TCS model. These interesting problems will be left for a future work. 

\newpage

\appendix

%
%

\section{Analytical properties of modified Bessel functions}  \label{App-A}
\setcounter{equation}{0}
In this appendix, we recall basic properties of the modified Bessel functions which were frequently used for the estimation of the RTCS model in main context, and  we discuss the derivation of the RTCS model for monatomic gases and polyatomic gases such as diatomic, triatomic, and tetratomic gases with Synge energy. First, we introduce modified Bessel functions and study their basic properties without proofs. \newline

Let $K_j(\gamma)$ be the Bessel function defined by
\begin{equation}\label{defini}	 
K_j(\gamma) :=\frac{(2^j)j!}{(2j)!}\frac{1}{\gamma^j}\int_{\gamma}^{\infty}e^{-\lambda}(\lambda^2-\gamma^2)^{j-1/2}d\lambda,\quad j\geq 0.
\end{equation}
\begin{lemma}\cite{Groot-Leeuwen-Weert-1980,Oliver-1974,Wa}\label{def-pro}
The following relations hold:
\begin{align}
\begin{aligned} \label{transform}
	& (i)~K_j(\gamma)=\frac{2^{j-1}(j-1)!}{(2j-2)!}\frac{1}{\gamma^j}\int_{\gamma}^{\infty}e^{-\lambda}\lambda(\lambda^2-\gamma^2)^{j-3/2}d\lambda,\quad j>0. \\
	& (ii)~K_{j+1}(\gamma)=2j\frac{K_j(\gamma)}{\gamma}+K_{j-1}(\gamma),\quad j\geq1. \\
	& (iii)~K_j(\gamma)<K_{j+1}(\gamma),\quad j\geq 0.  \\
	& (iv)~\frac{d}{d\gamma}\left(\frac{K_j(\gamma)}{\gamma^j}\right)=-\frac{K_{j+1}(\gamma)}{\gamma^j},
	\quad j \geq 0. \\
	&(v)~ K_{j}(\gamma)=\sqrt{\frac{\pi}{2\gamma}}e^{-\gamma}\left(\gamma_{j,n}(\gamma)\gamma^{-n}
	+\sum_{m=0}^{n-1}A_{j,m}\gamma^{-m}\right),\quad j\geq0,~n\geq1,
\end{aligned}
\end{align}
where the coefficients $A_{j,m}$ and $\gamma_{j,n}$  in (v) satisfy
\begin{align*} 
\begin{aligned}
	&A_{j,0}=1, \quad A_{j,m}=\frac{(4j^2-1)(4j^2-3^2)\cdots(4j^2-(2m-1)^2)}{8^m \cdot m!},\quad j\geq0, \quad m\geq1, \\
	&|\gamma_{j,n}(\gamma)|\leq2e^{[j^2-1/4]\gamma^{-1}}|A_{j,n}|,\quad j\geq0,~n\geq1.
\end{aligned}
\end{align*}
\end{lemma}	
In the following proposition, we recall point-wise estimates for the ratio $\frac{K_0(\gamma)}{K_1(\gamma)}$ given in \cite{Ruggeri-Xiao-Zhao-ARMA-2021}. 
\begin{proposition}\cite{Ruggeri-Xiao-Zhao-ARMA-2021}\label{K01p}  
The following assertions hold. 
\begin{enumerate}
\item
If $\gamma\in(\sqrt{2}, \infty)$,  $\frac{K_0(\gamma)}{K_1(\gamma)}$ satisfies
\[
1-\frac{1}{2\gamma}\leq\frac{K_0(\gamma)}{K_1(\gamma)}\leq 1-\frac{1}{2\gamma}+\frac{3}{8\gamma^2}+\frac{3}{16\gamma^3}.
\]
\item	
If $\gamma\in(2, \infty)$, $\frac{K_0(\gamma)}{K_1(\gamma)}$ satisfies
\[ 1-\frac{1}{2\gamma}+\frac{3}{8\gamma^2}-\frac{3}{8\gamma^3}+\frac{63}{128\gamma^4}-\frac{31}{20\gamma^5} \leq \frac{K_0(\gamma)}{K_1(\gamma)}\leq 1-\frac{1}{2\gamma}+\frac{3}{8\gamma^2}-\frac{3}{8\gamma^3}+\frac{63}{128\gamma^4}+\frac{7}{8\gamma^5}. \]
\end{enumerate}
\end{proposition}

\vspace{0.5cm}

\section{The RTCS Model for monatomic gases}\label{App-B}
\setcounter{equation}{0}
In this appendix, we derive the RTCS model with Synge energy for monatomic gases. Now, we begin with system \eqref{A-2} introduced in \cite{Ha-Kim-Ruggeri-ARMA-2020}. 
For the relativistic Euler equations related to the relativistic Boltzmann equations, as in \cite{Groot-Leeuwen-Weert-1980}, the fluid variables $p_a$, and $e_a$ satisfy the relations:
\begin{equation}  \label{ener-den}
p_a =n_aT_a=\rho_a T_a,\quad e_a=c^2\rho_a\frac{K_1(\gamma_a)}{K_2(\gamma_a)}+3p_a, \quad a \in [N].
\end{equation}
In what follows,  $\hbar$ and $k_B$ are the normalized Planck constant and the Boltzmann constant, respectively, and they are assumed to be unity for convenience. \newline

\noindent On the other hand, it follows from \eqref{suminterres} and \eqref{ener-den}  that 
\begin{equation}\label{p0}
p_a=\frac{\rho_a c^2}{\gamma_a}, \quad e_a= \rho_ac^2\left(\frac{K_1(\gamma_a)}{K_2(\gamma_a)}+\frac{3}{\gamma_a}\right), \quad \varepsilon_a=c^2\left(\frac{K_1(\gamma_a)}{K_2(\gamma_a)}+\frac{3}{\gamma_a}-1\right), \quad a \in [N].
\end{equation}
Now we substitute \eqref{p0} into \eqref{A-2} and use $\rho_a\Gamma_a=1$ to get
\begin{equation}\label{mona-rtcs}
\begin{cases}
\displaystyle \frac{\mathrm{d}\mathbf{x}_a}{\mathrm{d} t}=\mathbf{v}_a,\quad a \in [N],\\
\displaystyle \frac{\mathrm{d}}{\mathrm{d} t} \left[ \Gamma_a\mathbf{v}_a\left(\frac{K_1(\gamma_a)}{K_2(\gamma_a)}+\frac{4}{\gamma_a}\right)\right ]
=\frac{1}{N}\sum_{b=1}^{N}\phi_{ab}\Big(\frac{\Gamma_b \mathbf{v}_{b}}{T_b}-\frac{\Gamma_a \mathbf{v}_{a}}{T_a}\Big),\\
\displaystyle \frac{\mathrm{d}}{\mathrm{d} t}\left[ c^2\left(\Gamma_a\left(\frac{K_1(\gamma_a)}{K_2(\gamma_a)}+\frac{4}{\gamma_a}\right)
-1-\frac{1}{\gamma_a \Gamma_a}\right)\right]=\frac{1}{N}\sum_{b=1}^{N}\zeta_{ab}\Big(\frac{\Gamma_a }{T_a}-\frac{\Gamma_b }{T_b}\Big).
\end{cases}
\end{equation}
Note that 
\begin{align}
\begin{aligned} \label{K0/K1}
&\lim_{c\rightarrow \infty} \Gamma_a=1,\quad \lim_{c\rightarrow \infty} c^2(\Gamma_a-1)=\frac{|\mathbf{v}_a|^2}{2}, \quad \mbox{and}  \\
&\frac{K_0(\gamma_a)}{K_1(\gamma_a)}=1-\frac{1}{2\gamma_a}+\frac{3}{8\gamma_a^2}
+\mathcal{O}(\gamma_a^{-3}), 
\quad \frac{K_1(\gamma_a)}{K_2(\gamma_a)}=1-\frac{3}{2\gamma_a}+\frac{15}{8\gamma_a^2}+\mathcal{O}(\gamma^{-3}).
\end{aligned}
\end{align}
Then, it follows from \eqref{transform} and Proposition \ref{K01p} that 
\begin{align}
\begin{aligned} \label{energy-m(D=3)}
&\lim_{c\rightarrow \infty} \Gamma_a=1,\qquad \lim_{c\rightarrow \infty} c^2(\Gamma_a-1)=\frac{|\mathbf{v}_a|^2}{2}, \\
&\lim_{c\rightarrow \infty}c^2\left(\Gamma_a\left(\frac{K_1(\gamma_a)}{K_2(\gamma_a)}+\frac{4}{\gamma_a}\right)
-1-\frac{1}{\gamma_a \Gamma_a}\right)\\
&\hspace{1cm} =\lim_{c\rightarrow \infty} c^2(\Gamma_a-1)+\lim_{c\rightarrow \infty} c^2\Gamma_a\left(\frac{K_1(\gamma_a)}{K_2(\gamma_a)}+\frac{4}{\gamma_a}-1-\frac{1}{\gamma_a \Gamma_a^2}
\right) \\
&\hspace{1cm}  =\frac{|\mathbf{v}_a|^2}{2}+\lim_{c\rightarrow \infty} c^2\Gamma_a\frac{3}{2\gamma_a}=\frac{|\mathbf{v}_a|^2}{2}+\frac{3 T_a}{2}.
\end{aligned}
\end{align}
Therefore, system \eqref{mona-rtcs} reduces to the classical TCS model proposed by Ha and Ruggeri \cite{Ha-Ruggeri-ARMA-2017} in the classical limit ($c \to \infty$):
\begin{equation*}\label{CA-3}
\begin{cases}
\displaystyle \frac{\mathrm{d}\mathbf{x}_a}{\mathrm{d} t}=\mathbf{v}_a, \quad a \in [N],\\
\displaystyle \frac{\mathrm{d} \mathbf{v}_a}{\mathrm{d} t}
=\frac{1}{N}\sum_{b=1}^{N}\phi_{ab}\Big(\frac{ \mathbf{v}_{b}}{T_b}-\frac{\mathbf{v}_{a}}{T_a}\Big),\\
\displaystyle \frac{\mathrm{d}}{\mathrm{d} t}\left(\frac{|\mathbf{v}_a|^2}{2}+\frac{3T_a}{2}\right)=\frac{1}{N}\sum_{b=1}^{N}\zeta_{ab}\Big(\frac{1 }{T_a}-\frac{1 }{T_b}\Big).
\end{cases}
\end{equation*}

\vspace{0.5cm}

%
%
\section{The RTCS Model for polyatomic gases}\label{App-C}
\setcounter{equation}{0}
In this appendix, we derive the RTCS model for polyatomic gases such as diatomic, triatomic, and tetratomic gases.
For polyatomic gases, corresponding to \eqref{ener-den} in Appendix \ref{App-B}, the pressure $p_a$, the energy $e_a$, and the internal energy $\varepsilon_a$ are given as follows (see \cite{Pennisi_Ruggeri}):
\begin{equation}   \label{Synge}
\begin{cases}
\displaystyle  p_a = \frac{n_ac^2}{\gamma_a}=\rho_a T_a,\quad a\in [N], \\
\displaystyle  e_a = \frac{n_ac^2}{A(\gamma_a)K_2(\gamma_a)}\int_0^{\infty} \left[ K_3(\gamma^*)-\frac{1}{\gamma^*}K_2(\gamma^*) \right]\mathcal{I}^\frac{D-5}{2} \mathrm{d}\mathcal{I},\\
\displaystyle  \varepsilon_a = \frac{c^2}{ A(\gamma_a)K_2(\gamma_a)}\int_0^{\infty} \left[ K_3(\gamma^*)-\frac{1}{\gamma^*}K_2(\gamma^*) \right]\mathcal{I}^\frac{D-5}{2} \mathrm{d}\mathcal{I}-c^2,
\end{cases} 
\end{equation}
where $D, \gamma^*$ and $  A(\gamma_a)$ are given as follows.
\begin{equation} \label{AC-1}
D>3, \quad  \gamma^* :=\gamma_a+\frac{\mathcal{I}}{T_a},\quad
A(\gamma_a) :=\frac{\gamma_a}{K_2(\gamma_a)}\int_0^{\infty}\frac{K_2(\gamma^*)}{\gamma^*}
\mathcal{I}^\frac{D-5}{2}\mathrm{d}\mathcal{I}.
\end{equation}
It follows from Lemma \ref{def-pro} that 
\begin{equation}
\hspace{-0.5cm}
\begin{cases}\label{some properties of Bessel functions}
\displaystyle  K_2(\gamma_a)=\frac{2}{\gamma_a}K_1(\gamma_a)+K_0(\gamma_a),\quad K_3(\gamma_a)=\frac{4}{\gamma_a}K_2(\gamma_a)+K_1(y), \\ 
\displaystyle  \frac{\mathrm{d}}{\mathrm{d}\gamma_a}\left(\frac{K_1(\gamma_a)}{\gamma_a}\right) =\frac{1}{\gamma_a}\frac{\mathrm{d}}{\mathrm{d}\gamma_a}K_1(\gamma_a)-\frac{K_1(\gamma_a)}{\gamma_a^2}
=-\frac{K_2(\gamma_a)}{\gamma_a}, \\
\displaystyle \frac{\mathrm{d}}{\mathrm{d}\gamma_a}K_0(\gamma_a)=-K_1(\gamma_a), \\
\displaystyle \frac{\mathrm{d}}{\mathrm{d}\gamma_a}K_1(\gamma_a)=-K_2(\gamma_a)+\frac{K_1(\gamma_a)}{\gamma_a}
=-K_0(\gamma_a)-\frac{1}{\gamma_a}K_1(\gamma_a).
\end{cases}
\end{equation}
For $D>7$, we use \eqref{some properties of Bessel functions} to obtain
\begin{align} 
\label{(1)D=7,9,11}
& (i)~ \int_0^{\infty} \left[ K_3(\gamma^*)-\frac{1}{\gamma^*}K_2(\gamma^*) \right] \mathcal{I}^\frac{D-5}{2}\mathrm{d}\mathcal{I} \nonumber\\
& \hspace{2cm} =\frac{(D-5)(D-7)}{4}T_a^{\frac{D-3}{2}}\int_{\gamma_a}^{\infty}K_1(y)(y-\gamma_a)^\frac{D-9}{2}\mathrm{d}y \nonumber\\
& \hspace{2.2cm} +\frac{D-5}{2}T_a^{\frac{D-3}{2}}\int_{\gamma_a}^{\infty} \frac{2K_1(y)}{y} (y-\gamma_a)^\frac{D-7}{2}\mathrm{d}y. \nonumber\\
& (ii)~\int_0^{\infty}\frac{K_2(\gamma^*)}{\gamma^*}\mathcal{I}^\frac{D-5}{2}\mathrm{d}
\mathcal{I } =T_a^{\frac{D-3}{2}}\int_{\gamma_a}^{\infty}\frac{K_2(y)}{y} (y-\gamma_a)^\frac{D-5}{2}\mathrm{d}y  \\
& \hspace{2cm}  = \frac{D-5}{2}T_a^{\frac{D-3}{2}}\int_{\gamma_a}^{\infty}\frac{K_1(y)}{y} (y-\gamma_a)^\frac{D-7}{2}\mathrm{d}y.\nonumber\\
& (iii)~e_a= n_ac^2\left[\frac{\frac{(D-7)}{2}\int_{\gamma_a}^{\infty}K_1(y)(y-\gamma_a)^\frac{D-9}{2}\mathrm{d}y}
{\gamma_a\int_{\gamma_a}^{\infty}\frac{K_1(y)}{y} (y-\gamma_a)^\frac{D-7}{2}\mathrm{d}y}+\frac{2}{\gamma_a}\right].\nonumber
\end{align}
\subsection{Diatomic gases} \label{App-C-1}
In this part, we consider diatomic gases with $D=5$. In this case,  we use \eqref{some properties of Bessel functions} to obtain
\begin{align*}
\begin{aligned}
& (i)~\int_0^{\infty} \left[ K_3(\gamma^*)-\frac{1}{\gamma^*}K_2(\gamma^*) \right] \mathrm{d}\mathcal{I} =\int_0^{\infty} \left[ K_1(\gamma^*)+\frac{3}{\gamma^*}K_2(\gamma^*) \right] \mathrm{d}\mathcal{I} = T_a\left( K_0(\gamma_a)+\frac{3K_1(\gamma_a)}{\gamma_a}\right). \\ 
& (ii)~\int_0^{\infty}\frac{K_2(\gamma^*)}{\gamma^*}\mathrm{d}\mathcal{I}=T_a  \int_{\gamma_a}^{\infty}\frac{K_2(y)}{y}\mathrm{d}y=T_a\frac{K_1(\gamma_a)}{\gamma_a}.
\end{aligned}  
\end{align*}
Note that $\gamma^*$ is a function of ${\mathcal I}$ (see \eqref{AC-1}). Similar to \eqref{p0}, we obtain
\begin{align} 
\begin{aligned} \label{pei-5d}
p_a &= \frac{n_ac^2}{\gamma_a}=\frac{\rho_a c^2}{\gamma_a}, \quad \varepsilon_a = c^2\left(\frac{K_0(\gamma_a)}{K_1(\gamma_a)}
+\frac{3}{\gamma_a}-1\right), \\
e_a &=\frac{nc^2}{A(\gamma_a)K_2(\gamma_a)}\int_0^{\infty} \left[ K_3(\gamma^*)-\frac{1}{\gamma^*}K_2(\gamma^*) \right] \mathrm{d}\mathcal{I} =\rho_a c^2\left(\frac{K_0(\gamma_a)}{K_1(\gamma_a)}
+\frac{3}{\gamma_a}\right).
\end{aligned}
\end{align}
As in the derivation of \eqref{mona-rtcs}, we can obtain the RTCS model for diatomic gases:
\begin{equation*}
\begin{cases}  \label{RTCSmodel(D=5)}
\displaystyle \frac{\mathrm{d}\mathbf{x}_a}{\mathrm{d} t}=\mathbf{v}_a, \quad   a \in [N], \\
\displaystyle \frac{\mathrm{d}}{\mathrm{d} t}\left[ \Gamma_a\mathbf{v}_a\left(\frac{K_0(\gamma_a)}{K_1(\gamma_a)}+\frac{4}{\gamma_a}\right)\right ]
=\frac{1}{N}\sum_{b=1}^{N}\phi_{ab}\Big(\frac{\Gamma_b \mathbf{v}_{b}}{T_b}-\frac{\Gamma_a \mathbf{v}_{a}}{T_a}\Big),\\
\displaystyle \frac{\mathrm{d}}{\mathrm{d} t}\left[c^2\left(\Gamma_a\left(\frac{K_0(\gamma_a)}{K_1(\gamma_a)}+\frac{4}{\gamma_a}\right)
-1-\frac{1}{\gamma_a \Gamma_a}\right)\right]=\frac{1}{N}\sum_{b=1}^{N}\zeta_{ab}\Big(\frac{\Gamma_a }{T_a}-\frac{\Gamma_b }{T_b}\Big).
\end{cases}
\end{equation*}
It follows from  \eqref{K0/K1} that 
\begin{align}\label{energy-m(D=5)}
c^2\left(\Gamma_a\left(\frac{K_0(\gamma_a)}{K_1(\gamma_a)}
+\frac{4}{\gamma_a}\right)
-1-\frac{1}{\gamma_a \Gamma_a}\right)=\frac{|\mathbf{v}_a|^2}{2}+\frac{5 T_a}{2}+\mathcal{O}(c^{-2}).
\end{align}

\subsection{Triatomic gases} \label{App-C-2}
Consider triatomic gases with $D=7$. From  (v) in \eqref{transform} with $j=1$, we have
\begin{align}\label{int of K1(y)/y}
\int_{\gamma_a}^{\infty}\frac{K_1(y)}{y}\mathrm{d}y
=\sqrt{\frac{\pi}{2}}e^{-\gamma_a}\gamma_a^{-\frac{3}{2}}\left[1-\frac{9}{8\gamma_a}
+\frac{345}{128\gamma_a^2}+\mathcal{O}(\gamma_a^{-3}) \right].
\end{align}
This implies
\begin{align}\label{H_a(D=7)}
&\gamma_a\int_{\gamma_a}^{\infty}\frac{K_1(y)}{y}\mathrm{d}y=\sqrt{\frac{\pi}{2}}e^{-\gamma_a}
\gamma_a^{-\frac{1}{2}}\left[1-\frac{9}{8\gamma_a}+\frac{345}{128\gamma_a^2}+\mathcal{O}(\gamma_a^{-3}) \right],\nonumber
\\ & \frac{K_1(\gamma_a)}{\gamma_a\int_{\gamma_a}^{\infty}\frac{K_1(y)}{y}\mathrm{d}y}= 1+\frac{3}{2\gamma_a}-\frac{9}{8\gamma_a^2}+\mathcal{O}(\gamma_a^{-2}).    
\end{align}
On the other hand, we can use \eqref{some properties of Bessel functions} to get
\begin{align*}
&\int_0^{\infty} \left[ K_1(\gamma^*)+\frac{3}{\gamma^*}K_2(\gamma^*) \right] \mathcal{I}\mathrm{d}\mathcal{I}
=T^2K_1(\gamma_a)+2T_a^2\int_{\gamma_a}^{\infty}\frac{K_1(y)}{y}\mathrm{d}y,
\\&\int_0^{\infty}\frac{K_2(\gamma^*)}{\gamma^*}\mathcal{I}\mathrm{d}\mathcal{I}
=T_a^2\int_{\gamma_a}^{\infty}\frac{K_2(y)}{y}(y-\gamma_a)\mathrm{d}y =T_a^2\int_{\gamma_a}^{\infty}\frac{K_1(y)}{y}\mathrm{d}y.
\end{align*}
From \eqref{Synge} and the above computation for integrals, we can obtain explicit expressions of the pressure $p_a$, the energy $e_a$, and the internal energy $\varepsilon_a$ for triatomic gases:
\begin{equation} \label{pei-7d}
p_a= \frac{\rho_a c^2}{\gamma}, \quad 
e_a=\rho_a c^2\left(\frac{K_1(\gamma_a)}{\gamma_a\int_{\gamma_a}^{\infty}\frac{K_1(y)}{y}\mathrm{d}y}
+\frac{2}{\gamma_a}\right), \quad \varepsilon_a= c^2\left(\frac{K_1(\gamma_a)}{\gamma_a\int_{\gamma_a}^{\infty}
\frac{K_1(y)}{y}\mathrm{d}y}+\frac{2}{\gamma_a}-1\right).
\end{equation}
As in the derivation of \eqref{mona-rtcs}, we can derive the RTCS model related to triatomic gases:
\begin{equation*} 
\begin{cases} \label{RTCSmodel(D=7)}
\displaystyle \frac{\mathrm{d}\mathbf{x}_a}{\mathrm{d} t}=\mathbf{v}_a,~~t > 0, \quad  a \in [N], \\
\displaystyle \frac{\mathrm{d}}{\mathrm{d} t}\left(\Gamma_a\mathbf{v}_a\left(\frac{K_1(\gamma_a)}{\gamma_a\int_{\gamma_a}^{\infty}\frac{K_1(y)}{y}\mathrm{d}y}+\frac{3}{\gamma_a}\right)\right)
=\frac{1}{N}\sum_{b=1}^{N}\phi_{ab}\Big(\frac{\Gamma_b \mathbf{v}_{b}}{T_b}-\frac{\Gamma_a \mathbf{v}_{a}}{T_a}\Big),\\
\displaystyle \frac{\mathrm{d}}{\mathrm{d} t}\left[c^2\left(\Gamma_a\left( \frac{K_1(\gamma_a)}{\gamma_a\int_{\gamma_a}^{\infty}\frac{K_1(y)}{y}\mathrm{d}y}+\frac{3}{\gamma_a}\right)
-1-\frac{1}{\gamma_a \Gamma_a}\right)\right]=\frac{1}{N}\sum_{b=1}^{N}\zeta_{ab}\Big(\frac{\Gamma_a }{T_a}-\frac{\Gamma_b }{T_b}\Big).
\end{cases}
\end{equation*}
Moreover, we can obtain
\begin{equation}\label{energy-m(D=7)}
c^2\left[ \Gamma_a\left(\frac{K_1(\gamma_a)}{\gamma_a\int_{\gamma_a}^{\infty}\frac{K_1(y)}{y}\mathrm{d}y}+\frac{3}{\gamma_a}\right)
-1-\frac{1}{\gamma_a \Gamma_a}\right ]=\frac{|\mathbf{v}_a|^2}{2}+\frac{7 T_a}{2}+\mathcal{O}(c^{-2}).
\end{equation}
\subsection{Tetratomic gases} Consider the relativistic TCS model for tetratomic gases with $D=9$. We use \eqref{(1)D=7,9,11} with $D=9$, to obtain
\begin{align*}
\begin{aligned}
&  (i)~\int_0^{\infty} \left[ K_1(\gamma^*)+\frac{3}{\gamma^*}K_2(\gamma^*) \right] \mathcal{I}^2\mathrm{d}\mathcal{I}
= 2T_a^3K_0(\gamma_a)+4T_a^3\int_{\gamma_a}^{\infty}\frac{K_1(y)}{y}(y-\gamma)\mathrm{d}y.
\\ 
& (ii)~\int_0^{\infty}\frac{K_2(\gamma^*)}{\gamma^*}\mathcal{I}^2\mathrm{d}\mathcal{I}= 2T_a^3K_0(\gamma_a)
-2T_a^3\gamma_a\int_{\gamma_a}^{\infty}\frac{K_1(y)}{y}\mathrm{d}y.
\end{aligned} 
\end{align*}
These and \eqref{Synge} lead to the explicit forms of the pressure $p_a$, the energy $e_a$, and the internal energy $\varepsilon_a$ for tetratomic gases:
\begin{align}
\begin{aligned} \label{pei-9d}
&p_a=\frac{\rho_a c^2}{\gamma},\quad
e_a=\rho_a c^2\left(\frac{K_0(\gamma_a)}{\gamma_aK_0(\gamma_a)-\gamma_a^2
\int_{\gamma_a}^{\infty}\frac{K_1(y)}{y}\mathrm{d}y}+\frac{2}{\gamma_a}\right),
\\&   \varepsilon_a=c^2\left(\frac{K_0(\gamma_a)}{\gamma_aK_0(\gamma_a)-\gamma_a^2
\int_{\gamma_a}^{\infty}\frac{K_1(y)}{y}\mathrm{d}y}+\frac{2}{\gamma_a}-1\right).
\end{aligned}
\end{align}
As in the derivation of \eqref{mona-rtcs}, we get the RTCS model corresponding to the tetratomic gases:
\begin{equation*} 
\begin{cases} \label{RTCSmodel(D=9)}
\displaystyle \frac{\mathrm{d}\mathbf{x}_a}{\mathrm{d} t}=\mathbf{v}_a, \quad  a \in [N], \\
\displaystyle \frac{\mathrm{d}}{\mathrm{d} t}\left[ \Gamma_a\mathbf{v}_a\left(\frac{K_0(\gamma_a)}{\gamma_aK_0(\gamma_a)-\gamma_a^2
\int_{\gamma_a}^{\infty}\frac{K_1(y)}{y}\mathrm{d}y}+\frac{3}{\gamma_a}\right)\right ]
=\frac{1}{N}\sum_{b=1}^{N}\phi_{ab}\Big(\frac{\Gamma_b \mathbf{v}_{b}}{T_b}-\frac{\Gamma_a \mathbf{v}_{a}}{T_a}\Big),\\
\displaystyle \frac{\mathrm{d}}{\mathrm{d} t}c^2\left[\Gamma_a\left(\frac{K_0(\gamma_a)}{\gamma_aK_0(\gamma_a)-\gamma_a^2
\int_{\gamma_a}^{\infty}\frac{K_1(y)}{y}\mathrm{d}y}+\frac{3}{\gamma_a}\right)-\frac{1}{\gamma_a\Gamma_a}-1\right]=\frac{1}{N}\sum_{b=1}^{N}\zeta_{ab}\Big(\frac{\Gamma_a }{T_a}-\frac{\Gamma_b }{T_b}\Big).\nonumber
\end{cases}
\end{equation*}
Note that
\begin{align}\label{H_a(D=9)}
\begin{aligned}
&  \gamma\int_{\gamma_a}^{\infty}\frac{K_1(y)}{y}\mathrm{d}y=\sqrt{\frac{\pi}{2}}
e^{-\gamma_a}\gamma_a^{-\frac{1}{2}}\left[1-\frac{9}{8\gamma_a}
+\frac{345}{128\gamma_a^2}+\mathcal{O}(\gamma_a^{-3}) \right],
\\& \frac{K_0(\gamma_a)}{\gamma_aK_0(\gamma_a)-\gamma_a^2
\int_{\gamma_a}^{\infty}\frac{K_1(y)}{y}\mathrm{d}y}=1+\frac{5}{2\gamma_a}+\mathcal{O}(\gamma_a^{-2})    
\end{aligned}
\end{align}
By (v) in \eqref{transform} with $j=0$  and \eqref{int of K1(y)/y}, we have
\begin{align}\label{energy-m(D=9)}
c^2\left[\Gamma_a\left(\frac{K_0(\gamma_a)}{\gamma_aK_0(\gamma_a)-\gamma_a^2
\int_{\gamma_a}^{\infty}\frac{K_1(y)}{y}\mathrm{d}y}+\frac{3}{\gamma_a}\right)-\frac{1}{\gamma_a\Gamma_a}-1\right]
=\frac{|\mathbf{v}_a|^2}{2}+\frac{9 T_a}{2}+\mathcal{O}(c^{-2}).
\end{align}

	\begin{center}
{\bf Acknowledgment}
\end{center}  The work of Z. M. Bian and Q. H. Xiao was supported by the National Natural Science Foundation of China grant 12271506 and the National Key Research and Development Program of China (No.~2020YFA0714200). 
	The work of S.-Y. Ha was supported by NRF grant funded by the Korea government (MIST) (RS-2025-00514472). 
	The work of T. Ruggeri was carried out in the framework of the activities of the Italian National Group for Mathematical Physics [Gruppo Nazionale per la Fisica Matematica (GNFM/IndAM)].
	
\bigbreak
\noindent\textbf{Data Availability Statement:}
Data sharing is not applicable to this article as no datasets were generated or analysed during the current study.

\noindent\textbf{Conflict of Interest:}
The authors declare that they have no known competing financial interests or personal relationships that could have appeared to influence the work reported in this paper.
	

\end{document}